%
%
%
\let\ReallyTheEnd=\end
\def\end{\vfill\eject}
\ifx\RESOURCE\relax   \fi \let\RESOURCE=\relax



\newcount\Fnt 
\newcount\Ht 
\newcount\Wd 




\def \thefont{}


\def\Accent#1#2{\ifcase \Fnt {#1#2}%
	\or {\rm#1\thefont #2}%
	\or {\bf#1\thefont #2}\fi }

\def\`#1{\Accent{\accent 18 }{#1}}
\def\'#1{\Accent{\accent 19 }{#1}}
\def\^#1{\Accent{\accent 94 }{#1}}
\def\"#1{\Accent{\accent "7F }{#1}}


\def \Bbd#1{{\bf #1}}
\def \Calig#1{{\cal #1}}

\def \Smallfonts {}

\def \Headingfont {\bf }
\def \Subheadingfont {\bf }

\def \Theoremfont {\bf }
\def \Prooffont {\it }
\def \Remarkfont {\bf }
\def \Diagramfont {\bf }



\def \It #1{{\it{#1}\unskip\/}}


\def \Bold #1{{\bf #1}}

\def \Bf { \ifmmode \let \this\Bold \else \let\this\BF \fi \this} 
\def \BF #1{{\bf{\unskip#1}\unskip}}

\def \Admin #1{\begingroup\mathsurround=0 pt
	\leavevmode
	\ifmmode\hbox{$\rm #1$}\else$\rm #1$\fi
	\endgroup }

\def \Rm #1{\hbox{\kern 1pt \rm #1\kern 1pt}}



\def \Subheading#1{\medskip\bigskip \goodbreak \par \noindent {\Subheadingfont #1.} 
\nobreak \vskip 3pt \nobreak }%

\def \Remark#1{\bigskip \goodbreak\par\noindent{\Remarkfont  #1.}}

\def \endRemark{\medskip \goodbreak}

\def \Definition#1{\Remark {#1}}
\def \endDefinition{\medskip \goodbreak}

\def \Theorem #1{\goodbreak\bigskip\par\noindent\Theoremfont #1. \hskip 2pt 
plus 1pt minus 1pt \it}
\def \endTheorem {\rm \goodbreak \smallskip}

\def \Proof#1{\goodbreak \medskip 
\par\noindent \Prooffont #1\hskip .7pt:\hskip 3pt\rm}

\def\proofblock{\hbox{$\sqcap \unskip \kern -6.5pt \sqcup$}}
\def \endProof{\ifmmode \proofblock
		\else \nobreak\hfill\proofblock \hskip 5pt\smallskip \fi}

\def\Endphase {}

\def \Benchmark { }

\def \References#1{\begingroup \leftskip=25 pt \parskip=4 pt plus 2 pt
		\goodbreak \hbox to 1 pt{} \vskip 15 pt plus 10 pt minus 5 pt
\centerline{\Headingfont #1}
\frenchspacing \Smallfonts \def \Benchmark{\Refmark }
\def \Refmark##1##2{\par\noindent \llap {##1{##2}\kern 12 pt}\kern 0pt}
\nobreak\vskip 8pt \nobreak}

\def \endReferences {~\unskip\par\endgroup \medskip\goodbreak }

\def \Phantom{}

\catcode`\w=\active  \catcode`\h=\active

\def\preLinefigure[#1*#2]_{%
\Wd=#1\Ht=#2\catcode`\w=11  \catcode`\h=11 \LLinefigure}

\catcode`\w=11  \catcode`\h=11 

\def \LLinefigure#1{\mathsurround=0 pt%
\setbox1=\hbox{#1}%
\Phantom
{\hskip 0 pt\hbox{$%
	\vcenter{\hbox{%
		\vrule \vbox to \Ht pt{%
			\hrule \vfil \hbox to \Wd pt{%
				\hfil\unhbox1\hss}%
			\vfil\hrule }%
		\vrule }}\hskip 0 pt%
	$}}%
\endgroup}

\def\TexturesLinefigure[#1*#2scaled#3]_#4{\mathsurround=0 pt%
  \dimen1=#1 pc\dimen2=#2 pc%
  \divide\dimen1 by 1000 \multiply\dimen1 by #3%
  \divide\dimen2 by 1000 \multiply\dimen2 by #3%
\noindent\hbox{$%
	\vcenter{\hbox{%
		\vbox to \dimen2{%
			\hbox to \dimen1{%
				\special{picture #4 scaled #3}\hfil}%
			\vfil }%
		}}%
	$}}

\def \metaDiagram#1#2@{
	\def\SetHt##1{\def\test{##1}\def\Test{h} \ifx \test\Test \Ht=40 
		\else \Ht=##1 \fi}
	\SetHt{#1}\goodbreak\midinsert\vskip -8pt 
	\vbox to \Ht pt{\vfil \noindent\hfil\Diagramfont#2 \hfil }
\vskip-8pt\endinsert}

\def \Diagram#1{\metaDiagram#1@} 


\def \TexturesMetaDiagram#1#2@{\goodbreak\midinsert
	\vskip\abovedisplayskip
		\line{\hfil #1\hfil}%
		\vskip\bigskipamount
		\line{\hfil\Diagramfont #2\hfil }%
		\vskip\belowdisplayskip
\endinsert}


\def \Item #1{\item {$\bf {{#1}} $}}

\def \Cite #1{{\bf [#1]}}


\def \preXbox{\hbox{$
	\vcenter{\hbox{
		\vrule\vbox to 6.7 pt{
			\hrule \vfil \hbox to 12 pt{
				\hfil}%
			\vfil\hrule}%
		\vrule}}\hskip 4pt%
	$}}

\def \Xbox{\raise -.25pt\hbox{\preXbox}}

\def \Nonsense {{~\unskip \kern-3.5 pt \mathsurround=0 pt%
		\hbox{\Xbox \kern -16.5 pt $>\kern-3pt<$}}}

\def \Lguillemets {{$\mathsurround= 0 pt\raise 1.4 pt\hbox{$\scriptscriptstyle 
	\langle \kern -1 pt \langle\hskip 2 pt $}$}}
 
\def \Rguillemets {{\mathsurround= 0 pt$\hskip 2 pt \raise 1.4 pt\hbox{$\scriptscriptstyle 
	\rangle \kern -1 pt \rangle$}$}} 

\def \Eqno #1$${\eqno \Admin{#1}$$}

\def \Rparen {\right ) }
\def \Lparen {\futurelet\next \Lptaupe}
\def\Lptaupe{\ifx \next ^ \let\this\LLparen 
	\else \let\this\LLLparen  \fi\this} 
\def\LLparen {\left ( \Atop }
\def \Atop ^#1_#2{{#1\atop#2}}
\def\LLLparen {\left (}

\def \bigMidvert{\kern4pt \big \vert \kern4pt}

\def \Midvert{\kern3pt \vert \kern3pt}

\def \Sharp {{\mathsurround=0pt\kern1pt%
\hbox{$\vcenter{\hbox{$\scriptstyle \# $}\vskip.7pt}$}\kern1pt}}

\def \Cup{\bigcup }
\def \Cap{\bigcap }

\def \Otimes{\mathbin{\kern-2pt\raise 1.2pt%
	\hbox{$\scriptstyle \otimes$}\kern-2pt}}

\def \Oplus{\mathbin{\kern-2pt\raise 1.2pt%
	\hbox{$\scriptstyle \oplus$}\kern-2pt}}

\def \Amalg{\mathbin{\raise .5pt%
	\hbox{$\scriptstyle \amalg$}}}



\def \Circ {\circ}

\def \Coprod {\mathop{\raise 1.2pt \hbox{$\coprod$}}}


\def \Acc{\expandafter }

\def\swthat{\raise -1.1 ex\hbox{$\widehat{}$}}
\def\swttilde{\raise -1.2 ex\hbox{$\widetilde{}$}}
\def \overdot{{\raise .2 ex \hbox to 0pt {\hss\bf\smash{.}\hss}}}
\def \overcircle{{\raise .1 ex \hbox to 0pt
{\mathsurround=0pt$\scriptstyle\hss\circ\hss$}}}

\def \Mathaccent#1#2{\mathsurround=0 pt
\setbox4=\hbox{$\vphantom{#2}$}
\Ht=\ht4 
\setbox5=\hbox{${#1}$}
\setbox6=\hbox{${#2}$}
\setbox7=\hbox to .5\wd6{}
\copy7\kern .1\Ht \raise\Ht sp\hbox{\copy5}\kern-.1\Ht 
\copy7\llap{\box6}
}

\def \Overcircle #1{\Mathaccent {\overcircle} {#1}}

\def \SwtHat #1{\Mathaccent {\swthat} {#1}}

\def \SwtTilde #1{\Mathaccent {\swttilde} {#1}}

\def \ChOline#1{\mathsurround=0pt\setbox1=\hbox{${#1}$}
	\ifdim \wd1 > 7pt
		\kern .15\ht1 \kern .9 pt 
		\overline {\kern -.15\ht1 \kern -.9 pt#1\kern-.9 pt}
		\kern .9 pt 
	\else
		\ifdim \wd1 > 4pt
			\kern .3\ht1 
			\overline {\kern -.3\ht1 {#1}}
		\else
			\kern .3\ht1 \kern-.9 pt 
			\overline {\kern -.3\ht1 \kern .9 pt{#1}\kern .9 pt }
			\kern-.9 pt
		\fi
	\fi}

\def \ChUline#1{
	{\kern .5pt \underline {\kern -.5pt#1\kern-2.2pt}\kern1.3pt}
	}

\def \Uuline#1{\closerunderline{.9pt}{\closerunderline {-.3pt}{#1}}}
\def \ChUuline#1{
	{\kern .5pt \Uuline {\kern -.5pt#1\kern-2.2pt}\kern1.3pt}
	}

\def \Cdot{\mathbin{\raise .4 ex \hbox to 0pt {\hss\bf .\hss}}}


\def \Longrightarrow {\kern-2pt\mathop{\kern3pt\longrightarrow\kern3pt}\limits}

\def \Longleftarrow{\kern-2pt\mathop{\kern3pt\longleftarrow\kern3pt}\limits}
 
\def \Longtwoheadrarrow {\kern-2pt
	\mathop{\kern3pt\longrightarrow \kern-14pt \longrightarrow\kern3pt}\limits}

\def \Longleftrightarrow {\kern-2pt
	\mathop{\kern3pt\longleftrightarrow\kern3pt}\limits}

\def \Longmapsto {\kern-2pt
	\mathop{\kern3pt\longmapsto\kern3pt}\limits}

\def\rarrow{\rightarrow }

\def\LaTeX{{\rm L\kern-.34em\raise.47ex\hbox{\mathsurround=0pt 
$\scriptstyle\rm A$}\kern-.15em \TeX}} 


\def\Matrix#1{\matrix{#1}}




\hsize=350 pt
\vsize=520 pt
\tolerance=1000
\widowpenalty=5000
\mathsurround=1.2pt%

\abovedisplayskip=5pt plus3pt minus2pt
\belowdisplayskip=5pt plus3pt minus2pt
\predisplaypenalty=1000 %

\lineskip=1.5pt
\lineskiplimit=.8pt
\topskip=0pt
\parskip=2pt plus 1pt
\bigskipamount=10pt plus 4pt minus 1pt
\hfuzz=2\hsize 
\hoffset=10pt

\quad
\vskip150pt
\centerline{\bf A PARTIAL DESCRIPTION OF THE PARAMETER SPACE}
\centerline{\bf OF RATIONAL MAPS OF DEGREE TWO: PART 2}
\vskip50pt
\centerline{\bf Mary Rees}
\vskip30pt
\centerline{ Department of Pure Mathematics, University of Liverpool,}
\centerline{ P.O. Box 147, Liverpool L69 3BX, U. K.}
\vskip30pt
\centerline{ Address until 30th June 1991: }
\centerline{ Institute of Mathematical Sciences, SUNY at Stony Brook,}
\centerline{ Stony Brook NY 11794-3660, U.S.A.}
\vskip30pt
\centerline{\bf Abstract}
This continues the investigation of a combinatorial model for the variation of
dynamics in the family of rational maps of degree two, by concentrating on those varieties in which one critical point is periodic. We prove some general
results about nonrational critically finite degree two branched coverings,
and finally identify the boundary of the rational maps in the combinatorial 
model, thus completing the proofs of results announced in Part 1.
\end

\Subheading {Introduction}  

This is the second in a series of papers aimed at 
understanding the variation in dynamical 
behaviour within the family of rational maps of 
degree two. This introduction will be brief, as 
Chapter 1 of \Cite{R2} was intended as an 
introduction to this paper also, and includes 
statements of the results of this paper.

As in \Cite{R2}, we concentrate on maps which are 
\It{critically finite. } Recall that an 
orientation-preserving branched covering 
$g:\Acc\ChOline \Bbd C\rarrow \Acc\ChOline \Bbd 
C$ is called \It{critically finite } if 
${\Sharp}(X(g))<\infty $, where
$$X(g)=\lbrace g^{n}(c):c\Rm { critical, 
}n>0\rbrace ,$$
and that these play an important r\^ole in the 
understanding of rational maps. A branched 
covering of $\Acc\ChOline \Bbd C$ of degree two 
necessarily has two critical points $c_{1}$, 
$c_{2}$. More generally, we consider degree two 
branched coverings for which one critical point 
is periodic. Likewise, these play an important 
r\^ole. Douady and Hubbard have studied the 
\It{Mandelbrot set } for degree two polynomials 
(for which the critical point at $\infty $ is, of 
course, fixed). They provided a combinatorial 
model for the Mandelbrot set (\Cite{D-H1}, 
\Cite{D-H2}) which was reinterpreted by Thurston 
\Cite{T}. The work described in this paper, and 
its predecessor and successors, is a 
generalisation of the combinatorial model of the 
Mandelbrot set to maps with a critical point 
which is no longer fixed, but still periodic. It 
seems likely that this model can be completed to 
include maps for which both critical points have 
infinite forward orbits.

Recall that two critically finite branched 
coverings $f$, $g$ are \It{equivalent }, written 
$f\simeq g$, if there exists an 
orientation-preserving homeomorphism $\varphi 
:\Acc\ChOline \Bbd C\rarrow \Acc\ChOline \Bbd C$ 
and a path $g_{t}$ ($t\in [0,1]$) through 
critically finite branched coverings such that 
$X(g_{t})=X(g)$ for all $t\in [0,1]$, $\varphi 
\Circ f\Circ \varphi ^{-1}=g_{0}$, $g=g_{1}$. 
Thurston \Cite{T} gave a necessary and sufficient 
condition for a critically finite branched 
covering to be equivalent to a rational map. 
Almost always - and always in degree two - there 
is at most one rational map in an equivalence 
class, and the rational map can be determined up 
to topological conjugacy (and its Julia set, up 
to homeomorphism) by any map in the equivalence 
class. More detail is given on this in \Cite{R2}. 

We are interested in understanding rational maps, 
but, paradoxically, the subject of this paper is, 
very largely, critically finite branched 
coverings which are not equivalent to rational 
maps. The reason is basically the following. 
Roughly speaking, it is easier to describe 
dynamically the family of all critically finite 
degree two branched coverings with, say, one 
critical point of period $m$ than simply the 
family of all such rational maps. (We remark that 
there is not so much difference between the two 
families if $m=1$.) It is also easier to analyse 
nonrational branched coverings (up to 
equivalence) rather than rational maps, and then, 
of course, the rational ones are what remain when 
the nonrational ones are removed. 

The chapters of this paper give increasingly 
precise results about increasingly specified 
branched coverings. The last two chapters concern 
lamination maps, which are, for the present, our 
preferred models for rational maps. Their 
dynamics are easily computable, and a lamination 
map determines an equivalent rational map up to 
topological conjugacy. This is the substance of 
the Lamination Map Conjugacy Theorem in 
\Cite{R2}.
 The chapters, and their rough content, are as 
follows.

\par \noindent \It{Chapter 1: Homotopy 
Information about Branched Coverings of Degree 
Two. }

Thurston's theorem giving a necessary and 
sufficient condition for equivalence of a 
critically finite branched covering $f$ to a rational map (\Cite{T}, 
\Cite{D-H3}) is in terms of the action of 
$f^{-1}$ on isotopy classes of loops in 
$\Acc\ChOline \Bbd C\setminus X(f)$. This 
condition, and a couple of refinements in degree 
two, which are essentially isotopy-invariance 
conditions, are stated in the introduction of 
\Cite{R2}. This chapter examines more closely the 
structure of $f$ when $f$ is \It{not } equivalent 
to a rational map. In particular, a \It{critical 
branched covering } is associated to $f$ (1.9), 
essentially a restriction of an iterate of $f$ to 
a subsurface. It may be degree one or two, and is 
in some sense irreducible. The main result of the 
chapter is the Decomposition Proposition of 1.12.

\par \noindent \It{Chapter 2: The Fixed 
Holed-Sphere Theorem. }

The isotopy-invariance of the first chapter is 
translated into actual invariance, in the Fixed 
Holed-Sphere Theorem. Underlying this is an 
application of the common idea that a map with 
homotopy type satisfying certain ``expanding'' or 
``hyperbolic'' properties is determined up to a 
useful semiconjugacy by its homotopy type.

\par \noindent \It{Chapter 3:  
}$1$\It{-dimensional Topological Dynamics }.

No particular originality is claimed for this. 
The $1$-sided Neighbourhood Theorem is similar 
(though not identical, so far as I know) to many 
in the literature, and shows, roughly, that under 
certain mild conditions, a graph preserved by a 
map is a union of unstable manifolds.

\par \noindent \It{Chapter 4: The Nonrational 
Lamination Map Theorem. }

Now we specify to lamination maps, and show that 
an invariant lamination with a critically finite 
lamination map with a nondegenerate Levy cycle 
(roughly speaking, not equivalent to a rational 
map) contains periodic leaves - with a specified 
bound on the period - with certain additional 
properties. Finally, this result is used in

\par \noindent \It{Chapter 5: Proof of the 
Admissible Boundary Theorem. }

Lamination maps (at least, those which have been 
defined at any point in this work) have $0$ as a 
periodic critical point, and lie in countably 
many families, depending on the definition of the 
map on a set containing the orbit of $0$. The 
lamination maps in each family are identified 
with leaves or complementary components - called 
gaps - of a lamination on a disc, called a 
parameter lamination and denoted by $\Calig 
L_{r}$, where $r$ runs through the odd 
denominator rationals in $(0,1)$. For details of 
this structure, see Chapter 1 of \Cite{R2}. This 
chapter contains the proof of two results stated 
in the introduction of \Cite{R2} - the 
Admissibility Proposition and the Admissible 
Boundary Theorem. Identifying lamination maps 
with the corresponding leaves and gaps in the 
laminated disc, these two theorems together 
identify an open set of lamination maps (and its 
boundary). In this set, all critically finite 
branched coverings have at most only degenerate 
Levy cycles, which virtually means they are all 
equivalent to rational maps, and the set contains 
all critically finite rational maps, up to 
equivalence.
\end

\pageno=4
\Subheading {Chapter 1. Homotopy information 
about Branched Coverings of Degree Two}  

\Subheading {1.1}  

The information in this chapter is not very deep, 
and most of it is known anyway: it is repeated 
for convenience. Let $f:\Acc\ChOline \Bbd 
C\rarrow \Acc\ChOline \Bbd C$ be a critically 
finite degree two orientation-preserving branched 
covering. Then exactly one of the following 
occurs, as can be extracted from \Cite{T}, and is 
explained in 1.2 - 1.3 below.

\par \noindent  Up to equivalence, $f$ is a 
rational map.

\par \noindent  Up to equivalence, $f$ is 
$[x]\mapsto [Ax+b]:T^{2}/\sim \rarrow T^{2}/\sim 
$, where $T^{2}$ is the two-torus, $\sim $ is the 
relation $x\sim -x$, $[x]$ is the equivalence 
class of $x$, $b\sim -b$ and $A$ is an 
automorphism of $T^{2}$ with determinant 2 and 
eigenvalues $\lambda _{s}$, $\lambda _{u}$, 
$\Midvert \lambda _{s}\Midvert <1<2<\Midvert 
\lambda _{u}\Midvert $.

\par \noindent  There is a Levy cycle for $f$. 
(This term is recalled below.)

The third possibility is in some sense a 
decomposable case. We shall define a \It{complete 
Levy set } for $f$, consisting of a union of 
disjoint Levy cycles and some preimages of these. 
Given $f$ and a complete Levy set $\Gamma $, we 
shall see that components of $\Acc\ChOline \Bbd 
C\setminus \cup \Gamma $ are eventually periodic 
in a natural sense, and hence have maps 
associated to them, which are always either 
homeomorphisms or degree two critically finite 
branched coverings. In particular, such a map is 
associated to the component of $\Acc\ChOline \Bbd 
C\setminus \cup \Gamma $ containing both critical 
points, in the case when the forward orbits of 
both critical points of $f$ contain periodic 
critical points, and this map is of one of the 
first two indecomposable types if it is degree 
two. We shall summarise most of the information 
in the Decomposition Proposition (1.12). This 
information will be crucial in the semiconjugacy 
results for critically finite degree two branched 
coverings, in Chapter 2. It is clear that a more 
complete description of those critically finite 
degree two branched coverings with Levy cycles is 
possible, partly in terms of a natural tree 
obtained by treating Levy loops as vertices, and 
partly in terms of the maps associated with 
components of $\Acc\ChOline \Bbd C\setminus \cup 
\Gamma $. I suppose the description would be 
reminiscent of the Douady-Hubbard-Thurston 
description of critically finite quadratic 
polynomials, with some additional 
identifications. In fact, Shishikura \Cite{S} has 
some results in a more general framework: he has 
constructed a tree (with an expanding metric) 
using Thurston obstruction loops as vertices.

\Subheading {1.2. The cases when a Thurston 
obstruction does not exist}  

The Thurston obstruction is a necessary and 
sufficient condition  for nonequivalence to a 
rational map, when a certain orbifold $Q_{g}$, 
associated to an orientation-preserving 
critically finite branched covering $g$, is 
\It{hyperbolic } (\Cite{T} Ch. 16). Without going 
into too much detail about what this means, we 
identify those $g$ of degree two which are not 
equivalent to rational maps, and for which 
$Q_{g}$ is not hyperbolic. If $Q_{g}$ is not 
hyperbolic, then ${\Sharp}(X(g))\leq 4$ (\Cite{T} 
Ch. 16).

If $g$ and $h$ are critically finite 
orientation-preserving branched coverings and 
${\Sharp}(X(g))\leq 3$, then $g$ and $h$ are 
equivalent if and only if there is an isomorphism 
between $(g^{-1}X(g),g)$ and $(h^{-1}X(h),h)$ 
which maps critical points to critical points of 
the same multiplicity. This follows from the fact 
that any homeomorphism of the sphere which fixes 
three points is isotopic to the identity, via an 
isotopy fixing those three points. From this, it 
is easy to see that if ${\Sharp}(X(g))\leq 3$, 
$g$ is equivalent to a rational map.

If ${\Sharp}(X(g))=4$ and $Q_{g}$ is not a 
hyperbolic orbifold, then there is only one 
possibility for the universal orbifold cover 
$\Acc\SwtTilde  Q_{g}$ of $Q_{g}$ and for the 
covering map $\pi :\Acc\SwtTilde  Q_{g}\rarrow 
Q_{g}$: $\Acc\SwtTilde  Q_{g}=\Bbd C$, $Q_{g}$ is 
the quotient of $\Bbd C$ by the relation $z\sim 
\pm z+m+ni$ ($m$, $n\in \Bbd Z$) and $\pi $ is 
the quotient map. Since $g^{-1}$ lifts to $Q_{g}$ 
and induces a homomorphism of a subgroup of the 
covering group onto the covering group (as 
always: see \Cite{T} 16.6) and the covering group 
in this case is of the form $\Bbd Z^{2}\times 
_{s}\Bbd Z_{2}$, we see that $g^{-1}$ lifts 
injectively in this case, so that $g$ also lifts 
- to $\Acc\SwtTilde  g:\Bbd C\rarrow \Bbd C$, 
preserving equivalence classes in general, and 
those of critical orbits in particular. Hence, 
$g$ is equivalent to $[z]\mapsto [Az+\beta ]$, 
for some real-linear map $A$ and $\beta \in (\Bbd 
Z/2)^{2}$. In addition, $\Rm {det}(A)=\Rm { 
degree}(g)$. Let
$$A=\Lparen \Matrix {a&b \Endphase \cr  
                              c&d \Endphase \cr  
}\Rparen $$
 If $(a+d)^{2}<4(ad-bc)$, $A$ is conjugate to a 
skew-symmetric matrix, and $g$ is equivalent to a 
rational map. Now let $ad-bc=2=\Rm { degree}(g)$. 
If $(a+d)^{2}=9$, then the existence of an 
eigenvector of 
$A$ with eigenvalue $\pm 1$ implies the existence 
of a simple loop $\gamma $ in $\Acc\ChOline \Bbd 
C\setminus X(g)$ such that a component $\gamma '$ 
of $g^{-1}(\gamma )$ is isotopic to $\gamma $ in 
$\Acc\ChOline \Bbd C\setminus X(g)$, and 
$g\Midvert \gamma '$is a homeomorphism. If 
$(a+d)^{2}>9$, then
$A$ has eigenvalues $\lambda _{s}$, $\lambda 
_{u}$ with $\Midvert \lambda _{s}\Midvert 
<1<2<\Midvert \lambda _{u}\Midvert $.
\Subheading {1.3. Levy Cycles}  
Let $g:\Acc\ChOline \Bbd C\rarrow \Acc\ChOline 
\Bbd C$ be a critically finite degree two 
orientation-preserving branched covering.  Then 
we recall that a set $\lbrace \gamma 
_{1}...\gamma _{r}\rbrace $ of simple nontrivial 
disjoint loops in 
$\Acc\ChOline \Bbd C\setminus X(g)$ is a \It{Levy 
cycle for } $g$ if there exists $\gamma _{i}'$ in 
$g^{-1}(\gamma _{i+1})$ (writing $\gamma 
_{r+1}=\gamma _{1}$) such that $\gamma _{i}$ and 
$\gamma _{i}'$are isotopic in $\Acc\ChOline \Bbd 
C\setminus X(g)$ and $g\Midvert \gamma _{i}'$ is 
a homeomorphism. We recall that Thurston's 
theorem for a necessary and sufficient condition 
for \It{non }-equivalence to a rational map 
\Cite{T} was given three equivalent formulations 
in the case of degree two in 1.6 of \Cite{R2} 
(where further references are given). The middle 
condition was the existence of a Levy cycle. In 
\Cite{R2} the hypothesis of Thurston's theorem 
was slightly simplified. In fact, the theorem 
holds whenever the associated orbifold is 
hyperbolic. Non-hyperbolic  orbifolds have been 
investigated in 1.2, and, putting these facts 
together, we see that, if $g$ is an 
orientation-preserving degree two critically 
finite branched covering, exactly one of the 
three properties stated in 1.1 does indeed hold. 

So now we wish to obtain more information about 
$g$, where $g$ is a degree two critically finite 
orientation-preserving rational map with a Levy 
cycle. We start with two lemmas, which are, in 
fact, the basis of a proof of existence of a Levy 
cycle in the first place (\Cite{TL}, for 
example). The proofs are included for 
convenience, since they are quite short.

\Theorem {1.4. Lemma}  Let $g:\Acc\ChOline \Bbd 
C\rarrow \Acc\ChOline \Bbd C$ be a critically 
finite degree two orientation-preserving branched 
covering. Let $\Gamma $ be a (finite) set of 
simple disjoint isotopically distinct loops in 
$\Acc\ChOline \Bbd C\setminus X(g)$  with $\Gamma 
\subset g^{-1}\Gamma $ up to isotopy. Let $\gamma 
\in \Gamma $ separate critical values. Then 
$\gamma $ bounds a disc $D$ such that all 
components of $g^{-n}D$ are discs, for all $n\geq 
0$ .\endTheorem   

\Proof {Proof}  

Let $\Gamma _{1}$ be the set of loops $\gamma $ 
in $\Gamma $ such that $\gamma $ bounds a disc 
$D$ for which all components of $g^{-n}D$ are 
discs for all $n\geq 0$. Then $\Gamma \cap 
g^{-1}\Gamma \subset \Gamma _{1}$. (We mean ``up 
to isotopy in $\Acc\ChOline \Bbd C\setminus 
X(g)$'', but we omit saying this, for the rest of 
this proof.) Let $\Gamma _{2}=\Gamma \setminus 
\Gamma _{1}$. Then $\Gamma _{2}\subset 
g^{-1}\Gamma _{2}$. We claim that no loop in 
$\Gamma _{2}$ separates critical values. 
Otherwise, let $D$ be a disc containing exactly 
one critical value, bounded by a loop $\gamma $ 
from $\Gamma _{2}$, and containing no loop from 
$\Gamma _{2}$ in its interior. Then we can prove 
by induction on $n$ that all components of 
$g^{-n}D$ are discs, and none contains a loop 
from $\Gamma _{2}$ in its interior. For assume 
this is true for n, and $D_{1}$ is a disc 
component of $g^{-n}D$. Then $D_{1}$ cannot 
contain both critical values, since it does not 
contain $\gamma $ in its interior. So components 
of $g^{-1}D_{1}$ are all discs, and do not 
contain elements of $g^{-1}\Gamma _{2}\supset 
\Gamma _{2}$ in their interiors. The inductive 
step is completed. So $\gamma \in \Gamma _{1}$, 
giving a contradiction, and the lemma is proved.

\Subheading {1.5. Degenerate Levy Cycles}  

A Levy cycle $\Gamma =\lbrace \gamma 
_{1}...\gamma _{r}\rbrace $ for $g$ is 
\It{degenerate  } if each $\gamma _{i}$, $\gamma 
_{i}'$ (as in 1.3) bounds a disc $D_{i}$, 
$D_{i}'$, with $D_{i}'$ isotopic to $D_{i+1}$ via 
an isotopy preserving $X$, and $g:D_{i}'\rarrow 
D_{i}$ is a homeomorphism. Thus no disc $D_{i}$ 
contains a critical value of $g$. So if, for 
instance, the forward orbits of both critical 
points of $g$ contain periodic critical points, 
any Levy cycle for $g$ must be \It{nondegenerate 
}. The following lemma is then immediate.

\Theorem { Lemma}  

Let  $g$ be a critically finite degree two 
orientation-preserving branched covering, such 
that any Levy cycle is nondegenerate.  Let  
$\gamma $ be a loop in  $\Acc\ChOline \Bbd 
C\setminus X(g)$ bounding a disc $D$ such that 
all components of $g^{-n}D$ are discs for all 
$n\geq 0.$ Then  $D$  cannot contain any loop of 
a Levy cycle.\endTheorem   

\Subheading {1.6. Complete Levy sets}  

Let $f:\Acc\ChOline \Bbd C\rarrow \Acc\ChOline 
\Bbd C$ be a critically finite degree two 
orientation-preserving branched covering, and let 
$\Gamma $ be a finite set of disjoint simple 
nontrivial loops in $\Acc\ChOline \Bbd C\setminus 
X(f)$. Then $\Gamma $ is a \It{complete Levy set 
} (for $f$) if the conditions below hold, writing 
$X=X(f)$.

\par \noindent \It{The Levy Condition. } $\Gamma 
\subset f^{-1}\Gamma $, and for $\gamma \in 
f^{-1}\Gamma $, $\gamma \in \Gamma $ if and only 
if $\gamma $ does not bound a disc $D$ such that 
all components of $f^{-n}D$ are discs, for all 
$n$.

Thus, by Lemma 1.4, no loop in $\Gamma $ 
separates critical values, and $f$ is degree one 
restricted to any component of $f^{-1}\Gamma $. 
There are Levy cycles in $\Gamma $, but $\Gamma $ 
may not be a union of Levy cycles.

\par \noindent \It{The Gap Invariance Condition. 
} If $G$ is a component of $\Acc\ChOline \Bbd 
C\setminus \cup \Gamma $ such that $G\setminus X$ 
is not an annulus, then there is a unique 
component $G'$ of $\Acc\ChOline \Bbd C\setminus 
X$ - which we term $f_{*}G$ - and a unique 
component $G''$ of $f^{-1}G'$ such that 
$G''\subset G$ and $\partial G\subset \partial 
G''$.

\Definition { Definition}  $G$ is a \It{gap for } 
$\Gamma $ if $G$ is a component of $\Acc\ChOline 
\Bbd C\setminus \cup \Gamma $ such that 
$G\setminus X$ is not an annulus.\endDefinition   

\par \noindent \It{The Gap Restrictions. } No 
annulus gap of $\Acc\ChOline \Bbd C\setminus 
(X\cup (\cup \Gamma ))$ is mapped onto itself by 
$f^{n}$ (any $n>0$). If $G_{1}$, $G_{2}$ are 
distinct with $\Acc\ChOline G_{1}\cap 
\Acc\ChOline G_{2}\not =\phi  $, then at least 
one $G_{i}$ is not mapped homeomorphically by all 
$f^{n}$. Any Levy cycle disjoint from $\Gamma $ 
can be isotoped  in 
$\Acc\ChOline \Bbd C\setminus X$
either into $\Gamma $ itself or 
into a periodic cycle of gaps on which $f$ acts 
homeomorphically.

\Theorem {1.7. Lemma} Let $g$ have at least one 
nondegenerate Levy cycle, but no degenerate ones. 
Then there exists a complete Levy set $\Gamma $ 
for some $f$, where $f$ is homotopic to $g$ via a 
homotopy $\lbrace f_{t}\rbrace $ through 
critically finite branched coverings with 
$X(f_{t})=X(g)$ for all $t$.\endTheorem  

\Proof {Proof}   Let $\Gamma _{0}$ be a maximal 
union of disjoint Levy cycles for $g$. Let $k>0$ 
be such that, if $\gamma \in g^{-(k+1)}\Gamma 
_{0}$, and $\gamma $ does not bound a disc $D$ 
such that all components of $g^{-n}D$ are discs 
for all $n>0$, then $\gamma \in g^{-k}\Gamma 
_{0}$ up to isotopy. Let $\Gamma _{1}$ be a set 
of two copies, up to isotopy, of all such loops 
in $g^{-(k+1)}\Gamma _{0}$ (or, equivalently, in 
$g^{-k}\Gamma _{0}$). Then $g$, $\Gamma _{1}$ 
satisfy the Levy Condition except that $\Gamma 
_{1}\subset g^{-1}\Gamma _{1}$ up to isotopy 
only. Note that $\Gamma _{1}$ contains up to 
isotopy all loops in $g^{-1}\Gamma _{1}$ which 
separate two elements of $\Gamma _{1}$. Hence 
$G''$, $G'$ exist as in the Gap Condition, up to 
isotopy. Then $f$ can be chosen equivalent to $g$ 
(via a homotopy preserving $X$) so that the Levy 
and Gap Conditions are satisfied by $f$, $\Gamma 
_{1}$. If $G_{1}$, $G_{2}$ are distinct gaps for 
$\Gamma _{1}$ such that $f^{n}$ always maps both 
of them homeomorphically, then if one of them is 
periodic they both are, and the period of both is 
$\leq $ the period of the common boundary 
component. So obtain $\Gamma $ from $\Gamma _{1}$ 
by removing all such common boundary components 
and all their preimages under all $f^{n}$, and 
also removing redundant loops in annuli (possibly 
removing some annuli gaps altogether) so that 
$\Gamma $ satisfies the Gap Restrictions. The 
last Gap Restriction is satisfied because $\Gamma 
_{0}$ was a maximal union of Levy cycles.

\Subheading {1.8. Critical gaps for Complete Levy 
Sets}  

From now on, until the end of this chapter, let 
$f:\Acc\ChOline \Bbd C\rarrow \Acc\ChOline \Bbd 
C$ be a critically finite degree two 
orientation-preserving branched covering with no 
degenerate Levy cycle, and let $\Gamma $ be a 
complete Levy set for $f$. Then by Lemma 1.4, and 
the Levy Condition, both critical values of $f$ 
are contained in the same gap for $\Gamma $. 
Hence, both critical points are contained in the 
same gap, which we call the \It{critical gap }, 
usually denoted by $G_{0}$. The following lemma 
holds, under the hypotheses of this section. It, 
too, can be found essentially in \Cite{TL}.

\Theorem {Lemma}  The disc gaps for $\Gamma $ are 
$f_{*}^{i}G_{0}$ ,with $1\leq i\leq r$ for some 
$r\geq 2$ .\endTheorem   

\Proof {Proof}  There is a gap $G$ which is a 
disc. Then $G\cap X(f)\not =\phi  $. So 
$f^{n}c\in G$ for some least $n\geq 1$ and $c$ 
critical. If $n>1$ then both components of 
$f^{-1}G$ are discs. Then if we write $G_{n}=G$, 
we can find $G_{i}$ ($1\leq i\leq n$) with 
$f_{*}G_{i}=G_{i+1}$ for $i<n$ and $f^{i}c\in 
G_{i}$. Then $G_{i}$ is a disc for each $i$, and 
$G_{1}$ is the gap containing the critical 
values, as required. We must have $r\geq 2$, 
because $\Acc\ChOline \Bbd C$ has Euler 
characteristic 2, so at least two components of 
$\Acc\ChOline \Bbd C\setminus (\cup \Gamma )$ are 
discs.

\Subheading {1.9. Branched coverings associated 
to gaps}  

 Let $G$ be a gap for $\Gamma $ which is of 
period $n$ under $f_{*}$. (All gaps must be 
eventually periodic under $f_{*}$, since there 
are only finitely many gaps.) Let $G_{1}$ be the 
component of $f^{-n}G$ with $G_{1}\subset G$, 
$\partial G\subset \partial G_{1}$. Then 
$G_{1}=SG$ for an inverse branch $S$ of the 
branched covering $f^{n}$. Inductively we can 
define $G_{k+1}=SG_{k}$. Then for some least 
$k\geq 0$, all components of $\partial 
G_{k+1}\setminus \partial G_{k}$ are trivial in 
$\Acc\ChOline \Bbd C\setminus X(f)$.

Now $f^{n}$ maps components of $\partial G_{k+1}$ 
to components of $\partial G_{k}$, and maps 
components of $\partial G$ to components of 
$\partial G$. Also, $f^{n}\Midvert G_{k+1}$ is 
either a homeomorphism or a degree two branched 
covering, depending on whether or not the 
critical gap is in the forward orbit of $G$. Then 
we define $h=h_{G}$ by

\par \noindent a) $h=f^{n}$ on $G_{k+1}$,
\par \noindent b) $h(D)=D'$ if $D$, $D'$ are 
components of $\Acc\ChOline \Bbd C\setminus 
G_{k+1}$, $\Acc\ChOline \Bbd C\setminus G_{k}$ 
with $h(\partial D)=\partial D'$,
\par \noindent c) each component of $\Acc\ChOline 
\Bbd C\setminus G_{k+1}$ contains at most one 
point of $h^{-1}X(h)$.

Then if $G$  is in the forward orbit 
of the critical gap, $h$ is the 
\It{critical branched covering of } $f$. It is 
uniquely defined up to equivalence (given $\Gamma 
$), and $h$ has no Levy cycle, since $\Gamma $ is 
a complete Levy set. (See the Gap Restrictions of 
1.6.) If the critical gap is not in the forward 
orbit of $G$, then $k=1$ and $f^{n}=h:G\rarrow G$ 
is a homeomorphism.

\Theorem {1.10. Lemma}  Let the critical gap 
$G_{0}$ of $\Gamma $  be periodic and the 
critical branched covering $h$ be not equivalent 
to a rational map. Then  the period of $G_{0}$ is 
less than the period of any periodic critical 
point of $f.$ \endTheorem   

\Proof {Proof}  Let $n$ be the period of $G_{0}$. 
Let $G_{1}$ be the component of $f^{-2n}G_{0}$ 
with $G_{1}\subset G_{0}$, $\partial G_{0}\subset 
\partial G_{1}$. Now by 1.2, $h$ has no fixed 
critical points (in fact, no periodic critical 
points). Then for any critical point $c$, there 
is a disc component $D$ of $G_{0}\setminus G_{1}$ 
containing $c$ such that $f^{n}D$ is another 
component of $G_{0}\setminus G_{1}$ and 
$f^{i}D\cap G_{0}=\phi  $ for $0<i<n$. So if $c$ 
is periodic, $c$ must have period $>n$.

\Subheading {1.11. The Fixed Set P}  

By the Lefschetz Fixed Point Theorem, some loop 
or gap $P$ of $\Gamma $ is fixed by $f$. (This is 
the third equivalent condition of 1.6 of 
\Cite{R2} for nonequivalence to a rational map. 
This proof is due to Shishikura.) If $P$ is a 
loop, it must have orientation reversed, because 
otherwise one of the discs it bounds is fixed, 
contradicting the Levy Condition. Neither can $P$ 
be the critical gap, because if it is, any 
periodic boundary component is a degenerate Levy 
cycle for $f$. So if $P$ is a gap, $f:P\rarrow P$ 
is a homeomorphism, and boundary components are 
cyclically permuted, again because there are no 
degenerate Levy cycles. 

\Subheading {1.12}  

Now we can summarise the main homotopy-type 
information we shall use.

\Theorem {Decomposition Proposition}  

Let  $g:\Acc\ChOline \Bbd C\rarrow \Acc\ChOline 
\Bbd C$ be an orientation-preserving critically 
finite degree two branched covering which is not 
equivalent to a rational map, and with no 
degenerate Levy cycle. Let $f$ be an equivalent 
branched covering (via a homotopy which is 
constant on $X(g)=X(f)$) with complete Levy set 
$\Gamma $, critical gap $G_{0}$ and fixed loop or 
gap $P$. Let $\ell $ be the number of boundary 
components of $P$, or $\ell =2$ if $P$ is a loop.
Then one of the following holds.

\par \noindent \Rm {\Bf{Rational Critical 
Branched Covering. }} The critical gap $G_{0}$ is 
periodic under $f_{*},$   and the critical 
branched covering is equivalent to a rational 
map.

\par \noindent \Rm {\Bf{Periodic Separating Set 
}.} A gap $Q$ is mapped homeomorphically onto 
itself by $f^{\ell }$, $Q$ separates $P$ and
 $f_{*}G_{0}$, and 
$\Acc\ChOline Q\cap \Acc\ChOline P=\Acc\ChOline 
Q\cap \Acc\ChOline {f_{*}G_{0}}=\phi  $.

\par \noindent \Rm {\Bf{Adjacent Periodic Set }.} 
The gap $f_{*}G_{0}$ is periodic, and adjacent to 
some periodic gap or loop of lesser period 
than itself. \endTheorem   

\Subheading {1.13. Outline of proof of the 
Decomposition Proposition}  

Now we fix some notation. Let $\Gamma $ be a 
complete Levy set for $f$, and let $P$, $\ell $ be as 
in the statement of the proposition. Let $R$ be 
the set bounded by $P$ and $f_{*}G_{0}$ - which 
is either a closed annulus or a loop of $\Gamma 
$. Clearly, if $R$ is a loop of $\Gamma $, 
$G_{0}$ has period $\ell $, and the converse is also 
true, because no annulus component of 
$\Acc\ChOline \Bbd C\setminus \cup \Gamma $ can 
be mapped homeomorphically to itself by $f^{\ell }$, 
by the definition of complete Levy set. Let 
$\gamma =\partial R\cap \partial P$. Let $T$ be 
the local inverse of $f^{\ell }$ defined on $R$ such 
that $T$ preserves $\gamma $. Let $S$ be the 
local inverse of $f^{\ell }$ defined on $R$ which is 
$\not =T$ and such that $SR$ is in the same 
component of $\Acc\ChOline \Bbd C\setminus P$ as 
$f_{*}G_{0}$ and $R$. Now we prove the following 
lemma, under the assumptions that $R$ is a closed 
annulus and $G_{0}$ does not have period $\ell $. It 
clearly completes the proof of the Decomposition 
Proposition.

\Theorem {1.14. Lemma}  
 $TR\subset R,$  $SR\cap R\not =\phi  ,$ and 
there is a set $Q$ in $SR\cap R$ which is either 
a gap or loop of $\Gamma ,$ such that $Q$ has 
period $\ell $ under $f$ and is mapped 
homeomorphically by $f^{\ell },$ and $Q$ separates 
$f_{*}G_{0}$ from $TR.$ \endTheorem  

We split the proof of this lemma itself into a 
number of stages.

 \Theorem {1.15. Lemma}  There are least integers 
$N$, $N'\geq 0$ with
$$f_{*}^{N\ell +1}G_{0}\cap SR=\phi  \Rm {, 
}f_{*}^{N'\ell +1}G_{0}\cap TR=\phi  .$$\endTheorem   
\Proof {Proof}  The property derived from the 
definition of complete Levy set that we use 
repeatedly is:

\Item {(1)}   for all $i\geq 0$, $f_{*}^{i}G_{0}$ 
is the connected union of a component of 
$f^{-1}f_{*}^{i+1}G_{0}$ and some (possibly none) 
discs which are disjoint from $\cup \Gamma $.   

If $N$ does not exist (for example) then 
$f_{*}^{i}G_{0}$ is in the component of 
$\Acc\ChOline \Bbd C\setminus P$ bounded by 
$f^{i-1}\gamma $ for all $i$. (If $N$ does exist, 
this is true for $i\leq N\ell $.) If $G_{0}$ is 
periodic, this implies $G_{0}=f_{*}^{p\ell }G_{0}$ 
for some $p\geq 2$. But $G_{0}\cap f^{-1}R=\phi  
$, so then $f_{*}^{(p-1)\ell +1}G_{0}\cap SR=\phi  $, 
and $N$ exists. If $G_{0}$ is not periodic, there 
are integers $i$, $j\geq 1$ such that 
$f_{*}^{i\ell +1}G_{0}=f_{*}^{j\ell +1}G_{0}$ but 
$f_{*}^{i\ell }G_{0}$, $f_{*}^{j\ell }G_{0}$ 
are distinct and separated by $G_{0}$. 
Then $f_{*}^{i\ell }G_{0}$, 
$f_{*}^{j\ell }G_{0}$ cannot intersect the same 
component of $f^{-1}R$. Then
$f_{*}^{(i-1)\ell +1}G_{0}$ or 
$f_{*}^{(j-1)\ell +1}G_{0}$ is disjoint from $SR$, 
and, again, $N$ exists. The proof of the 
existence of $N'$ is exactly similar.

Now $f_{*}^{i}G_{0}\not =f_{*}^{N\ell +1}G_{0}$ for 
$1<i<N\ell +1$, so $f_{*}^{i}G_{0}$ 
($1<i\leq N\ell +1$) 
are all distinct. Hence they are all disjoint, 
and disjoint from $P$ for $0\leq i<N\ell +1$. A 
similar statement holds with $N$ replaced by 
$N'$.

\Theorem {1.16. Lemma}   $N'=0,$ so that 
$TR\subset R$ and $SR\cap R\not =\phi  ,$ and 
$S\partial f_{*}G_{0}$ separates $f_{*}G_{0}$ 
from $P.$ \endTheorem   

 \Proof {Proof}  It follows from 1.15 that

\Item {(1)}   $f_{*}^{i\ell +1}G_{0}$ is the union of 
$Tf_{*}^{(i+1)\ell +1}G_{0}$ and some discs disjoint 
from $X(f)$, for $0\leq i\leq N'-2$.    

Let $R_{0}'=R$. Inductively, define 
$R_{i+1}'=TR_{i}'\setminus f_{*}G_{0}$. Then 
$R_{i}'\not =\phi  $ for all $i$ by induction, 
since $\gamma \subset TR_{i}'$, and 
$R_{i+1}'\subset R_{i}'$. But by 1.15, using (1),
$$TR_{j}'\cap f_{*}^{i\ell +1}G_{0}\not =\phi  \Rm {, 
but }TR_{j}'\cap f_{*}^{(N'-j)\ell +1}G_{0}=\phi  $$
for $0\leq j\leq N'$, and $0\leq i<N'-j$. In 
particular, $TR_{i}'\cap f_{*}G_{0}\not =\phi  $ 
for $0\leq i<N'$, but $TR'_{N'}\cap 
f_{*}G_{0}=\phi  $, so $R'_{N'+1}=TR'_{N'}$, and 
$$\partial R_{j}'=\gamma \cup \Cup 
_{i=0}^{j}T^{i}\partial f_{*}G_{0}\Rm { for 
}j\leq N',$$
$$\partial R'_{N'+1}=\gamma \cup \Cup 
_{i=1}^{N'+1}T^{i}\partial f_{*}G_{0}.$$
See Diagram 1. The shaded region denotes 
$TR_{N'}$ (drawn in the case $N'=1$, which does 
not actually happen).

\Diagram {{250}Diagram 1}

Then $T^{N'+1}\partial f_{*}G_{0}$ must bound a 
disc disjoint from $R'_{N'+1}=TR'_{N'}$ 
containing $\partial f_{*}G_{0}$. Then all 
boundary components of $T^{k}R'_{N'}$ are 
nontrivial in $\Acc\ChOline \Bbd C\setminus X$ 
for $k\geq 0$, $T^{k+1}R'_{N'}\subset 
T^{k}R'_{N'}$ and for some $p\geq 0$, 
$T^{p+1}R'_{N'}$ and $T^{p}R'_{N'}$ are isotopic 
in $\Acc\ChOline \Bbd C\setminus X(g)$. If $N'>0$ 
then $T^{p}R'_{N'}$ has $N'+2\geq 3$ boundary 
components, all in the full orbit of $\partial 
R$. Since $\Gamma $ is a complete Levy set, 
$\partial T^{p}R'_{N'}\subset \Gamma $, and there 
is a gap $G$ for $\Gamma $ isotopic to 
$T^{p}R'_{N'}$. But then $\Acc\ChOline G\cap 
\Acc\ChOline P\not =\phi  $ with $f_{*}^{n}G\not 
=G_{0}$ for any $n\geq 0$, which is a 
contradiction to the definition of a complete 
Levy set. So $N'=0$. In particular, $TR\subset 
R$. Then $TR$ and $SR$ are separated by a 
component of $f^{-\ell }f_{*}G_{0}$ which is disjoint 
from $f_{*}G_{0}$. So $SR\cap R\not =\phi  $, and 
$S\partial (f_{*}G_{0})$ separates $f_{*}G_{0}$ 
from $P$. 

\Subheading {1.17. Proof of Lemma 1.14}  

Let $R=R_{0}$, and let $R_{1}$ be the annulus 
bounded by $S\partial (f_{*}G_{0})$ and $\partial 
f_{*}G_{0}$, so that $R_{1}\subset R=\Rm { 
domain}(S)$, and $R_{1}=R\cap SR$.
Then define
$$R_{i+1}=SR_{i}\setminus f_{*}G_{0}.$$
As in 1.16, we see that
$$SR_{j}\cap f_{*}^{i\ell +1}G_{0}\not =\phi  \Rm {, 
but }SR_{j}\cap f_{*}^{(N-j)\ell +1}G_{0}=\phi  $$
for $0\leq j\leq N$ and $0\leq i<N-j$. Hence, 
$R_{j}\not =\phi  $ for $0\leq j\leq N+1$, 
$R_{N+1}=SR_{N}$. Then if $N\geq 1$,
$$\partial R_{j}=\cup _{i=0}^{j}S^{i}\partial 
G_{0}\Rm { for }j\leq N,$$
$$\partial R_{N+1}=\cup _{i=1}^{N+1}S^{i}\partial 
G_{0}.$$
If $N=0$, $\partial R_{1}=S\partial R$.
Proceeding as in the case of $TR$, we obtain $Q$. 
If $N\geq 2$, $Q$ is isotopic to $S^{p}R_{N}$ for 
some $p\geq 0$ and is a gap for $\Gamma $ of 
period $\ell $  with boundary components cyclically 
permuted. If $N=0$ or $1$, $Q$ is isotopic to 
$\partial S^{p}R_{N}$ for $p\geq 0$,  $Q$ is part 
of a Levy cycle of period $\ell $ and has oriented 
period $2\ell $, up to isotopy. So $Q$ has the 
required properties.

\Remark {1.18 Note}  Tan Lei's proof of the 
existence of a minimal Levy cycle \Cite{TL} 
(that is, the 
existence of $P$ in the current language) was 
obtained by taking the least integer $n>1$ with 
$f_{*}^{n}G_{0}\cap P'=\phi  $, where $P'$ is the 
annulus separating $f_{*}G_{0}$ and 
$f^{-1}f_{*}G_{0}$, and then taking successive 
inverse images under $f^{n}$ of $P'\setminus \cup 
_{i<n}f_{*}^{i}G_{0}$. The proof above was 
obtained by generalising this idea.\endRemark   

\Subheading {1.19. Examples}  

It might be as well to give some examples to show 
that the various alternatives of the 
Decomposition proposition are all quite common. 
We use matings to obtain the examples, and refer 
to 1.19 of \Cite{R2} for the notations. So let 
$p$, $q$ be odd denominator rationals.Then 
$f=s_{p}\Amalg s_{q}:\Acc\ChOline \Bbd C\rarrow 
\Acc\ChOline \Bbd C$ is a type IV branched 
covering which preserves a lamination $L_{p}\cup 
L_{q}^{-1}$. If $f$ is not a rational map, then 
some leaves of $L_{p}\cup L_{q}^{-1}$ join up to 
form closed loops.

\par \noindent \It{Example 1. } If $p=q=1/3$, the 
two leaves joining the points $1/3$, $2/3$ join 
to give a loop which is, in itself, a Levy cycle, 
with the critical gap adjacent to it. We can 
easily identify the critical branched covering by 
examining $(s_{p}\Amalg s_{q})^{2}$ restricted to 
the union of the gaps containing $0$, $\infty $, 
and we see that the critical branched covering is 
equivalent to $z\mapsto z^{2}$. This gives the 
Rational Critical Branched Covering alternative 
of the Decomposition Proposition. 
\par \noindent \It{Example 2. } If $p$ and $1-q$ 
are separated from $0$ by a finite-sided gap of 
$QML$ - such as the triangle with vertices at 
${9/56}$, ${11/56}$, ${15/56}$ - but $p$ and 
$1-q$ are not in the same component of the 
finite-sided gap complement, then the critical 
branched covering is a homeomorphism.

\par \noindent \It{Example 3. } Let $p$, $q$ be 
on minor leaves which are separated from $0$ by 
the leaf with endpoints ${17\over 48}$, ${19\over 
48}$. See Diagram 2.

\Diagram {{250}Diagram 2}
As in example 1, we can identify the critical 
branched covering by examining $(s_{p}\Amalg 
s_{q})^{2}$ restricted to the union of the gaps 
containing $0$, $\infty $, and in this case, the 
critical points $c_{1}$, $c_{2}$ of the critical 
branched covering $h$ satisfy $h(c_{1})\not 
=h(c_{2})$, but $h^{2}(c_{1})=h^{2}(c_{2})$, and 
$h^{3}(c_{i})\not =h^{2}(c_{i})$ but 
$h^{4}(c_{i})=h^{3}(c_{i})$. Then $h$ is actually 
equivalent to a rational map. But we can modify 
this example as follows. First, perturb 
$s_{p}\Amalg s_{q}$ to $g$ having a complete Levy 
cycle, so that the critical gap $G$ contains a 
four-holed sphere component $G'$ of $g^{-2}G$ not 
intersecting $X(g)$. Then we can change $g$ to 
another branched covering $g_{1}$ with 
$X(g)=X(g_{1})$ by changing the definition only 
on $G'\cup gG'$. Then we can choose $g_{1}$ so 
that the new critical branched covering is 
equivalent to $[x]\mapsto [Ax]:T^{2}/\sim \rarrow 
T^{2}/\sim $, where $A$ is any integral matrix 
with determinant $2$ and all matrix entries odd, 
to ensure that the critical orbits have the same 
pattern as for $h$.

\end

\pageno=15
\Subheading {Chapter 2. The Fixed Holed-Sphere 
Theorem}  

\Subheading {2.1}  

In this chapter we prove a semiconjugacy result. 
We recall that there are many instances in 
dynamical sytems of homotopy information (usually 
containing some sort of ``expanding'' or 
``hyperbolic'' criterion) for a map $g$ yielding 
semiconjugacy information, that is, the existence 
of a nontrivial map $g_{1}$ and continuous 
$\varphi $ such that $\varphi \Circ g=g_{1}\Circ 
\varphi $. The reference which comes to mind, and 
which is actually relevant here, is \Cite{Fr}. 
There it is proved (for example) that if $T^{2}$ 
is the two-torus, and $g:T^{2}\rarrow T^{2}$ is 
continuous and homotopic to a linear automorphism 
$g_{1}$ of $T^{2}$ with both eigenvalues off the 
unit circle, then a $\varphi $ homotopic to the 
identity exists. The same result holds - and can 
be proved more simply - if $g_{1}$ is, instead, 
an expanding linear endomorphism. However, 
somewhat surprisingly, the result does not hold 
if $g_{1}$ is simply an endomorphism, even if 
both eigenvalues are off the unit circle. We can 
see this as follows. Let $g_{1}$ be such that the 
associated matrix $A$ is in $GL(2,\Bbd Z)$ and 
has determinant $2$ and eigenvalues $\lambda 
_{u}$, $\lambda _{s}$ with $\Midvert \lambda 
_{s}\Midvert <1<2<\Midvert \lambda _{u}\Midvert 
$, with corresponding eigenvectors $v_{u}$, 
$v_{s}$. Identify the two-torus with $\Bbd 
R^{2}/\Bbd Z^{2}$. We are now going to construct 
a covering map $g:T^{2}\rarrow T^{2}$ by 
constructing the well-defined lift, $\Acc\ChOline 
g$, of the multivalued function $g^{-1}$ to $\Bbd 
R^{2}$ (so that notation agrees with what we 
shall do later). In fact, we can take 
$\Acc\ChOline g$ to be an arbitrarily small 
$C^{0}$ perturbation of $A^{-1}$. We choose 
$\Acc\ChOline g$ so that $\Acc\ChOline 
g(0,0)=(0,0)$, and for some $(m,n)\in \Bbd 
Z^{2}$, $(m,n)\not =(0,0)$, the set
$$\lbrace (x,y):\lim_{{k\rarrow \infty }} 
\Acc\ChOline g^{k}(x,y)=(0,0)\rbrace $$
contains $(m,n)$. Note that the set
$$\lbrace (x,y):\lim_{{k\rarrow \infty }} 
A^{-k}(x,y)=(0,0)\rbrace $$
is simply $\Rm {sp}(v_{u})$, hence is not 
rationally defined and does not contain $(m,n)$. 
Then if Franks' theorem held for $g$, there would 
exist $\Psi :\Bbd R^{2}\rarrow \Bbd R^{2}$ 
satisfying
\par \noindent a) $\Psi ((x,y)+(p,q))=\Psi 
(x,y)+(p,q)$ for all $(x,y)\in \Bbd R^{2}$ and 
all $(p,q)\in \Bbd Z^{2}$,
\par \noindent b) $\Psi \Circ \Acc\ChOline 
g=A^{-1}\Circ \Psi $.

Note that b) yields $\Psi (0,0)=(0,0)$. Then 
$\Psi (\Acc\ChOline 
g^{k}(m,n))=A^{-k}(m,n)\rarrow 0$ as $k\rarrow 
\infty $, giving a contradiction to the existence 
of $\Psi $. Despite the fact that Franks' theorem 
does not hold in such circumstances, the 
technique of the proof can still be used to 
obtain weaker results. We shall do this in the 
case we are interested in, namely, branched 
coverings of $\Acc\ChOline \Bbd C$. Recall that, 
in Chapter 1, torus endomorphisms came up in the 
description of branched coverings of 
$\Acc\ChOline \Bbd C$. This is the reason for the 
rather involved statement of the Fixed 
Holed-Sphere Theorem below. Note that the theorem 
has two alternative conclusions. Awkward torus 
maps (such as the one above) are more involved in 
the second conclusion. The theorem can be applied 
to branched coverings which are also lamination 
maps. For the definitions of invariant 
laminations and lamination maps, and their 
purpose, see the extended introduction of 
\Cite{R2}.

\Theorem {2.2. The Fixed Holed-Sphere Theorem}  
Let $g:\Acc\ChOline \Bbd C\rarrow \Acc\ChOline 
\Bbd C$ be an orientation-preserving critically 
finite degree two branched covering, with at 
least one nondegenerate Levy cycle, and no 
degenerate Levy cycle. In a neighbourhood of any 
periodic critical point in $X=X(g)$ fixed by 
$g^{p}$, let there exist a local topological 
conjugacy to $z\mapsto z^{n}$ for suitable $n$.

Then there is an open subset $Y$ of $\Acc\ChOline 
\Bbd C\setminus X,$ and an integer $m,$ such that 
$\Acc\ChOline \Bbd C\setminus Y$ has $\ell \geq 2$ 
components, each intersecting $X$, $g^{m}\Midvert 
\Acc\ChOline Y$ is a homeomorphism, $Y$ and 
$g^{m}Y$ are isotopic in $\Acc\ChOline \Bbd 
C\setminus X$, with $g^{m}$ cyclically permuting 
the boundary components, and one of the following 
holds.

\par \noindent \Rm {\Bf{Expanding Fixed Set Case 
}}: $m=1$ and $\Acc\ChOline Y\subset gY$.

\par \noindent \Rm {\Bf{Nonrational Adjacent 
Case: }}  $m\ell $ is less than the period of any 
point in $X$ and $\partial Y=\partial _{1}\cup 
\partial _{2}$, where $\partial _{1}$ and 
$\partial _{2}$ have the following properties. 
Each $\partial _{i}$ is a finite disjoint union 
of closed arcs, two on each component of 
$\partial Y,$ $\partial _{1}\cap \partial _{2}$ 
is the set of endpoints of these arcs, and  
$\partial Y\cap g^{m}\Acc\ChOline Y\subset 
\partial _{1},$  $\partial g^{m}Y\cap 
\Acc\ChOline Y\subset g^{m}\partial _{2}.$ Each 
component of $Y\setminus g^{m}Y$ has one 
component of intersection with each of $\partial 
_{1}$ and $g^{m}\partial _{1}$, and 2 components 
of intersection with $\partial _{2}$. Each 
component of $g^{m}Y\setminus Y$ has one 
component of intersection with each of $\partial 
_{2}$ and $g^{m}\partial _{2}$, and 2 components 
of intersection with $g^{m}\partial 
_{1}$.\endTheorem   

See Diagram 3 for what $Y$ looks like in the 
nonrational adjacent case.

\Diagram {{200}Diagram 3}

\Subheading {2.3. Different Cases of, and 
Additions to, the Fixed Holed-Sphere Theorem}  
 
The different cases of the Fixed Holed-Sphere 
Theorem are dependant on the different 
conclusions of the Decomposition Proposition of 
Chapter 1. All the proofs start with the same 
basic idea, but then develop differently. Let $f$ 
be the branched covering  of the Decomposition 
Proposition which is equivalent to $g$ via a 
homotopy $\lbrace f_{t}\rbrace $, such that $f$ 
has a complete Levy set $\Gamma $, fixed loop or 
gap $P$. The Expanding Fixed Set Case of the 
Fixed Holed-Sphere Theorem holds when either of 
the first two cases of the Decomposition 
Proposition hold, that is, the Rational Critical 
Branched Covering Case, or the  Separating 
Periodic Set Case. The Nonrational Adjacent Case 
of the Fixed Holed Sphere Theorem holds when the 
critical branched covering of $f$ is nonrational 
and the Adjacent Periodic Set Case of the 
Decomposition Theorem holds.

Some unfortunately complicated variations on the 
Fixed Holed Sphere Theorem will be needed later, 
when we come to apply this result to lamination 
maps. We state these now.

\par \noindent \Bf{Variation in the Fixed Set 
Case.  }\It{ If the condition } $\Acc\ChOline 
Y\subset gY$ \It{is weakened to  } $\Acc\ChOline 
Y\subset g^{n}Y$ \It{for some  }$n$, \It{then, 
given  }$M>0$, \It{we can ensure that, for any 
arc  }$\gamma $ \It{such that  }$g^{k}\gamma 
$ \It{is of length } $\leq M$, \It{for all 
sufficiently large  }$k$, \It{then }  $\gamma 
\subset Y$.
\par \noindent \Bf{Variations in the Nonrational 
Adjacent Case. } 
\par \noindent 1. \It{Given an integer } $p$, 
\It{we can ensure that } $g^{im}(g^{m}Y\setminus 
Y)\cap Y=\phi  $ \It{for  } $0\leq i\leq p$.
\par \noindent 2. \It{Given  }$M''>0$, $Y$ 
\It{can be chosen so that the length of any path 
in  }$Y\setminus g^{m}Y$ \It{joining  }$\partial 
_{2}$ \It{and } $g^{m}\partial _{2}$\It{ is  
}$>M''$.
\par \noindent 3. \It{If } $g$ \It{is a 
lamination map preserving a lamination } $L$\It{, 
and  }$I_{i}$ $(1\leq i\leq r)$ \It{is any finite 
collection of  }$S^{1}$\It{-segments between 
noncompact leaves in gap polygons, then  
}$\partial _{1}\cup \partial _{2}$ \It{can be 
chosen to be a union of leaf segments and  
}$S^{1}$\It{-segments not including any  
}$I_{i}$.

\Subheading {2.4. Outline of proof of the Fixed 
Holed-Sphere Theorem}  

We can assume that the homotopy $f_{t}$ from $g$ 
to $f$ satisfies $f_{t}=g$ on a neighbourhood $U$ 
of the forward orbits of periodic critical 
points. Then as in the semiconjugacy proof of 
Chapter 4 of \Cite{R2}, we can find 
homeomorphisms $\Psi _{t}$ ($t\in [0,1]$) with
$$\Psi _{0}=\Rm {identity},$$
$$\Psi _{t}=\Rm {identity on }U\cup X,$$
$$f_{t}\Circ \Psi _{t}=g.$$
Then $f\Circ \Psi _{1}=g$, and $\Psi 
_{1}(g^{-1}U)=f^{-1}U$. Then define $\Psi _{t+n}$ 
($n\geq 1$) by
$$f^{n}\Circ \Psi _{t+n}=\Psi _{t}\Circ g^{n},$$
$$\Psi _{t+n}=\Rm {identity on }U\cup X.$$
Then $\Psi _{t+n+1}=\Psi _{t+n}$ on $g^{-n}U$ for 
$t\in [0,1]$. We also have
$$f\Circ \Psi _{n+1}=\Psi _{n}\Circ g,$$
since $f^{n}\Circ f\Circ \Psi _{n+1}=f^{n+1}\Circ 
\Psi _{n+1}=g^{n+1}=g^{n}\Circ g=f^{n}\Circ \Psi 
_{n}\Circ g$. Our set $Y$ in the Rational 
Critical Branched Covering and Periodic 
Separating Set Cases will be obtained from sets 
of the form $\Psi _{n}(Z)$, where $Z$ is a 
neighbourhood of $P$.

\Subheading {2.5. Proof of the Fixed Holed-Sphere 
Theorem: Rational Critical Branched Covering 
Case}  

This is the most satisfactory case to deal with. 
We define a continuous semimetric $d$ on 
$\Acc\ChOline \Bbd C$ with the following 
properties.

\par \noindent 1. There is $\delta _{0}>0$ such 
that, given $\delta _{1}>0$, there are $C$ and 
$\lambda >1$ such that:
\par \noindent  if $d(f^{i}x,f^{i}y)<\delta _{0}$ 
for $0\leq i<n$ and $d(f^{i}x,X)>\delta _{1}$ for 
$0\leq i<n$ then $d(f^{n}x,f^{n}y)\geq C\lambda 
^{n}d(x,y)$.

Note, however, that we might have $d(x,y)=0$ and 
$d(fx,fy)>0$.

\par \noindent 2. If a set $A=\lbrace 
x:d(x,x_{0})=0\rbrace $ is noncontractible, then 
it is mapped by $f$ to a set of the form $\lbrace 
y:d(y,y_{0})=0\rbrace =B$. If $A$ is contractible 
then $\partial A$ is mapped into $\partial B$ for 
some such  set $B$.

\par \noindent 3. $d(x,y)=0$ for all $x$, $y\in 
P$. If $\varepsilon >0$ is sufficiently small, 
then $\lbrace x:d(x,P)<\varepsilon \rbrace $ can 
be homotoped in $\Acc\ChOline \Bbd C\setminus X$ 
into an arbitrarily small neighbourhood of $P$.

\par \noindent 4. $\lbrace y:d(x,y)<\varepsilon 
\rbrace $ is always connected.

Suppose this semimetric has been defined. Then as 
in Chapter 4 of \Cite{R2} (that is, using 
property 1 above) we obtain from the equation 
$f^{n}\Circ \Psi _{t+n}=\Psi _{t}\Circ g^{n}$ 
that
$$d(\Psi _{n+1}(x),\Psi _{n}(x))\leq C_{1}\lambda 
^{-n},$$
for all $x$ and for some $C_{1}>0$ and $\lambda 
>1$. We use the fact that $\Psi _{n}=\Psi _{n+1}$ 
on $g^{-n}U$. Similarly, if $\gamma $ is such 
that $g^{n}\gamma $ always has length $\leq M$, 
then $\lim_{{n\rarrow \infty }} \Rm {diameter 
}\Psi _{n}(\gamma )=0$. (See the variations on 
the Fixed Holed Sphere Theorem.) Then we can find 
$\delta _{2}>0$, $N$ and $\mu >1$ such that, 
$d(X,P)>2\delta _{2}$, and if
$$d_{N}(x,y)=\sum _{i=0}^{N}d(f^{i}x,f^{i}y),$$
and $d_{N}(x,y)<\delta _{2}$, $d(x,X)>\delta 
_{0}$, then $d_{N}(fx,fy)\geq \mu d_{N}(x,y)$.
Then for $0<\varepsilon <\delta _{0}$, consider
$$Y'=\lbrace x:\lim _{{n\rarrow \infty 
}}d_{N}(\Psi _{n}(x),P)<\varepsilon \rbrace ,$$
which is automatically open. Then since, if 
$Z=\lbrace z:d_{N}(z,P)<\varepsilon \rbrace $, 
$\Acc\ChOline Z\subset fZ$, and $f\Midvert Z$ is 
a homeomorphism for $\varepsilon $ small enough, 
$\Acc\ChOline Y'\subset gY'$ and $g\Midvert Y'$ 
is a homeomorphism.
Let
$$Y=\Rm {Interior}(Y'\cup (\Rm {components of 
}\Acc\ChOline \Bbd C\setminus Y'\Rm { disjoint 
from }X)).$$
Then $Y$ has the required properties.

We can, in fact, do much more than construct $d$: 
we can construct a semiconjugacy $\pi $ between 
$g$ and $[f]$ on a certain subset of 
$\Acc\ChOline \Bbd C$, with dense image under 
$\pi $, where $[f]$ is as follows. Let $\sim $ 
denote the equivalence relation $x_{1}\sim x_{2}$ 
if and only if $d(x_{1},x_{2})=0$. Let $\pi 
:\Acc\ChOline \Bbd C\rarrow \Acc\ChOline \Bbd 
C/\sim $ be the quotient map. Let 
$[f]([x])=[f(x)]$ if $f([x])=[f(x)]$ (where $[x]$ 
is the equivalence class of $[x]$), and in 
general $[f]([x])=[f(x')]$ for any $x'\in 
\partial [x]$. By property 2 of $d$, $[f]$ is 
well-defined. Then $\pi \Circ f=[f]\Circ \pi $ on 
a subset of $\Acc\ChOline \Bbd C$. It follows 
that $\pi \Circ \Psi _{n}\Circ g=[f]\Circ \pi 
\Circ \Psi _{n+1}$ on a subset $A_{n}$ of 
$\Acc\ChOline \Bbd C$. Then we can show that $\pi 
\Circ \Psi _{n}$ converges to a map $[\Psi ]$ 
satisfying $[\Psi ]\Circ g=[f]\Circ [\Psi ]$ on 
$\limsup_{{n\rarrow \infty }}A_{n}$.

\Subheading {2.6. Construction of the 
semi-metric}  

Recall that the critical gap (see 1.8) $G_{0}$ is 
periodic under $f_{*}$, and that the critical 
branched covering is equivalent to a rational 
map, since we are in the Rational Critical 
Branched Covering case. Let us call the critical 
branched covering $h_{0}$. The critical points of 
$h_{0}$ may not be periodic. But the forward 
orbit of any periodic critical point is in 
$G_{0}$. We originally chose $h_{0}$ (in 1.8) 
with $h_{0}^{-1}G_{0}\subset G_{0}$, $\partial 
G_{0}\subset \partial h_{0}^{-1}G_{0}$. This is 
consistent with choosing $h_{0}$ so that any 
periodic critical point is locally conjugate to 
$z\mapsto z^{p}$ for suitable $p$, as we now do. 
We can also assume that $h_{0}=f^{n}\Midvert 
f^{-n}G_{0}\cap G_{0}$, if $G_{0}$ has period 
$n$. Now by the Semiconjugacy Proposition of 4.1 
of \Cite{R2}, if $h_{1}$ is the rational map 
equivalent to $h_{0}$, there is a semiconjugacy 
$\psi :\Acc\ChOline \Bbd C\rarrow \Acc\ChOline 
\Bbd C$ with $\psi \Circ h_{0}=h_{1}\Circ \psi $, 
and $\psi $ must collapse each component of 
$\Acc\ChOline \Bbd C\setminus G_{0}$ to a point. 
(We do not really need the Semiconjugacy 
Proposition. We could just be more careful in the 
choice of $h_{0}$, so that it would be clear that 
$h_{1}$ existed.) Let $d_{s}$ denote the 
spherical metric on $\Acc\ChOline \Bbd C$, and 
define
$$d_{G_{0}}(x,y)=d_{s}(\psi x,\psi y).$$
Then $d_{G_{0}}$ satisfies property 1 required 
for $d$.

Now let $G$ be any component of $\cup _{n\geq 
0}f^{-n}G_{0}$ which is not contained in a disc 
$D$ such that all components of $f^{-m}D$ are 
discs for all $m\geq 0$. Let $\Calig G_{n}$ be 
the set of those components $G$ for which $n$ is the 
least integer $\geq 0$ with $f^{n}G=G_{0}$. 
Then for $G\in 
\Calig G_{n}$, define
$$d_{G}(x,y)=(1/3^{n})d_{G_{0}}(f^{n}x,f^{n}y).$$
Note that, from the definition of $\Calig G_{n}$, 
$f^{n}:G\rarrow G_{0}$ is a homeomorphism for 
$G\in \Calig G_{n}$.

Now define $d$ as follows. For any $x$, $y\in 
\Acc\ChOline \Bbd C$, choose a path $\alpha 
:[0,1]\rarrow \Acc\ChOline \Bbd C$ with $\alpha 
(0)=x$, $\alpha (1)=y$, and such that $\alpha $ 
has at most one component of intersection with 
any $G\in \cup _{n\geq 0}\Calig G_{n}$. Then 
write
$$\alpha ^{-1}(G)=(s_{G},t_{G}),$$
if $\alpha ^{-1}(G)\not =\phi  $. Then define
$$d(x,y)=\sum _{\alpha ^{-1}(G)\not =\phi  
}d_{G}(\alpha (s_{G}),\alpha (t_{G})).$$
Then $d$ is well-defined and finite, since 
${\Sharp}(\Calig G_{n})\leq 2^{n}$, so that
$$\sum _{n=0}^{\infty }\sum _{G\in \Calig 
G_{n}}(1/3^{n})<\infty .$$
It is then easy to check that $d$ has the 
required properties.

\Subheading {2.7. Proof of the Fixed Holed-Sphere 
Theorem: A preliminary lemma in the Periodic 
Separating Set Case}  

In this case, $Q$ exists of period $\ell $ under 
$f_{*}$, separating $P$ from $f_{*}G_{0}$. See 
1.14 for the following. There is a local inverse 
$T$ of $f^{\ell }$ defined on $P\cup R$, where $R$ is 
the set bounded by $f_{*}G_{0}$ and $P$, with 
$TP=P$ and $TR\subset R$. There is a local 
inverse $S\not =T$ of $f^{\ell }$ defined on a 
neighbourhood of $Q$ with $SQ=Q$, such that $S$ 
either reverses orientation on $Q$, if $Q$ is a 
loop, or cyclically permutes the boundary 
components of $Q$ (if $Q$ is a gap) and $S^{k}$ 
fixes them, for some $k\geq 3$. Thus, $TQ$ 
separates $P$ and $Q$, and $T(\partial R\cap 
\partial f_{*}G_{0})$ - in the boundary of a 
component $G_{1}$ of $\cup _{n>0}f^{-n}G_{0}$ - 
separates $Q$ and $TQ$. 

\Theorem {Lemma}   There are $G_{i}$, $G_{i}'$ 
($i\geq 1$) between $Q$ and $TQ$, such that

\par \noindent a) $G_{i}'$ separates $G_{j}$ from 
$TQ$ for all $i$, $j$,
\par \noindent b) $G_{i}$ separates $G_{j}$ from 
$Q$, and $G_{i}'$ separates $G_{j}'$ from $TQ$, 
if $i>j$,
\par \noindent c) there is $N_{1}$ such that
$$G_{i}\Rm {, }G_{i}'\subset \Cup _{j\leq 
N_{1}+i\ell k}f^{-j}G_{0}.$$\endTheorem   
See Diagram 4.

\Diagram {{200}Diagram 4}

\Proof {Proof}  $S^{p}G_{1}$ separates $Q$ and 
$f_{*}G_{0}$ for some $1\leq p<k$. Then we can 
take $G_{i}=S^{k(i-1)}G_{1}$ and 
$G_{i}'=TS^{p+k(i-1)}G_{1}$. Thus, $G_{i}$, 
$G_{i}'$ are in the disjoint sets $SR$, $TR$ 
respectively, and have all the required 
properties.

\Subheading {2.8. Proof of the Fixed holed Sphere 
Theorem in the Periodic Separating Set Case}  

Let $Z_{0,i}$, $Z_{1}$ $Z_{2,i}$, $Z_{3}$, 
$Z_{4,i}$ $Z_{5,i}$, $Z_{6}$ be the open sets 
containing $P$ and bounded by the images under 
$f^{j}$ ($0\leq j<\ell $) of $G_{i}'$, $TQ$, 
$TG_{i}$, $T^{2}Q$, $T^{2}G_{i}$, $T^{2}G_{i}'$ 
$T^{3}Q$ respectively. Then
$$\Acc\ChOline Z_{6}\subset Z_{5,i+1}\subset 
Z_{5,i}\subset \Acc\ChOline Z_{5,i}\subset 
Z_{4,i}\subset \Acc\ChOline Z_{4,i}\subset 
Z_{4,i+1}\subset Z_{3}\subset \Acc\ChOline 
Z_{3}\subset Z_{2,i}\subset \Acc\ChOline 
Z_{2,i}$$
$$\subset Z_{2,i+1}\subset Z_{1}\subset 
\Acc\ChOline Z_{1}\subset 
Z_{0,i+1}\subset \Acc\ChOline Z_{0,i+1}\subset 
Z_{0,i}\eqno{(1)}$$
 We take
$$Y_{j,i,n}=\Psi _{n}^{-1}Z_{j,i},$$
if $j=0$, $2$, $4$ or $5$, and
$$Y_{j}=\lbrace z:\Acc\ChOline {\lbrace \Psi 
_{n}(z):n\geq N\rbrace }\subset Z_{j}\Rm { for 
some }N\rbrace $$
if $j=1$, $3$ or $6$. We shall find $N_{2}$ such 
that, if $n\Rm {, }n'\geq N_{2}+i\ell k$, for any 
$i$,
$$Y_{6}\subset \Acc\ChOline Y_{6}\subset 
Y_{5,i,n}\subset Y_{4,i,n}\subset \Acc\ChOline 
Y_{4,i,n}\subset Y_{1}\subset \Acc\ChOline 
Y_{1}\subset Y_{0,i,n'}.\eqno{(2)}$$
Then we take $Y$ to be the union of the component 
of $\Rm {Interior}(Y_{6})$ which is nontrivial in 
$\Acc\ChOline \Bbd C\setminus X$, together with 
all complementary components which do not 
intersect $X$. Then  $g^{2\ell }Y$ is the component 
of $\Rm {Interior}(Y_{1})$ which is nontrivial in 
$\Acc\ChOline \Bbd C\setminus X$, together with 
all complementary components which do not 
intersect $X$. Since $f^{2\ell }:\Acc\ChOline 
Z_{5,i}\rarrow \Acc\ChOline Z_{0,i}$ is a 
homeomorphism, so is $g^{2\ell }:\Acc\ChOline 
Y_{5,i,n+2\ell }\rarrow \Acc\ChOline Y_{0,i,n}$, and 
if (2) holds,
$$\Acc\ChOline Y\subset \Acc\ChOline 
Y_{5,i,n+2\ell }\subset Y_{4,j,n+2\ell }\subset 
g^{2\ell }Y.\eqno{(3)}$$
This is enough to show that $Y$ is nonempty and 
homotopic to $P$ in $\Acc\ChOline \Bbd C\setminus 
X$. This is the set $Y$ that we need for the 
variation 2.3 on the Fixed Holed-Sphere Theorem, 
rather than the Fixed Holed-Sphere Theorem 
itself, for which we would need to replace $Y$ by 
a union of $\cap _{i=o}^{2\ell }g^{i}Y$ and 
appropriate complementary components. 

It remains to show (2).  For this, it suffices to 
show that, for $n\geq N_{2}+i\ell k$, and $n\leq 
n'\leq n+2\ell k$,

$$Y_{5,i+2,n'}\subset Y_{5,i,n}\subset 
Y_{4,i,n}\subset Y_{4,i+2,n'}\Rm {,  
}Y_{0,i+2,n'}\subset Y_{0,i,n}.\eqno{(4)}$$

For then, 
$$\Acc\ChOline {\lbrace \Psi _{n'}(z):n'\geq 
n\rbrace }\subset \Acc\ChOline Z_{3}\subset 
Z_{1}$$
whenever $z\in Y_{4,i,n}$, and so $\Acc\ChOline 
Y\subset Y_{1}$. Also,
$$\Cap _{n\geq N}\Cap _{i}Y_{0,i,n}=\Cap _{n\geq 
N_{2}}\Cap _{i}Y_{0,i,n}$$
if $N\geq N_{2}$. So $Y_{1}\subset \Cap 
_{i}Y_{0,i,n}$ if $n\geq N_{2}$, and 
$\Acc\ChOline Y_{1}\subset Y_{0,i,n}$. To show 
(4), it suffices to find $N_{3}$ such that, for 
$n\geq N_{3}+i\ell k$, and any $x$, $\lbrace \Psi 
_{n+t}x:0\leq t\leq 2\ell k\rbrace $ crosses at most 
one set of the form $G_{j}'$ or $T^{2}G_{j 
}$ or $T^{2}G_{j}'$, for $j \leq i+2$.

To complete the proof of the Fixed Holed-Sphere 
Theorem and its variations in this case, assuming 
(4) holds, we have to show that if $\gamma $ is 
any arc such that for all sufficiently large $n$, 
has $g^{n}\gamma $ arc length $\leq M$, and 
$\gamma \cap Y\not =\phi  $, then $\gamma \subset 
Y$. It suffices to show this for $Y_{1}$ 
replacing $Y$, since $g^{2\ell }Y$ is obtained by 
adding some discs to $Y_{1}$. So fix $x\in \gamma 
$, and suppose $\Acc\ChOline {\lbrace \Psi 
_{n}(x):n\geq N\rbrace }\subset Z_{1}$ for some 
$N$. We need to show there is $N'$ such that 
$\Psi _{n}(\gamma ) \subset Z_{1}$ for all $n\geq 
N'$. But $\Acc\ChOline {\lbrace \Psi 
_{n}(x):n\geq N\rbrace }\subset Z_{2,i}$ for some 
$i$. So it suffices to find $N_{3}\geq N$ such 
that $\Psi _{n}(\gamma )$ crosses at most one set 
$TG_{s}$ for $s\leq j$ and $n\geq 
N_{3}+j\ell k$. For then $\Psi _{n}(\gamma )\subset 
Z_{2,i+2}\subset Z_{1}$ for $n\geq N_{3}+i\ell k$, 
and $\gamma \subset Y_{1}$, as required.

Now $\Psi _{n}(\gamma )$ and $\lbrace \Psi 
_{n+t}x:0\leq t\leq \ell k\rbrace $ are both 
components of the inverse image under $f^{n}$ of 
arcs of bounded length. So the proof of both (4) 
and the necessary fact about $\Psi _{n}(\gamma 
)$, and hence the proof of the theorem, is 
completed by the following lemma.

\Theorem {2.9. Lemma}  

Given $M>0,$ there is  $N_{4}$ such that the 
following holds. If $\alpha $ is any arc of 
length  $\leq M,$ then, for $n\geq N_{4},$ any 
component of $f^{-n}\alpha $ intersects at most 
one component of $\Acc\ChOline \Bbd C\setminus 
\cup _{i\geq 0}f_{*}^{i}G_{0}.$ \endTheorem   

\Proof {Proof} Recall that $\Gamma $ is a 
complete Levy set for $f$, and let $\Gamma 
_{1}=f^{-r}\Gamma $ for some suitable $r$, such 
that a nonperiodic loop of $\Gamma _{1}$ 
separates $Q_{1}$ and $Q_{2}$, whenever these are 
disjoint closed sets such that $Q_{i}$ is either 
a periodic loop in $\Gamma $ or the closure of a 
periodic gap for $\Gamma $ which is always mapped 
homeomorphically by $f^{n}$. This is clearly 
possible, because of the definition of complete 
Levy set (see Chapter 1) and because if $G_{0}$ 
is periodic, only one boundary component is 
periodic. Then it is clear that there exists $L$ 
for which the following is true. For any $\alpha 
$ of length $\leq M$, there exist $L'=L'(\alpha 
)\leq L$ and $\alpha _{i}$ ($0\leq i\leq L'$) 
with $\alpha =\cup _{i=1}^{L'}\alpha _{i}$, such 
that the arc $\alpha _{i}$ is homotopic in 
$\Acc\ChOline \Bbd C\setminus X$, via a homotopy 
fixing endpoints, to one with at most one 
intersection with a loop in $\Gamma _{1}$. Then 
there is $N_{4}$ (independent of $\alpha $) such 
that for $n\geq N_{4}$, any component of 
$g^{-n}\alpha $ intersects at most one element of 
$\Gamma _{1}$, since it intersects at most $L'$ 
elements of $g^{-n}\Gamma _{1}$. Then each one 
can intersect at most one component of 
$\Acc\ChOline \Bbd C\setminus \cup _{i\geq 
0}f_{*}^{i}G_{0}$, since the boundary of each 
$f_{*}^{i}G_{0}$ is contained in $\Gamma \subset 
\Gamma _{1}$.
\Subheading {2.10. Proof of the Fixed 
Holed-Sphere Theorem: a simpler analogue of the 
Nonrational Adjacent Case}  

Let $g:\Bbd R^{2}/\Bbd Z^{2}\rarrow \Bbd 
R^{2}/\Bbd Z^{2}$ be a covering of the two-torus 
homotopic to an endomorphism given by a matrix 
$A$ in $GL(2,\Bbd Z)$, where $A$ has determinant 
$2$ and eigenvalues $\lambda _{u}$, $\lambda 
_{s}$ with $\Midvert \lambda _{s}\Midvert 
<1<2<\Midvert \lambda _{u}\Midvert $. Let 
$v_{u}$, $v_{s}$ be the corresponding 
eigenvectors. Let $\Acc\ChOline g$ denote a 
(well-defined) lift of the multivalued function 
$g^{-1}$ to $\Bbd R^{2}$ (so that notation agrees 
with what we shall do later). Then let
$$V=\lbrace xv_{u}+yv_{s}:\Midvert x\Midvert 
<M'\rbrace ,$$
$$W=\lbrace xv_{u}+yv_{s}:\Midvert y\Midvert 
<M'\rbrace ,$$
for some suitable $M'>0$. We have $AW\subset W$, 
because $A(xv_{u}+yv_{s})=x\lambda 
_{u}v_{u}+y\lambda _{s}v_{s}$, and similarly 
$A^{-1}V\subset V$. In fact, the distance between 
the boundaries of $W$ and $AW$ increases with 
$M'$, and similarly with $V$, $A^{-1}V$. But 
there is $M$ such that $\Rm 
{distance}(\Acc\ChOline gx,A^{-1}x)\leq M$ for 
all $x$. So if $M'$ is sufficiently large given 
$M$, we can ensure that $\Acc\ChOline gV\subset 
V$, and similarly $W\subset \Acc\ChOline gW$.

Now we consider the set $\Acc\ChOline g^{n}(V\cap 
W)=Y_{n}$, so that $Y_{n+1}=\Acc\ChOline gY_{n}$. 
Let $$B_{n}=\lbrace z\in \Bbd Z^{2}:x\Rm {, 
}x+z\in Y_{n}\cup Y_{n+1} {\rm \ for\ some\ }
x\rbrace .$$
 Then $B_{n}=A^{-n}B_{0}\cap \Bbd Z^{2}$. We 
claim that $\cap _{n=0}^{\infty }A^{n}\Bbd 
Z^{2}=\lbrace 0\rbrace $. For $A$ is an 
automorphism of the subgroup $\cap _{n=0}^{\infty 
}A^{n}\Bbd Z^{2}$, which must then be either 
$\lbrace (0,0)\rbrace $ or rank two, since no 
nontrivial rationally defined subspace is 
invariant under $A$. But $A$ has determinant $2$, 
so this subgroup must indeed be $\lbrace 
(0,0)\rbrace $. (I would like to thank my 
colleagues for helping me clarify this rather 
trivial point.) So, since $B_{0}$ is finite, 
$B_{n}=\lbrace (0,0)\rbrace $ for all 
sufficiently large $n$. So for sufficiently large 
$n$, $Y_{n}\cup Y_{n+1}$ projects 
homeomorphically to $\Bbd R^{2}/\Bbd Z^{2}$, and 
$g$ maps the projection of $Y_{n+1}$ 
homeomorphically to the projection of $Y_{n}$. 
Moreover, if $Y$ denotes the projection of 
$Y_{n+1}$, we can write $\partial Y=\partial 
_{1}\cup \partial _{2}$, where each $\partial 
_{j}$ is the union of two disjoint arcs, and 
$g\partial _{1}\cap Y=\phi  $, $gY\cap \partial 
_{2}=\phi  $. This is about as far as we can take 
the arguments of \Cite{Fr}. We now want to 
generalise this argument, to prove the 
Nonrational Adjacent case of the Fixed 
Holed-Sphere Theorem.

\Subheading {2.11.Proof of the Fixed Holed-Sphere 
Theorem in the Nonrational Adjacent Case: lifting 
the branched covering inverse to the unit disc}  
We are working under the following assumptions. 
The map $g:\Acc\ChOline \Bbd C\rarrow 
\Acc\ChOline \Bbd C$ is a degree two critically 
finite orientation-preserving branched covering, 
with nonrational critical branched covering and 
an adjacent periodic set. (See the Decomposition 
Proposition.) Thus, there is a homotopy $\lbrace 
f_{t}:t\in [0,1]\rbrace $ through critically 
finite branched coverings with $f_{0}=g$, 
$f_{1}=f$, $X(f_{t})=X(f)=X$ for all $t\in 
[0,1]$, and $f$ has the following properties. 

\par \noindent \It{The set } $Q$. There is a set 
$Q$ such that 
$f^{m}Q=Q$, for some $m\geq 1$. Either $Q$ is a 
loop, in which case $f^{m}$ reverses orientation 
on $Q$ and we write $\ell =2$, or $Q$ is a $\ell $-holed 
sphere (for some $\ell \geq 3$). In the latter case, 
$f^{m}$ cyclically permutes the components of 
$\partial Q$, and $f^{q}$ fixes the boundary 
components, where $q=m\ell $. 

\par \noindent \It{The sets  }$G$, $\Acc\SwtHat 
G$. Either way, there is a set $G$ adjacent to 
$Q$ with the following properties. The sets 
$f^{i}G$ ($0\leq i\leq q$) are disjoint for 
$i<q$, and adjacent to $Q$ when $m\Midvert i$, 
with $G\subset f^{q}G$, $\partial f^{q}G\subset 
\partial G$. Also, $G$ has a holed torus as a 
double covering. There is a holed plane 
$\Acc\SwtHat G$ covering this torus, and 
$f^{-q}\Midvert G$ lifts to $\Acc\SwtHat 
f:\Acc\SwtHat G\rarrow \Acc\SwtHat G$ having the 
following properties. The map $\Acc\SwtHat f$ is 
injective, but not surjective, and there are a 
$C^{1}$ $\varphi :\Acc\SwtHat G\rarrow \Bbd 
R^{2}$ and $A\in GL(2,\Bbd Z)$, with determinant 
$2$, such that 
$$\varphi = A\Circ \varphi \Circ Acc\SwtHat f .$$
Moreover, $\varphi $ conjugates the action of the 
covering group of the torus to the action of 
$\Bbd Z^{2}$ on $\Bbd R^{2}$ by translation. The 
eigenvalues $\lambda _{u}$, $\lambda _{s}$ of $A$ 
satisfy $\Midvert \lambda _{s}\Midvert 
<1<2<\Midvert \lambda _{u}\Midvert $. We fix 
corresponding eigenvectors $v_{u}$, $v_{s}$.

\par \noindent \It{Lifting Inverses, and Paths 
Representing Points in the Universal Cover. } Now we 
recall how to lift $f_{t}^{-n}$ ($n\geq 1$) to 
the universal cover $\Acc\Overcircle {D}$ of 
$\Acc\ChOline \Bbd C\setminus X$ (This lifting of 
the inverses of critically finite branched 
coverings is fundamental in \Cite{T}.) We fix a 
basepoint $x_{1}$ in $\Acc\ChOline \Bbd 
C\setminus X$ - which, for convenience, we take 
in $Q$. Fix $n$, and fix $x_{1}'$ with 
$f_{1}^{n}x_{1}'=x_{1}$, hence fixing a path 
$\lbrace x_{t}'\rbrace $ with 
$f_{t}^{n}x_{t}'=x_{1}$. We also fix a path 
$\alpha _{1}$ from $x_{1}$ to $x_{1}'$. Then 
$\alpha _{t}$ is the path $\alpha _{1}$ followed 
by $\lbrace x_{u}'\rbrace $ from $x_{1}'$ to 
$x_{t}'$. Then $\Acc\Overcircle D$ can be 
identified with homotopy classes of paths in 
$\Acc\ChOline \Bbd C\setminus X$ based at 
$x_{1}$. If $\beta $ is any such path, let $\beta 
_{t}$ be the component of $f_{t}^{-n}\beta $ 
based at $x_{t}'$. We write
$$\Acc\ChOline f_{t}[\beta ]=[\alpha _{t}*\beta 
_{t}],$$ 
where $[\beta ]$ is the homotopy class of $\beta 
$. Then $\Acc\ChOline f_{t}:\Acc\Overcircle 
D\rarrow \Acc\Overcircle D$ is well-defined.

\par \noindent \It{The group  }$\Delta $.  If 
$\Acc\ChOline f_{t}[\beta ]=\Acc\ChOline 
f_{t}[\gamma ]$, and $[\beta ]\not =[\gamma ]$ 
then $[\beta ]=[\tau *\gamma ]=[\tau ].[\gamma 
]$, where $\tau $ is a closed loop in 
$\Acc\ChOline \Bbd C\setminus X$ based at $x_{1}$ 
(so that $[\tau ]\in \pi _{1}(\Acc\ChOline \Bbd 
C\setminus X)$). The converse is not necessarily 
true. The set of $\tau $ which arises is 
independent of $t$. Similarly, if $\Acc\ChOline 
f_{t}[\beta ]=[\tau ].\Acc\ChOline f_{t}[\gamma 
]$ for some $\tau \in \pi _{1}(\Acc\ChOline \Bbd 
C\setminus X)$, then $[\beta ]=[\sigma ].[\gamma 
]$, where $[\sigma ]\in \pi _{1}(\Acc\ChOline 
\Bbd C\setminus X)$ is again independent of $t$. 
Let $\Delta $ be the subgroup of $\pi 
_{1}(\Acc\ChOline \Bbd C\setminus X)$ consisting 
of all $[\delta ]$ such that $\Acc\ChOline 
f^{n}([\delta ].x)\in \pi _{1}(\Acc\ChOline \Bbd 
C\setminus X).\Acc\ChOline f^{n}(x)$ for all $x$, 
for all $n$. Then $\pi _{1}(Q)$, $\pi 
_{1}(\Acc\SwtHat G)\leq \Delta $. From the 
definition of $\Acc\ChOline f_{t}[\beta ]$, if 
$\Acc\ChOline f_{t}[\beta ]=[\gamma ]$, there is 
$\delta $ with initial point at $x_{t}'$, and
$$[\alpha _{t}'\cup \delta ]=[\gamma ],$$
$$[f_{t}^{n}\delta ]=[\beta ].$$

\Subheading {2.12. The construction of $V_{t}$ 
and $W_{t}$}  

From now on, we take $\Acc\ChOline f_{t}$ to be a 
lift of $f_{t}^{-m}$, with $\Acc\ChOline 
f=\Acc\ChOline f_{1}$ and $\Acc\ChOline 
g=\Acc\ChOline f_{0}$. Let $\Acc\SwtTilde  Q$ be 
a fixed lift of $Q$, and let $\Acc\SwtTilde  G$ 
be a lift of $G$ adjacent to $\Acc\SwtTilde  Q$. 
By taking $x_{1}$, $x_{1}'\in Q$ and $\alpha 
_{1}\subset Q$, we can ensure that $\Acc\ChOline 
f\Acc\SwtTilde  Q=\Acc\SwtTilde  Q$. 

Now we are ready to construct $\partial V_{t}$ by 
constructing the $\Delta $-orbit $\partial _{t}$ 
of its boundary, and then defining $V_{t}$ to be 
the component of $\Acc\Overcircle D\setminus 
\partial _{t}$ containing $\Acc\SwtTilde  Q$. We 
define $\partial _{t}=\cup 
_{i=0}^{\ell -1}\Acc\ChOline f_{t}^{i}\partial '$ 
where $\partial '$ is $\Delta $-invariant, and is 
defined by its intersection with $\Acc\SwtTilde  
G$. Since $\pi _{1}(\Acc\SwtHat G)\leq \Delta $, 
it suffices to construct the projections of the 
components to $\Acc\SwtHat G$. To do this, it 
suffices to define the images under $\varphi $ in 
$\Bbd R^{2}$, provided these are disjoint from 
$\varphi (\partial \Acc\SwtHat G)$.

We take the images under $\varphi $ of 
projections of components of $\partial '$ to 
$\Acc\SwtHat G$ to be
$$\lbrace xv_{u}+yv_{s}:x=\pm M'\rbrace .$$
The definition of $W_{t}$ is similar, with 
$v_{u}$ and $v_{s}$ interchanged.

\Theorem {2.13. Lemma}  

$\Acc\ChOline fV_{1}\subset V_{1}$ and given 
$M>0$, $M'$ can be chosen so that $d_{P}(\partial 
V_{1},\Acc\ChOline f\partial V_{1})>M$, where the 
$d_{P}$ is the Poincar\'e metric. Similarly, 
$W_{1}$ contains a component of $\Acc\ChOline 
f^{-1}W_{1}$ containing $\Acc\SwtTilde  Q$, and 
$d_{P}(\partial W_{1},\Acc\ChOline f\partial 
W_{1})>M$. \endTheorem   
\Proof {Proof} We represent points in 
$\Acc\Overcircle D$ by homotopy classes of paths 
in $\Acc\ChOline \Bbd C\setminus X$ with first 
endpoints at $x_{0}\in Q$. We shall use the same
letters to represent the lifts of the paths to the
universal cover $\Acc\Overcircle D$.
We consider the case 
of $\partial V_{1}$ - the proof for $\partial 
W_{1}$ is exactly similar. So first, let $\alpha 
$ be a path representing an element of $V_{1}$, 
with first endpoint in $Q$. Thus, $\alpha \subset 
V_{1}$. We have to show $\Acc\ChOline f\alpha 
\subset V_{1}$. By Lemma 2.14 below, for all 
sufficiently large $n$, $\Acc\ChOline f^{n}\alpha 
$ has only one component of intersection with $Q$ 
and none with any other component of 
$\Acc\ChOline \Bbd C\setminus \cup 
_{i=0}^{\ell -1}f^{im}G$ (assuming that $\alpha $ has 
only essential intersections with components of 
$\Acc\ChOline \Bbd C\setminus \cup 
_{i=0}^{\ell -1}f^{im}G$).
 For $n$ divisible by $\ell$ 
(without loss of generality), $\Acc\ChOline 
f^{n}\alpha \subset Q\cup G$, and $\Acc\ChOline 
f^{n}\alpha $ does not cross $\Acc\ChOline 
f^{n}\partial _{1}$, because $\partial 
_{1}=\Acc\ChOline f^{-n}(\Acc\ChOline 
f^{n}\partial _{1})$. In $\Acc\SwtTilde  G$ 
itself, $\Acc\ChOline f^{n}\partial _{1}$ 
separates $\Acc\ChOline f^{n-1}\partial _{1}$ 
from $\Acc\SwtTilde  Q$. So $\Acc\ChOline 
f^{n}\alpha $ also does not cross $\Acc\ChOline 
f^{n-1}\partial _{1}$, and so $\Acc\ChOline 
f\alpha $ does not cross $\partial _{1}$. Hence 
$\Acc\ChOline f\alpha \subset V_{1}$, as 
required.

To show $d_{P}(\partial V_{1},\Acc\ChOline 
fV_{1})>M$, suppose not, and let $\alpha $ be a 
path of length $\leq M$ between the boundary 
components. Again by 2.14, for some $n$ and 
$M_{1}$ depending only on $M$, $\Acc\ChOline 
f^{n}\alpha $ is a path of length $\leq M_{1}$ 
which may have one component of intersection with 
$Q$, but intersects no other component of 
$\Acc\ChOline \Bbd C\setminus \cup 
_{i=0}^{\ell -1}f^{im}G$ (again, assuming only 
essential intersections for $\alpha $). 
Suppose without loss of generality that part of
$\Acc\ChOline f^{n}\alpha$ with endpoint lifting 
to $\Acc\ChOline f^{n}\partial V_{1}$ lies in$G$,
and hence lifts to $\Acc\SwtHat G$, with endpoint
in $\varphi (\Acc\SwtHat f^{n}\partial ')$. But 
because $\Delta \cap \pi _{1}(G)=\pi 
_{1}(\Acc\SwtHat G)$, $\Acc\ChOline f^{n}\alpha $ 
cannot intersect $Q$ in just one component. So we 
can choose $n$ so that $\Acc\ChOline f^{n}\alpha 
\subset G$, and then obtain a bound (in terms of 
$M$) on the distance between components of 
$\Acc\ChOline f^{n}\partial _{1}\cap 
\Acc\SwtTilde  G$ and $\Acc\ChOline 
f^{n+1}\partial _{1}\cap \Acc\SwtTilde  G$, which 
can be contradicted by choice of $M'$.

\Theorem {2.14. Lemma}  Let $M$ be given. Then 
there is $r$ such that, if $\gamma :[0,1]\rarrow 
\Acc\ChOline \Bbd C\setminus X$ is any path of 
length $\leq M$, the following holds after 
homotopy of $\gamma $ in $\Acc\ChOline \Bbd 
C\setminus X$ fixing the endpoints. For $n\geq 
r$, a lift of $\gamma $ under $f^{n}$ has at most 
one component of intersection with $Q$, and none 
with any other component of $\Acc\ChOline \Bbd 
C\setminus \cup _{i=0}^{\ell -1}f^{im}G$.\endTheorem  

\Proof {Proof}  We can homotope $\gamma $, fixing 
endpoints, to have $m_{0}<\infty $ intersections 
with $\cup _{i=0}^{\ell -1}f^{im}\partial G$, all 
essential in $\Acc\ChOline \Bbd C\setminus X$. 
Let $m_{n}$ be the maximal number of 
intersections for a lift of $\gamma $ under 
$f^{n}$. Then since the intersections for a lift 
under $f^{n}$ are inverse images under $f^{n}$ of 
the intersections for $\gamma $, $\lbrace 
m_{n}\rbrace $ is a decreasing sequence, and 
there is $r_{1}$ such that, for $n\geq r_{1}$, 
the intersections are only with components of 
$\partial Q$. Now if $\gamma '$ is any path of 
length $\leq M_{2}$ with endpoints in $\partial 
Q$ but otherwise lying in some $f^{im}G$, there 
is $r_{2}$ (depending only on $M_{2}$) such that, 
for $n\geq r_{2}$, no lift of $\gamma '$ under 
$f^{n}$ has both endpoints in $\partial Q$. 
Applying this repeatedly, we obtain the required 
$r$.

\Subheading {2.15. Properties of $V_{t}$, 
$W_{t}$}  
Note that there is $M$ such that, for $0\leq 
i\leq \ell $,
$$d_{P}(\Acc\ChOline f^{i}x,\Acc\ChOline 
f_{t}^{i}x)<M\eqno{(1)}.$$
Then if $M'$ is sufficiently large given $M$, we 
have the following corollary of 2.13.

\Theorem {Corollary}  

 $\Acc\ChOline f_{t}^{n+1}\Acc\ChOline 
V_{t}\subset \Acc\ChOline f_{t}^{n}V_{t},$   
$\Acc\ChOline f_{t}^{n}\Acc\ChOline W_{t}\subset 
\Acc\ChOline f_{t}^{n+1}W_{t},$  $\partial 
\Acc\ChOline f_{t}^{n}V_{t}=\Acc\ChOline 
f_{t}^{n}\partial V_{t},$ and $\partial 
\Acc\ChOline f_{t}^{n}W_{t}=\Acc\ChOline 
f_{t}^{n}\partial W$  for all $t\in [0,1]$ and 
$n\geq 0.$  \endTheorem   
\Proof {Proof}  By 2.13 and (1) above, 
$\Acc\ChOline f_{t}\partial V_{t}\subset 
\Acc\ChOline V_{t}$. The last two statements 
follow from the fact that $\partial (\Delta 
.\Acc\ChOline f_{t}^{n}V_{t})=\Delta 
.\Acc\ChOline f_{t}^{n}\partial V_{t}$, and 
similarly for $W_{t}$. This is easily proved by 
induction: the point is that $\Acc\ChOline f_{t}$ 
is a covering from $\Delta .\Acc\ChOline 
f_{t}^{n}\Acc\ChOline V_{t}$ onto $\Delta 
.\Acc\ChOline f_{t}^{n+1}\Acc\ChOline V_{t}$, and 
similarly for $W_{t}$. Hence $\Acc\ChOline 
f_{t}V_{t}\subset V_{t}$ and by induction 
$\Acc\ChOline f_{t}^{n+1}V_{t}\subset 
\Acc\ChOline f_{t}^{n}V$. The proof that 
$\Acc\ChOline f_{t}^{n}W\subset \Acc\ChOline 
f_{t}^{n+1}W_{t}$ is similar, and actually 
slightly easier.

\Subheading {2.16. Proof of the Fixed 
Holed-Sphere Theorem in the Nonrational Adjacent 
Case}  

Note that $V_{1}\cap W_{1}$ has the following 
property: for some integer $N$, $\Acc\ChOline 
f^{n}(\Acc\ChOline V_{1}\cap \Acc\ChOline W_{1})$ 
projects to a compact subset of $\Acc\Overcircle 
D/\pi _{1}(Q)$ for $n\geq N$. This is simply 
because, for some $N$ the set 
$$A^{-n}(\lbrace xv_{u}+yv_{s}:\Midvert x\Midvert 
\leq M'\Rm {, }\Midvert y\Midvert \leq M'\rbrace 
)$$
does not intersect $\Bbd Z^{2}\setminus \lbrace 
(0,0)\rbrace $ for $n\geq N$. We claim that 
$\Acc\ChOline f_{t}^{n}(\Acc\ChOline V_{t}\cap 
\Acc\ChOline W_{t})$ also projects to a compact 
subset for all $t$. For
$$\Acc\ChOline f_{t}^{n}(\Acc\ChOline V_{t}\cap 
\Acc\ChOline W_{t})\subset \Acc\ChOline 
f^{n}(\Acc\ChOline V_{1}\cap \Acc\ChOline 
W_{1})\cup \Cup _{0\leq i<\ell ,t\leq u\leq 
1}\Acc\ChOline f_{u}^{i}(\psi _{u}\Acc\ChOline 
f^{n}\partial (V_{1}\cap W_{1})).$$
Here, $\psi _{u}$ is a lift of the homeomorphism 
$f_{u}^{-n}\Circ f^{n}$ to $\Acc\Overcircle D/\pi 
_{1}(Q)$. So the claim is proved. So, in 
particular, $\Acc\ChOline g^{N}(V_{0}\cap 
W_{0})=[\tau ].\Acc\ChOline g^{N}(V_{0}\cap 
W_{0})$ for $[\tau ]$ lying in only finitely many 
cosets of $\pi _{1}(Q)$. Then 2.14 implies that 
by enlarging $N$, we can assume this is true only 
for $[\tau ]\in \pi _{1}(Q)$. Take $Y$ to be the 
projection of $\Acc\ChOline g^{N}(V_{0}\cap 
W_{0})$ to $\Acc\ChOline \Bbd C\setminus X$. Let 
$\partial _{1}$ be the projection of 
$\Acc\ChOline g^{N}(\partial V_{0}\cap \partial 
(V_{0}\cap W_{0}))$, and $\partial _{2}$ the 
projection of $\Acc\ChOline g^{N}(\partial 
W_{0}\cap \partial (V_{0}\cap W_{0}))$. If we 
isotope $\partial '$ (see 2.12) a distance $\leq 
M$ in a $\Delta $-invariant way, and redefine 
$V_{t}$ using the new definition of $\partial '$, 
then $V_{t}$ still has the properties claimed 
above, and similarly for $W_{t}$. So we can 
assume that each component of $\partial V_{0}$, 
or $\Acc\ChOline g\partial V_{0}$ intersects each 
component of $\partial W_{0}$ at most once, by 
choice of $\partial '$ and the corresponding set 
for $\partial W_{0}$. Then $Y$ has all the 
properties stated in the Fixed Holed-Sphere 
Theorem.

\Subheading {2.17. Proof of the Variations in the 
Nonrational Adjacent Case}  
 For any $p>0$, the set $\cup 
_{i=0}^{p}\Acc\ChOline g^{N-i}(V_{0}\cap W_{0})$ 
intersects its image under $[\tau ]$ only for 
$[\tau ]\in \pi _{1}(Q)$ if $N$ is large enough. 
We can also make the distance in $\Acc\ChOline 
gW_{0}$ between $\partial W_{0}$ and 
$\Acc\ChOline g\partial W_{0}$ arbitrarily large 
by choice of $M'$ - and then the distance in 
$\Acc\ChOline g^{n+1}W_{0}$ between $\Acc\ChOline 
g^{n}W_{0}$ and $\Acc\ChOline g^{n+1}\partial 
W_{0}$ is even larger, for any $n>0$. If $g$ is a 
lamination map, we can choose $\partial '$ (and 
the analogue used to define $\partial W_{0}$) to 
be a union of lifts of $S^{1}$-segments and leaf 
segments, and by tracing round the boundaries of 
gaps, can avoid the lifts of the forward orbits 
of the $S^{1}$-segments $I_{i}$ (see variation 3 
of 2.3). Then $Y$ has all the required 
properties.
\end

\pageno=28
\Subheading {Chapter 3. One-dimensional 
Topological Dynamics}  

\Subheading {3.1. Definitions}  

Let $T$ be a locally finite graph. If $x\in T$, 
$U\subset T$ is a $1$-\It{sided neighbourhood of 
} $x$ if $U\setminus \lbrace x\rbrace $ is 
connected, nonempty, contains no vertex of $T$, 
and $x\in U$. (To fit in with what we need to 
prove, this definition and the following one may 
be a bit unorthodox.) Let $h:T\rarrow T$ be 
continuous with $hx=x$. Then we define 
$W^{u}(x,h,U)$ by
$$W^{u}(x,h,U)=\lbrace x\rbrace ,$$
if arbitrarily small $1$-sided neighbourhoods $V$ 
of $x$ in $U$ satisfy $hV\setminus U\not =\phi  
$, and
$$W^{u}(x,h,U)=\Cap \lbrace \Cup _{n\geq 
0}h^{n}V:V\Rm { is a }1\Rm {-sided neighbouhood 
of }x\Rm {, }V\subset U\rbrace $$
otherwise. Thus 
$W^{u}(x,h,U_{1})=W^{u}(x,h,U_{2})$ if $U_{1}$, 
$U_{2}$ are $1$-sided neighbourhoods with 
$U_{1}\cap U_{2}\not =\phi  $.

\Subheading {3.2}  

The main purpose of this chapter is to prove the 
following.

\Theorem {The $1$-sided Neighbourhood Theorem}  

Let $Y\subset \Acc\ChOline \Bbd C$ be an $\ell $ 
-holed sphere with boundary $(\ell \geq 2).$ Let 
$T\subset Y$ be a graph such that 
$T\hookrightarrow Y$ is a homotopy equivalence. 
Let $\chi :Y\rarrow Y$ be an 
orientation-preserving homeomorphism which 
cyclically permutes the components of $\partial 
Y.$ Let $h:T\rarrow T$ be continuous and 
homotopic to $\chi \Midvert T$ as a map from $T$ 
to $Y.$ 

Then there is a graph $T_{1}\subset T$ with no 
free vertices such that $\Acc\ChOline \Bbd 
C\setminus T_{1}$ has $k$ components (for some 
$2\leq k\leq \ell ),$  $h\Midvert T_{1}$ is homotopic 
as a map from $T_{1}$ to $T$ to some 
$h_{0}:T_{1}\rarrow T_{1}$ which cyclically 
permutes the boundaries of the components of 
$\Acc\ChOline \Bbd C\setminus T_{1}$ up to 
homotopy, and
 $$T_{1}\subset \Cup \lbrace 
W^{u}(x,h^{k},U):h^{k}x=x\Rm {, }x\in T_{1}\Rm {, 
}U\Rm { is a }1\Rm {-sided }$$  $$\Rm 
{neighbourhood of }x\Rm { in }T_{1}\rbrace $$ 
\endTheorem   

\Subheading {3.3}  

The main tool in the proof of the $1$-sided 
Neighbourhood Theorem is the Pseudo-Anosov 
Proposition below. Some of the same ideas 
occurred in the proof of the Tuning Proposition 
in Chapter 7 of \Cite{R2}, and, as there, we rely 
heavily on Thurston's description of surface 
homeomorphisms up to isotopy, as described in 
\Cite{F-L-P}.

\Theorem {Pseudo-Anosov Proposition}  

Let $Y_{1}$ be a $k$ -holed sphere  $(k\geq 2).$ 
Let $\chi :Y_{1}\rarrow Y_{1}$ be a homeomorphism 
which cyclically permutes the boundary components 
of $Y_{1}$ and does not leave isotopically 
invariant any nontrivial nonboundary disjoint 
loop set in $Y_{1}.$ Suppose also that $\chi $ 
fixes a point $x_{0}$ in $Y_{1}.$ Then $\chi $ is 
isotopic, via an isotopy fixing $x_{0},$ to $\chi 
_{0}:Y_{1}\rarrow Y_{1}$ satisfying one of the 
following.

\par \noindent 1.  $\chi _{0}^{k}=\Rm 
{identity}.$ 

\par \noindent 2.  $\chi _{0}$ fixes a disc $D$ 
round $x_{0},$ and $\chi _{0}:Y_{2}\rarrow Y_{2}$ 
is pseudo-Anosov, where $Y_{2}=Y_{1}\setminus D.$ 
In particular, a measured foliation  on $Y_{2}$ 
with singularities is preserved by $\chi _{0},$ 
and this foliation has  singular leaf $\ell $ 
with endpoint at $x_{1}\in \partial Y_{1}\subset 
\partial Y_{2}$ such that $x_{1},$  $\ell $ are 
preserved by $\chi _{0}^{k}.$ Moreover, there is 
a length measurement defined on $\ell ,$ and 
$\lambda >1,$ such that $\chi _{0}^{k}$ 
multiplies length on $\ell $ by $\lambda ,$ and 
the length of any segment of $\ell $ is 
proportional to the diameter of a lift to 
$\Acc\SwtTilde  Y_{2},$ where $\Acc\SwtTilde  
Y_{2}$ is the universal cover of $Y_{2}.$ 
\endTheorem   

\Proof {Proof}  To follow the outline of 
\Cite{F-L-P} more closely, we assume (as we may 
do) that $\chi $ preserves $D$. The Teichmuller 
space of $Y_{2}=Y_{1}\setminus D$ is then the 
space of (non-complete) hyperbolic structures on 
$Y_{2}$ for which all the boundary components are 
geodesics of length $1$. This space is  
homeomorphic to the open unit ball in $\Bbd 
R^{2k-4}$ (or a single point if $k=2$) and its 
Thurston compactification is homeomorphic to the 
closed unit ball in $\Bbd R^{2k-4}$. Exactly as 
in \Cite{F-L-P}, $g$ induces a homeomorphism of 
this Thurston compactification (the only 
difference is that, here, $\chi $ does not fix 
boundary components) and must have a fixed point. 
Since no simple disjoint nontrivial nonboundary 
loop set in $Y_{1}$ is left invariant by $\chi $ 
up to isotopy, and $\chi $ cyclically permutes 
the components of $\partial Y_{1}$, there is 
similarly no such loop set in $Y_{2}$. So there 
is $\chi _{0}$ isotopic to $\chi $ such that 
either

\par \noindent 1. $\chi _{0}$ is an isometry, for 
some hyperbolic structure on $Y_{2}$, or

\par \noindent 2. $\chi _{0}$ is pseudo-Anosov, 
in a sense made more precise below (and see 
\Cite{F-L-P}).

\par \noindent \It{Case 1. } We claim that $\chi 
_{0}^{k}=\Rm {identity on }Y_{2}$, and we can 
then extend $\chi _{0}$ to $Y_{1}$ with $\chi 
_{0}^{k}=\Rm {identity}$ on $Y_{1}$. For let 
$\gamma $ be a geodesic arc of minimal length 
joining $\partial D$ to $\partial Y_{1}$. Then 
$\gamma $ is disjoint from all other such arcs. 
It follows that $\chi _{0}^{i}\gamma $ ($0\leq 
i<k$) are disjoint, and $\chi _{0}^{k}\gamma 
=\gamma $. Since $Y_{1}\setminus (D\cup \Cup 
_{i\geq 0}\chi _{0}^{i}\gamma )$ is a disc, $\chi 
_{0}^{k}=\Rm {identity}$ on $Y_{2}$.

\par \noindent \It{Case 2. } We are only really 
interested in the unstable foliation preserved by 
$\chi _{0}$. This singular foliation has finitely 
many singularities, at least one on each 
component of $\partial Y_{2}$, and the complement 
of the singularities in $\partial Y_{2}$ is a 
finite set of leaves. The set of singularities is 
also preserved by $\chi _{0}$. There is a length 
measurement on leaves in $\Rm {Interior}(Y_{2})$ 
(induced by the unstable foliation being 
transverse to the stable), and there is $\lambda 
>1$ such that $\chi _{0}$ multiplies length by 
$\lambda $. Let $x_{i}$ ($1\leq i\leq r$) be the 
singularities. There is an integer $p_{i}\geq 3$ 
and a chart
$$\lbrace z:\Midvert z\Midvert <1\rbrace \Rm { if 
}x_{i}\in \Rm {Interior}(Y_{2}),$$
$$\lbrace z:\Midvert z\Midvert <1\Rm {, 
Re}(z)\geq 0\rbrace \Rm { if }x_{i}\in \partial 
Y_{2},$$
with $x_{i}$ corresponding to $0$ in the chart, 
such that the leaves of the foliation are locally 
the level sets of the multivalued function
$$\Rm {Re}(z^{p_{i}/2})\Rm { if }x_{i}\in \Rm 
{Interior}(Y_{2}),$$
$$\Rm {Re}(z^{p_{i}-1})\Rm { if }x_{i}\in 
\partial Y_{2}.$$
Then by the Poincar\'e-Hopf Theorem,
$$\sum _{i=1}^{r}({p_{i}\over 
2}-1)=k-1\eqno{(1)}$$
Since the number of singularities on each 
component of $\partial Y_{1}$ is the same 
(because $\chi _{0}$ preserves the set of 
singularities) it follows from (1) that there 
must be exactly one singularity on each component 
of $\partial Y_{1}$, and these are of the form 
$\chi _{0}^{i}x_{1}$ ($0\leq i<k$) for some 
$x_{1}\in \partial Y_{1}$ with $\chi 
_{0}^{k}x_{1}=x_{1}$, and the $p_{i}$ associated 
to $x_{1}$ must be $3$, so that exactly one leaf 
$\ell $ in $\Rm {Interior}(Y_{2})$ ends at 
$x_{1}$, and $\chi _{0}^{k}\ell =\ell $. Also, 
from (1) and the preservation of the foliation by 
$\chi _{0}$, we deduce that $\ell $ meets no 
other singularities. Also, since all infinite 
half-leaves of the unstable foliation in $\Rm 
{Interior}(Y_{2})$ are dense, $\Acc\ChOline \ell 
=Y_{2}$. Since any measured foliation is 
recurrent, there is $M_{0}$ such that every 
segment of $\ell $ of length $M_{0}$ meets an 
element of $\Gamma $, where $\Gamma $ is a fixed 
collection of disjoint simple closed curves on 
$Y_{2}$ cutting $Y_{2}$ into pairs of pants. It 
follows that the length of any segment $\ell 
_{1}$ of $\ell $ is proportional to the diameter 
of any lift of it to $\Acc\SwtTilde  Y_{2}$.

\Subheading {3.4. Proof of the $1$-sided 
Neighbourhood Theorem: finding $Y_{1}$, $T_{1}$, 
$x_{0}$, $\chi _{0}$}  

Since, by the Lefschetz Theorem, $\chi $ has a 
fixed point in $Y$, we can choose a connected 
subsurface $Y_{1}$ of $Y$ which is invariant 
under $\chi $ up to isotopy, and such that no 
nontrivial nonboundary simple disjoint loop set 
in $Y_{1}$ is preserved by $\chi $ up to isotopy. 
Then we can assume $\chi Y_{1}=Y_{1}$. Then $\chi 
$ cyclically permutes the components of $\partial 
Y_{1}$,  since $\chi $ cyclically permutes the 
components of $Y\setminus Y_{1}$. We can also 
assume there is a subgraph $T_{1}\subset T$ with 
no free vertices such that $T_{1}\subset Y_{1}$ 
and $T_{1}\hookrightarrow Y_{1}$ is a homotopy 
equivalence. So there is $p:Y_{1}\rarrow T_{1}$ 
such that $p\Midvert T_{1}$ is homotopic to the 
identity. Let $\Acc\SwtTilde  Y$, $\Acc\SwtTilde  
T$, $\Acc\SwtTilde  Y_{1}$, $\Acc\SwtTilde  
T_{1}$ be the universal covers of $Y$, $T$, 
$Y_{1}$, $T_{1}$ with $\Acc\SwtTilde  
T_{1}\subset \Acc\SwtTilde  Y_{1}\subset 
\Acc\SwtTilde  Y$ and $\Acc\SwtTilde  
T_{1}\subset \Acc\SwtTilde  T\subset 
\Acc\SwtTilde  Y$. Let $d$ be a fixed metric on 
$\Acc\SwtTilde  Y$ (lifted from a metric on $Y$) 
invariant under the action of the fundamental 
group $\pi _{1}(Y)$. Let $\Acc\SwtTilde  h$, 
$\Acc\SwtTilde  \chi $, $\Acc\SwtTilde  p$ be the 
lifts of $h$, $\chi $, $p$ so that a homotopy 
between $p\Circ \chi $ and $h$ (as maps from 
$T_{1}$ to $T$) lifts to a homotopy between 
$\Acc\SwtTilde  h$ and $\Acc\SwtTilde  p\Circ 
\Acc\SwtTilde  \chi $. Then there is $M_{1}>0$ so 
that 
$$d(\Acc\SwtTilde  h(x),\Acc\SwtTilde  p\Circ 
\Acc\SwtTilde  \chi (x))\leq M_{1}\Rm { for all 
}x\in \Acc\SwtTilde  T.$$
Since $\Acc\SwtTilde  h(\Acc\SwtTilde  T_{1})$ is 
connected, for any $x_{1}$, $x_{2}\in 
\Acc\SwtTilde  T_{1}$, $\Acc\SwtTilde  
h(\Acc\SwtTilde  T_{1})$ contains the unique arc 
between $\Acc\SwtTilde  hx_{1}$ and 
$\Acc\SwtTilde  hx_{2}$. Then we claim that
$$\Acc\SwtTilde  T_{1}\subset \Acc\SwtTilde  
h\Acc\SwtTilde  T_{1}.$$
For given $y\in \Acc\SwtTilde  T_{1}$, we can 
find $x_{1}$, $x_{2}\in \Acc\SwtTilde  T_{1}$ 
such that $\Acc\SwtTilde  p\Circ \Acc\SwtTilde  
\chi (x_{1})$, $\Acc\SwtTilde  p\Circ 
\Acc\SwtTilde  \chi (x_{2})$ are separated by $y$ 
and both distance $>M_{1}$ from $y$. Then $y$ 
also separates $\Acc\SwtTilde  hx_{1}$ and 
$\Acc\SwtTilde  hx_{2}$. So $y\in \Acc\SwtTilde  
hT_{1}$. Then a continuous map $\Acc\SwtTilde  
h_{0}:\Acc\SwtTilde  T_{1}\rarrow \Acc\SwtTilde  
T_{1}$ is uniquely determined by the properties
$$\Acc\SwtTilde  h_{0}=\Acc\SwtTilde  h\Rm { on 
}\Acc\SwtTilde  T_{1}\cap \Acc\SwtTilde  
h^{-1}\Acc\SwtTilde  T_{1},$$
$$\Acc\SwtTilde  h_{0}\Rm { maps each component 
of }\Acc\SwtTilde  T_{1}\setminus \Acc\SwtTilde  
h^{-1}\Acc\SwtTilde  T_{1}\Rm { to a point.}$$
Then $\Acc\SwtTilde  h_{0}$ descends to a map 
$h_{0}:T_{1}\rarrow T_{1}$ which is homotopic to 
$\chi :T_{1}\rarrow Y_{1}$. So by the Lefschetz 
Fixed Point Theorem, $h_{0}$ has a fixed point 
$x_{0}$. Then we can assume that $\chi 
(x_{0})=x_{0}$, and that $\chi \Midvert T_{1}$ 
and $h_{0}:T_{1}\rarrow Y_{1}$ are homotopic via 
a homotopy fixing $x_{0}$. Then we can find 
$\chi _{0}:Y_{1}\rarrow Y_{1}$ as in the 
Pseudo-Anosov Proposition.

\Subheading {3.5. Proof of the $1$-sided 
Neighbourhood Theorem: a reformulation}  

It suffices to prove, for some lift $H$ of 
$h^{k}$ to $\Acc\SwtTilde  T$,
 $$\Acc\SwtTilde  T_{1}\subset \Cup \lbrace \tau 
.W^{u}(x,H,U):x\in \Acc\SwtTilde  T_{1}\Rm {, 
}Hx=x\Rm {, }U\Rm { is a }1\Rm 
{-sided}\leqno{(3.5)}$$ 
$$\Rm {neighbourhood of }x\Rm { in }\Acc\SwtTilde 
 T_{1}\Rm {, and }\tau \in \pi _{1}(T_{1})\rbrace 
.$$

For, projecting down, this proves the theorem.

\Theorem {3.6. Lemma}  

It suffices to prove: 
\Item {(3.6)}   for some closed arc  $I$ in 
$\Acc\SwtTilde  T_{1}$ projecting down to cover 
$T_{1}$ and some lift $H'$ of $h_{0}^{k},$  
$I\subset \Cup \lbrace W^{u}(x,H',U):x\in I\Rm {, 
}H'x=x\Rm { and }U\Rm { is a }1\Rm {-sided}$   
$\Rm {neighbourhood of }x\Rm { in }I\rbrace .$   
\endTheorem   

\Proof {Proof}  Suppose (3.6) holds. Then let 
$\tau \in \pi _{1}(T_{1})$ be such that $H'=\tau 
.\Acc\SwtTilde  h_{0}^{k}$, and let $H=\tau 
.\Acc\SwtTilde  h^{k}$. We claim that a), b) and 
c) hold. Clearly they suffice to prove (3.5), as 
is required.

\par \noindent a) $W^{u}(x,H',U)\subset 
W^{u}(x,H,U)$ if $H'x=x$ and $\Acc\SwtTilde  
h_{0}^{i}x=\Acc\SwtTilde  h^{i}x$ for $1\leq 
i\leq k$,

\par \noindent b) $W^{u}(x,H',U)=\lbrace x\rbrace 
$ if $H'x=x$ and $\Acc\SwtTilde  h_{0}^{i}x\not 
=\Acc\SwtTilde  h^{i}x$ for some $1\leq i\leq k$,

\par \noindent c) if, again, $\Acc\SwtTilde  
h_{0}^{i}x\not =\Acc\SwtTilde  h^{i}x$ for some 
$1\leq i\leq k$, then $\sigma .x\in 
W^{u}(x',H,U')$ for some $x'\not =\sigma .x$, 
$Hx'=x'$ and $1$-sided neighbourhood $U'$ of $x'$ 
in $I$, and $\sigma \in \pi _{1}(T_{1})$.

To see a), note that if $V$ is any connected set 
in $\Acc\SwtTilde  T$, $\Acc\SwtTilde  
h_{0}V\subset \Acc\SwtTilde  hV$. We see b) as 
follows. If $x\in \Acc\SwtTilde  T_{1}$ and 
$H'x=x$ but $\Acc\SwtTilde  h^{i}x\not 
=\Acc\SwtTilde  h_{0}^{i}x$ for some $0<i\leq k$, 
then we can find inductively $x_{j}\in 
\Acc\SwtTilde  T_{1}$ ($0\leq j<k$) with 
$x_{0}=x$, $\Acc\SwtTilde  hx_{j+1}=\Acc\SwtTilde 
 h_{0}x_{j+1}=x_{j}$. Then for some $\zeta \in 
\pi _{1}(T_{1})$, $\zeta .x_{k}\not =x$ but 
$H(\zeta .x_{k})=H'(\zeta .x_{k})=x$. Here,  
$\zeta $ is determined by $\tau $, and by the 
automorphism of $\pi _{1}(T_{1})$ which is itself 
determined by the homotopy class of 
$\Acc\SwtTilde  h$. Moreover, $W^{u}(x,H',U)$ is 
a point (for any $1$-sided neighbourhood $U$ of 
$x$), because, if $i$ is the least integer $>0$ 
with $\Acc\SwtTilde  h^{i}x\not =\Acc\SwtTilde  
h_{0}^{i}x$, then $i$ is also the least integer 
$>0$ with $\Acc\SwtTilde  h^{i}x\notin 
\Acc\SwtTilde  T_{1}$, and there is a 
neighbourhood $V$ of $x$ with $\Acc\SwtTilde  
h_{0}^{i}V=\Acc\SwtTilde  h_{0}^{i}x$, and then 
$H'V=x$. 
Now for c): by (3.6) there is $x'\in I$ and $\eta 
\in \pi _{1}(T_{1})$ such that $\eta .x_{k}\in 
W^{u}(x',H',U')$ for some $1$-sided neighbourhood 
of $x'$ in $I$. Then $H(\eta .x_{k})=H'(\eta 
.x_{k})=\sigma .x\in W^{u}(x',H',U')$ for some 
$\sigma \in \pi _{1}(T_{1})$. Now $\eta 
.x_{k}\not =\sigma .x$ because $\Acc\SwtTilde  
h^{i}(\eta .x_{k})\in \Acc\SwtTilde  T_{1}$ for 
$0\leq i\leq k$, but this is not true for $\sigma 
.x$, because it is not true for $x$. So $H'(\eta 
.x_{k})\not =\eta .x_{k}$ and $\eta .x_{k}\not 
=x'$. So $W^{u}(x',H',U')$ is not a point, so by 
b) $\Acc\SwtTilde  h_{0}^{i}x'=\Acc\SwtTilde  
h^{i}x'$ for $0\leq i\leq k$ and 
$W^{u}(x',H',U')\subset W^{u}(x',H,U')$ by a).
  
\Subheading {3.7}  There are two cases of the 
$1$-sided Neighbourhood Theorem to prove, 
corresponding to cases 1 and 2 of the 
Pseudo-Anosov Proposition. In case 2, recall that 
$\chi _{0}$ fixes a disc $D$ round $x_{0}$, and 
that $Y_{2}=Y_{1}\setminus D$. Then we choose a 
graph $T_{2}\subset Y_{2}$ to be of the form 
$(T_{1}\setminus D)\cup \partial D$, so that 
$T_{2}\hookrightarrow Y_{2}$ is a homotopy 
inclusion, and we choose maps $h_{2}:T_{2}\rarrow 
T_{2}$ and $p_{2}:T_{2}\rarrow T_{1}$ so that 
$p_{2}\Circ h_{2}=h_{0}\Circ p_{2}$, and $p_{2}$ 
is injective except for mapping $\partial D$ to 
$x_{0}$. Thus, $p_{2}$ extends to a map 
$Y_{1}\rarrow Y_{1}$ homotopic to the identity, 
and $h_{2}$ is homotopic to $g\Midvert T_{2}$. 
Let $\Acc\SwtTilde  p_{2}$ be any lift of $p_{2}$ 
to $\Acc\SwtTilde  Y_{1}$. In case 2, it suffices 
to prove 
\Item {(3.7)}   for some closed arc  $I$ in 
$\Acc\SwtTilde  T_{2}$ projecting down to cover 
$T_{2}$, and some lift $H'$ of $h_{2}^{k}$,  
$I\subset \Cup \lbrace W^{u}(x,H',U):x\in I\Rm {, 
}H'x=x\Rm { and }U\Rm { is a }1\Rm {-sided}$   
$\Rm {neighbourhood of }x\Rm { in }I\rbrace .$  

For suppose this is true. Then, by projecting by 
$\Acc\SwtTilde  p_{2}$ we find that (3.6) is true 
also for a suitable $I$ and a suitable lift $H'$ 
of $h_{0}^{k}$.

\Subheading {3.8. A further reduction}  

We are now reduced to proving either (3.6) or 
(3.7), depending on whether we are in case 1 or 
2. We claim that it suffices to prove

\Item {(3.8)}   There is a closed arc $I$ in 
$\Acc\SwtTilde  T_{1}$ (case 1) or $\Acc\SwtTilde 
 T_{2}$ (case 2) which projects to cover all of 
$T_{1}$ (case 1) or $T_{2}$ (case 2) and such 
that $I$ has endpoints $a$, $b$, and there is a 
lift $H'$ of $h_{0}^{k}$ (case 1) or $h_{2}^{k}$ 
(case 2) such that either $H'a=a$ or $a$ 
separates $H'a$ and $b$, and either $H'b=b$ or 
$b$ separates $a$ and $H'b$, so that, in 
particular, $I\subset H'I$.   

We shall complete the proof of the $1$-sided 
Neighbourhood Theorem by proving that (3.8) holds 
in both cases 1 and 2. For the justification of 
this claim, we use the following lemma. Part (i) 
shows that $I$ contains at least one fixed point. 
Then parts (ii) and (ii) are applied to subarcs 
of $I$ between fixed points.

\Theorem {3.9. Lemma}  Let $J$ be a closed arc in 
a tree with endpoints $a,$  $b,$ and let 
$H:T\rarrow T$ be continuous. Then the following 
hold.

\par \noindent (i) If $H(a)$ is separated from 
$b$ by $a,$ and $H(b)$ is separated from $a$ by 
$b,$ then $J$ contains a fixed point of $H.$ 

\par \noindent (ii)  If $H(a)=a$ and $H(b)$ is 
separated from  $a$ by $b,$ and $J$ has no fixed 
interior point, then $J\subset W^{u}(a,H,J).$ 

\par \noindent (iii) If $H(a)=a,$  $H(b)=b$ and 
$J$ has no fixed interior point, then $J\subset 
W^{u}(a,H,J)\cup W^{u}(b,H,J).$ \endTheorem   

\Proof {Proof}  (i) Inductively we can take 
$J=J_{0}$, and find intervals $J_{n+1}\subset 
J_{n}$ with $J_{n}\subset H(J_{n+1})$, and such 
that either $J_{n}$ has a fixed endpoint for some 
$n$, or the conditions of (i) are satisfied by 
$J_{n}$ for all $n$. In this latter case, $\cap 
J_{n}$ is either an interval with fixed 
endpoints, or is a fixed point.

\par \noindent (ii) By (i), for $a'\in \Rm 
{Interior}(I)$ sufficiently near $a$, $a'$ is 
always separated from $b$ by $H(a')$. Let $c$ be 
the endpoint of $W^{u}(a,H,J)\cap J$ in 
$J\setminus \lbrace a\rbrace $. If $c\not =b$, 
then $c$ separates $H(c)$ from $b$, since 
$H(c)\not =c$. By a), this cannot happen. So 
$c=b$. Then $b\in W^{u}(a,H,J)$, since $b\in 
H(\Rm {Interior}(J))\subset W^{u}(a,H,J)$ if 
$c=b$. So $J\subset W^{u}(a,H,J)$, as required.

\par \noindent (iii) By (i), we can assume 
without loss of generality that $a'$ is separated 
from $b$ by $H(a')$ for $a'\in \Rm {Interior}(J)$ 
sufficiently near $a$. Let $c$ be the endpoint of 
$W^{u}(a,H,J)\cap J$ in $J\setminus \lbrace 
a\rbrace $. If $c=b$, we are done. If not, then 
as in (ii), $c$ separates $H(c)$ from $b$. Then 
we can complete the proof by applying (ii) with 
$J$ replaced by the interval between $c$ and $b$, 
and $a$ replaced by $b$ and $b$ replaced by $c$.

\Subheading {3.10. Proof of (3.8) in Case 1}  

By choice of $\Acc\SwtTilde  h$, we can ensure 
that $\Acc\SwtTilde  h_{0}$ fixes a preimage of 
$x_{0}$ in $\Acc\SwtTilde  T_{1}$. Take 
$H'=\Acc\SwtTilde  h_{0}^{k}$, and let the 
homotopy between $\chi _{0}^{k}$ and $h_{0}^{k}$ 
lift to one between a map $K$ and $H'$. Then $K$ 
also fixes the same preimage of $x_{0}$, so 
$K=\Rm {identity}$, since $\chi _{0}^{k}\Midvert 
Y_{1}=\Rm {identity}$. So $H'$ fixes all 
preimages of $x_{0}$. We can take $I$ to be any 
arc in $\Acc\SwtTilde  T_{1}$ with endpoints at 
preimages of $x_{0}$ and projecting to cover 
$T_{1}$. 

\Subheading {3.11. Proof of (3.8) in Case 2}  

Let $\ell $, $x_{1}$ be as in case 2 of the 
Pseudo-Anosov Proposition. Let $M>0$, yet to be 
chosen. Let $\ell _{1}$ be a segment of $\ell $ 
with endpoint at $x_{1}$, such that
$$\Rm {length}(\chi _{0}^{k}(\ell _{1})\setminus 
\ell _{1})>M.$$
Let $\beta $ be the component of $\partial Y_{1}$ 
containing $x_{1}$. We can choose a lift 
$\Acc\SwtTilde  \beta $ of $\beta $ in 
$\Acc\SwtTilde  Y_{2}$ and a lift $K$ of $\chi 
_{0}^{k}$ so that all preimages of $x_{1}$ in 
$\Acc\SwtTilde  \beta $ are fixed by $K$. Let 
$x_{1}'$, $x_{1}''$ be two distinct preimages of 
$x_{1}$ in $\Acc\SwtTilde  \beta $, let $\ell 
_{1}'$, $\ell _{1}''$ be the lifts of $\ell _{1}$ 
with endpoints at $x_{1}'$, $x_{1}''$, let $a'$, 
$b'$ be the second endpoints of $\ell _{1}'$, 
$\ell _{1}''$, and let $J$ be the arc in $\ell 
_{1}'\cup \ell _{1}''\cup \Acc\SwtTilde  \beta $ 
between $a'$ and $b'$. Then $J\subset K(J)$, $a'$ 
separates $K(a')$ from $b'$, $b'$ separates 
$K(b')$ from $a'$, and by choice of $M$ given 
$M_{1}>0$, we can ensure that
$$d(K(a'),a')>M_{1}\Rm {, }d(K(b'),b')>M_{1}.$$
Since $T_{2}\hookrightarrow Y_{2}$ is a homotopy 
equivalence, there is $q:Y_{2}\rarrow T_{2}$ 
which is homotopic to the identity as a map from 
$Y_{2}$ to $Y_{2}$. Then the homotopy lifts to 
one between a map $\Acc\SwtTilde  q$ and the 
identity. Now we choose $\ell _{1}$ so large
that $q(J)=T_{2}$, which is possible by the
Pseudo-Anosov Proposition. 
Then there is $M_{2}$ such that
$$d(Kx,H'\Acc\SwtTilde  qx)\leq M_{2}\Rm { for 
all }x\in \Acc\SwtTilde  Y_{2},$$
where the homotopy between $\chi _{0}^{k}$ and 
$h_{0}^{k}$ lifts to one between $K$ and $H'$. We 
can choose $M_{1}>M_{2}$. Then the arc $I$ 
between $a=\Acc\SwtTilde  qa'$ and 
$b=\Acc\SwtTilde  qb'$ satisfies (3.8), as 
required.
\end

\pageno=35
\Subheading {Chapter 4. The Nonrational 
Lamination Map Theorem}  

\Subheading {4.1}  

The aim of this chapter is to prove the following 
theorem. For the definitions of invariant 
lamination, lamination map and so on, see Chapter 
6 of \Cite{R2}. The definition of primitive can 
be found in Chapter 7 of \Cite{R2}, and is quite 
complicated, but can mostly be forgotten.

\Theorem {The Nonrational Lamination Map Theorem} 

Let $\rho $ be a critically finite branched 
covering with a nondegenerate Levy cycle, which 
is also a lamination map, preserving a primitive 
invariant lamination $L.$ Then there exist 
integers $m,$  $k,$  finite unions $A_{0}\subset 
A_{1}\subset \rho ^{-1}A_{0}$  of leaf segments 
of $L,$ and finite unions $B_{0},$  $B_{1}\subset 
\rho ^{-1}B_{0}$  of $S^{1}$ -segments, such that 
the following hold.

\par \noindent a) Either $m=1$, or $mk$ is less 
than the minimum period of a periodic critical 
point of $\rho .$ 

\par \noindent b) Both $A_{0}\cup B_{0}$ and 
$A_{1}\cup B_{1}$  have $m$ components, each 
homotopic in $\Acc\ChOline \Bbd C\setminus X(\rho 
)$ to a $k$ -holed sphere with boundary 
components nontrivial in $\Acc\ChOline \Bbd 
C\setminus X(\rho ).$ 

\par \noindent c) Each component of  $A_{0}$ 
contains a segment mapped into itself by $\rho 
^{mk},$ with orientation preserved.

\par \noindent d) For all $n\geq 0,$ $\rho 
^{-n}B_{0}\subset S^{1}$ and diameters of 
components tend to $0$ as $n\rarrow \infty .$ 

\par \noindent e) The sets $A_{0}\cup B_{0}$ and 
$A_{1}\cup B_{1}$ are homotopic in $\Acc\ChOline 
\Bbd C\setminus X(\rho )$ via a homotopy which is 
constant on $A_{0},$ and homotopes each component 
of $(A_{1}\setminus A_{0})\cup B_{1}$ into 
$B_{0}.$ \endTheorem   

\Remark {4.2. Remark}  If the theorem holds, we 
can inductively define $A_{n}$ with $A_{n}\subset 
A_{n+1}$ and $A_{n+1}$ a union of components of 
$\rho ^{-1}A_{n}$, and $A_{n}$ converges to a 
finite union $A$ of leaves, each of which is 
fixed by $\rho ^{mk}$, and $\Acc\ChOline {\cup 
A}$ is noncontractible in $\Acc\ChOline \Bbd 
C\setminus X(\rho )$. It has $m$ components, 
unless $A$ contains compact geodesics.\endRemark  

\Subheading {4.3. Outline of Proof} 

 The theorem will be proved by choosing $Y$ for 
$\rho =g$ as in the conclusion of the Fixed 
Holed-Sphere Theorem with certain additional 
properties, and then applying the $1$-sided 
Neighbourhood Theorem to $Y$ and some embedded 
graph $T$, where the isotopy class of the 
homeomomorphism $\chi $ is given by a local 
inverse of $\rho ^{m}$. Thus we have to choose 
$Y$, $T$, $h:T\rarrow T$ with certain properties 
and deduce the theorem from this.

\Subheading {4.4. Recall of some basic 
definitions}  

We shall need to refer to \Cite{R2} to some 
extent, particularly Chapters 6 and 7. Throughout 
the present chapter, $L$ is a $\beta $-invariant 
lamination on $\Acc\ChOline \Bbd C\setminus K$, 
with lamination map $\rho $, which therefore 
preserves $L$ and the Cantor set $K$, which is 
contained in $S^{1}$. Recall that a \It{leaf 
segment } is the closure of a component of $\ell 
\setminus S^{1}$, for a leaf $\ell $ of $L$. A 
\It{polygon } is the closure $P$ of a subset of 
$\Acc\ChOline \Bbd C\setminus S^{1}$ with 
$\partial P\subset (\cup L)\cup S^{1}$. A \It{gap 
} is a component of $\Acc\ChOline \Bbd C\setminus 
(K\cup (\cup L))$. A \It{gap polygon } is a 
polygon $P$ such that $\Rm {Interior}(P)$ is a 
component of $G\setminus S^{1}$ for a gap $G$. A 
\It{side } of a gap $G$ is a leaf $\ell $ with 
$\ell \cap P\not =\phi  $ for some gap polygon 
$P$ of $G$.

It was shown in \Cite{R2} - in particular, in 
6.10 - that there are gap polygons $Q_{i}$ 
($1\leq i\leq r$) and points $a_{i}$ in $K$ 
($1\leq i\leq q$) having the following 
properties, and so that $\rho $ can be chosen to 
satisfy the following condition, given 
$\varepsilon '>0$.
$$\rho ^{-n}(S^{1}\setminus ((\Cup 
_{i=1}^{q}I_{i})\cup (\Cup _{i=1}^{r}\Rm 
{Interior}(Q_{i}))\subset S^{1}\Rm { for all 
}n\geq 0,$$
where $a_{i}$ is strictly preperiodic under $\rho 
$, and $I_{i}$ is an interval in $S^{1}\setminus 
K$ of length $<\varepsilon '$ and with endpoint 
at $a_{i}$, and $\rho (\cup Q_{i})\subset \cup 
Q_{i}$. We shall need this because we shall need 
to consider families of leaf segments and 
$S^{1}$-segments, which we shall want mapped to 
unions of leaf segments and to $S^{1}$-segments 
by inverse images of $\rho $.

\Subheading {4.5. Essential Leaf Segments}  

For suitable $\zeta >0$, we say that a leaf 
segment $\ell $ is \It{essential } if some 
components of $\rho ^{-n}\ell $ have diameter 
$>\zeta $ for infinitely many $n$. By the Length 
Lemma (\Cite{R2}, 6.15), for suitable $\zeta $, 
this is only possible if either

\par \noindent   $\Rm {diameter}(\ell )>\zeta $, 
or

\par \noindent $\Rm {diameter}(\ell )\leq \zeta $ 
and $\ell $ has endpoints on opposite sides of, 
but close to, some $a_{i}$ ($1\leq i\leq q$).

From the definition (and for suitable $\zeta $), 
if $\ell $ is not essential, then all components 
of $\rho ^{-n}\ell $ are nonessential leaf 
segments, for all $n\geq 0$.

\Subheading {4.6. Properties of $\Bf{Y }$}  

As we have said, we need to construct $Y$, $T$, 
$h$. We summarise the properties of $Y$ we shall 
need, and which are proved in 4.7-4.10. That the 
first property is possible is simply a 
restatement of part of the Fixed Holed-Sphere 
Theorem.

\par \noindent 1. $Y$ is an open $\ell $-holed 
sphere in $\Acc\ChOline \Bbd C\setminus X(\rho )$ 
for some $\ell \geq 2$, with each boundary 
component nontrivial, and for some local inverse 
$S$ of $\rho ^{m}$, $Y$ and $SY$ are isotopic in 
$\Acc\ChOline \Bbd C\setminus X(\rho )$, and $S$ 
cyclically permutes boundary components up to 
isotopy in $\Acc\ChOline \Bbd C\setminus X(\rho 
)$.

\par \noindent 2. Either (i) or (ii) holds, 
depending on whether the Expanding Fixed Set or 
Nonrational Adjacent case of the Fixed 
Holed-Sphere Theorem holds.

\par \noindent (i) $m=1$ and $SY\subset Y$.

\par \noindent (ii)  Here, $m$ is such that 
$m\ell $ less than the minimum period of a 
periodic critical point of $\rho $, and $\partial 
Y=\partial _{1}Y\cup \partial _{2}Y$, where 
$S\Acc\ChOline Y\cap \partial _{1}Y=\phi  $, and 
$S\partial _{2}Y\cap \Acc\ChOline Y=\phi  $. Each 
component $C$ of $SY\setminus Y$ intersects only 
one component of each of $\partial _{2}Y$, 
$S\partial _{2}Y$.

For (i) see 4.7, and (ii) is a restatement of the 
Fixed Holed-Sphere Theorem.

\par \noindent 3. Either (i) or (ii) holds, 
depending, again, on the case of the Fixed 
Holed-Sphere Theorem.

\par \noindent (i) $\partial Y$ is a finite union 
of leaf segments and $S^{1}$-segments.

\par \noindent (ii) $\partial _{1}Y$ is a finite 
union of leaf segments and $S^{1}$-segments, 
$\partial _{2}Y$  is a finite union of 
$S^{1}$-segments.

For (i), see 4.7. For (ii), see 4.9. 

\par \noindent \It{Note. } We shall write 
$S^{n}Y$ even when $SY\setminus Y\not =\phi  $, 
meaning the local inverse defined on $S^{n-1}Y$ 
by lifting the inductively defined homotopy 
between $S^{n-2}Y$ and $S^{n-1}Y$.

\par \noindent 4. $S^{i}Y\cap \lbrace 
a_{1}...a_{q}\rbrace =\phi  $ for all $i\geq 0$, 
and if $b$ is a component of $\partial Y\cap 
S^{1}$, then $b\cap \Cup _{i=1}^{r}Q_{i}=\phi  $ 
so that $S^{n}b\subset S^{1}$ for all $n$.

See 4.7.

\par \noindent 5. $$\Acc\ChOline Y=\Cup 
_{i=1}^{v}Y_{i},$$
where the sets $\Rm {Interior}(Y_{i})$ are 
disjoint, and, for some $0\leq t\leq u\leq v$, 
the following hold.

\par \noindent a) $\partial Y_{i}\setminus 
\partial Y\subset S^{1}$, and if $b$ is a 
component of $\partial Y_{i}\cap S^{1}$, $b\cap 
\Cup _{i=1}^{r}Q_{i}=\phi  $, so that 
$S^{n}b\subset S^{1}$ for all $n$.

\par \noindent b) For $i>t$, $Y_{i}$ is a single 
polygon not containing any periodic gap polygon 
such that $\partial Y_{i}$ contains exactly two 
essential leaf segments for $t<i\leq u$, and none 
for $u<i\leq v$.

\par \noindent c) For $i\leq t$, $Y_{i}$ is a 
finite connected union of polygons, and there is 
a finite connected union $Y_{i}'\subset Y_{i}$ of 
gap polygons such that each maximal polygon of 
$Y_{i}$ contains exactly one polygon of $Y_{i}'$, 
$(Y_{i}\setminus Y_{i}')\cap S^{1})\cap \Cup 
_{i=1}^{r}Q_{i}=\phi  $, and there is a bijection 
$\tau :\lbrace 1...t\rbrace \rarrow \lbrace 
1...t\rbrace $ with
$$Y_{\tau (i)}'\subset SY_{i}'.$$

\par \noindent d) $SY_{j}\cap Y_{i}=\phi  $ for 
$i$, $j\leq t$ and $i\not =\tau (j)$, and 
$SY_{j}\cap Y_{\tau (j)}$ is connected.

See 4.10 (and 4.7-4.9).

\Theorem {4.7. Lemma}  

 $Y$ can be chosen so that either 1 or 2 holds, 
where $S$ is a local inverse of $\rho ^{m}.$ 

\par \noindent 1. $m=1,$  $\partial Y$ is a 
finite union of leaf segments and $S^{1}$ 
-segments, and $SY\subset Y.$ 

\par \noindent 2.  $\partial Y=\partial _{1}Y\cup 
\partial _{2}Y,$ where $\partial _{i}Y$ is a 
finite union of leaf segments and $S^{1}$ 
-segments,
 $$S\Acc\ChOline Y\cap \partial _{1}Y=\phi  \Rm 
{, }S\partial _{2}Y\cap \Acc\ChOline Y=\phi  ,$$
and each component of $SY\setminus Y$ has one 
component of intersection with each of $\partial 
_{2}Y$, $S\partial _{2}Y$. 

Furthermore, we can ensure that $\partial Y\cap 
S^{1}\cap Q_{i}=\phi  $ $(1\leq i\leq r)$ and 
$S^{n}Y\cap \lbrace a_{1}...a_{q}\rbrace =\phi  $ 
for all $n\geq 0.$   \endTheorem   

\Proof {Proof}  

Suppose the Expanding Fixed Set case of the Fixed 
Holed-Sphere Theorem holds. Then the variation of 
2.3 implies $Y$ can be chosen to contain all 
leaves and finite-sided gaps of $L$ it 
intersects, and so that $SY\subset Y$ and 
$S^{p}\Acc\ChOline Y\subset Y$ for some $p$. For 
any finite connected union of leaf segments, or 
arc in a finite-sided gap has the property of 
$\gamma $ in the statement of the variation. Then 
we can approximate $S^{i}Y$ ($0\leq i<p$) by a 
set $Y_{i}$ whose boundary is a finite union of 
leaf segments, components of intersection with 
gap polygons in infinite-sided gaps and their 
sides, and $S^{1}$-segments, such that 
$\Acc\ChOline Y_{i+1}\subset Y_{i}$ for $i\leq 
p-2$, $S\Acc\ChOline Y_{p-1}\subset Y_{0}$ and 
$\partial Y_{i}\cap \partial Q_{j}\cap S^{1}=\phi 
 $. Then $Y'=\cap _{i=0}^{p-1}Y_{i}$ satisfies 
$S\Acc\ChOline Y'\subset Y'$. Let $Y$ be obtained 
from $Y'$ by removing all components of 
intersection with gap polygons in infinite-sided 
gaps and their sides which are not complete gap 
polygons. Then $SY\subset Y$. We can ensure 
$S^{n}Y\cap \lbrace a_{1}...a_{q}\rbrace =\phi  $ 
for all $n\geq 0$ by replacing $Y$ by $S^{n}Y$ 
for a suitable $n\geq 0$, because the $a_{i}$ are 
strictly preperiodic.

Suppose the Nonrational Adjacent case of the 
Fixed Holed-Sphere Theorem holds. Then the 
variations of 2.3 ensure that we can choose $Y$ 
so that $\partial _{1}Y$, $\partial _{2}Y$ 
consist of finitely many leaf and 
$S^{1}$-segments, and do not intersect $S^{1}\cap 
Q_{i}$ for any $i$. Then $\partial _{1}Y$, 
$\partial _{2}Y$ have all the required 
properties, provided, as in case 1, we replace 
$Y$ by $S^{n}Y$ if necessary.

\Theorem {4.8. Lemma}  

 $Y$ can be chosen so that the periodic gap 
polygons contained in $S^{n}Y$ are the same for 
all $n\geq 0,$ and such that no two periodic gap 
polygons $P,$  $Q$ in $Y$ lie in the same 
component of $Y\setminus S^{1}.$ \endTheorem

\Proof {Proof}  Since the lamination $L$ is 
primitive, the total number of periodic gap 
polygons is finite. (See 6.15, 7.2 and 7.3 of 
\Cite{R2}. A periodic thin gap polygon lies in a 
periodic gap bounded by periodic segments on some 
of the finitely many periodic sides. But in a 
primitive lamination, there is a bound on the 
number of periodic segments on any leaf.) If 
$SY\subset Y$, the number of periodic gap 
polygons in $S^{n}Y$ obviously decreases with 
$n$,  so replacing $Y$ by a suitable $S^{n}Y$, we 
can assume the number is always constant. If $Y$ 
satisfies the Nonrational Adjacent case of the 
Fixed Holed-Sphere Theorem, so that $SY\setminus 
Y\not =\phi  $, let $p$ be the maximal period of 
a periodic gap polygon. If some periodic gap 
polygon is in $SY\setminus Y$, then for some 
$1\leq n\leq p$, $S^{n}(SY\setminus Y)\cap Y\not 
=\phi  $ - but, by the variations of 2.3  of the 
Fixed Holed-Sphere Theorem, $Y$ can be chosen so 
that this does not happen.

Finally, suppose for contradiction that $P$, $Q$ 
are in the same component of $S^{n}Y\setminus 
S^{1}$ for all $n$. Then $P$, $Q$ are 
finite-sided, since $S^{n}\Midvert Y$ is a 
homeomorphism for all $n$. We can assume there 
are no periodic gap polygons  between $P$ and $Q$ 
in the same component of $Y\setminus S^{1}$. But 
the region between $P$ and $Q$ must be mapped 
back into itself by some $S^{n}$, so it must be a 
single periodic leaf segment, which is 
impossible, because gaps are never immediately 
adjacent in a primitive lamination. So, again, 
replacing $Y$ by $S^{n}Y$ if necessary, $Y$ has 
the required properties.

\Theorem {4.9. Lemma}  

In the case when $\partial Y=\partial _{1}Y\cup 
\partial _{2}Y,$  $Y$ can be chosen so that 
$\partial _{2}Y\subset S^{1}.$ \endTheorem   

\Proof {Proof} Let $C$ be a component of 
$SY\setminus \Acc\ChOline Y$. Then by the 
variations of 2.3, the distance between 
$\Acc\ChOline C\cap \partial _{2}Y$ and 
$\Acc\ChOline C\cap S\partial _{2}Y$ can be 
assumed so large that no maximal polygon in $C$, 
or union of polygons intersecting some $Q_{i}$, 
has sides in both $\partial _{2}Y$ and $S\partial 
_{2}Y$. Then we can find an $S^{1}$-arc in $C$, 
which intersects no $Q_{i}$, separating 
$\Acc\ChOline C\cap \partial _{2}Y$ and 
$\Acc\ChOline C\cap S\partial _{2}Y$. Then we 
replace $Y$ by $SY$ cut off by $S^{1}$-arcs in 
all such components $C$.

\Subheading {4.10}  

4.7 to 4.9 show that $Y$ can be chosen satisfying 
1 to 4 of 4.6, and part of 5 of 4.6: we can write 
$\Acc\ChOline Y=\cup _{i=1}^{r}Y_{i}$ so that 
$Y_{i}$ is a single polygon not containing a 
periodic gap polygon for $i>t$, and $Y_{i}$ has 
the properties listed in 5 of 4.6 for $i\leq t$, 
except for 5b) and 5d). Now we claim that, if we 
replace $Y$ by $S^{n}Y$ for a suitable $n\geq 0$, 
and rename and renumber the $Y_{i}$, but keeping 
the name $Y_{i}$ ($i\leq t$) for the union of 
polygons of $\Acc\ChOline Y$ intersecting 
$Y_{i}'$,  then 5 of 4.6 is satisfied in its 
entirety. For gap polygons are dense and all 
eventually periodic, so we can ensure that, for 
$i>t$, $Y_{i}$ does not contain any specified gap 
polygon, by choosing $n$ large enough. So then, 
given $\zeta >0$, we can ensure that $\partial 
Y_{i}$ has at most two leaf segments of length 
$>\zeta /2$, and if at most one has length 
$>\zeta /2$, they all have length $<\zeta $. Thus 
we can ensure that $\partial Y_{i}$ contains 
either exactly two essential leaf segments, or 
none, for $i>t$, and, renumbering, 5a)-c) of 4.6 
are satisfied. we can also assume that $\cup 
_{i=1}^{t}(Y_{i}\setminus Y_{i}')\cap \cup 
_{i=1}^{t}SY_{i}'=\phi  $, so that 5d) is also 
satisfied, as required.

\Subheading {4.11. Properties of $\Bf{U }$, 
$\Bf{U' }$}  

Let $U\subset \Acc\ChOline U\subset U'$ be open 
sets in $\Acc\ChOline Y$ satisfying the following 
- which are possible.

$$\Cup _{i>u}Y_{i}\subset U,\leqno{a)} $$
and, in the Nonrational Adjacent case,
$$\partial _{2}Y\subset U.$$
\par \noindent b) For $t<i<u$, $Y_{i}\cap U$ has 
two components, and $\partial (Y_{i}\setminus 
U)\setminus \partial U$ has two components, which
 contain $\Acc\ChOline  Y_{i}\cap S^{1}$, and 
the essential leaf segments in $\partial Y_{i}$
run between these two components.
\par \noindent c) For $i\leq t$, each component 
of $Y_{i}\cap U$ contains exactly one component 
of $\partial Y_{i}\cap S^{1}$.

\par \noindent d) Similar properties to a), b) 
and c) hold for $U'$.
$$S\Acc\ChOline U\cap Y\subset U'.\leqno{e)}$$
\par \noindent f) Each component of $SY\setminus 
Y$ intersects a unique component of $U'$.

\Subheading {4.12. Properties of $\Bf{T }$}  

We can define a graph $T\subset Y$ satisfying the 
following.

\par \noindent a) $T\hookrightarrow Y$ is a 
homotopy equivalence.

\par \noindent b) $T$ has no vertex in 
$Y_{i}\setminus U$ for $t<i\leq u$, but has one 
edge crossing between the components of $U$ 
bounding $Y_{i}$.

\par \noindent c) For $i\leq t$, $T\cap 
Y_{i}\subset Y_{i}'$, $T$ has exactly one vertex 
$v_{i}$ in $Y_{i}'\setminus U$ with $Sv_{i}\in 
Y_{\tau (i)}'\setminus U'$, and exactly one edge 
of $T$ runs from $v_{i}$ into each bounding 
component of $U$. 

\par \noindent d) No arc of $ST\setminus U'$ has 
both endpoints in the same component of $U$. 

\Subheading {4.13. Definition of $\pi $, $\Bf{h 
}$}  

We can define $\pi :ST\rarrow T$ satisfying the 
following.

\par \noindent a) $\pi $ maps each component of 
$ST\cap U\setminus \cup _{i=1}^{t}Y_{i}$ into the 
component of $U'$ it is contained in, and each 
component of $ST\setminus Y$ into the component 
of $U'$ it intersects.

\par \noindent b) $\pi $ maps $ST\setminus 
(U'\cup \Cup _{i=1}^{t}Y_{i})$ to $T\setminus U$, 
and is monotone restricted to each arc in 
$ST\setminus (U\cap \Cup _{i=1}^{t}Y_{i})$.
$$\pi Sv_{i}=v_{\tau (i)}.\leqno{c)} $$
\par \noindent d) Any component of $ST\cap (\Cup 
_{i=1}^{t}Y_{i}\setminus U)$ whose closure 
intersects only one component of $U'$  is mapped 
by $\pi $ into that component of $U'$.
For $i\leq t$, $\pi $ is monotone restricted to 
each arc in $ST\cap Y_{i}\setminus U'$ which has 
endpoints in different components of $U'$, and 
maps it to an arc in $Y_{i}\cap T$ with endpoints 
in the same components of $U$.

Then we define
$$h=\pi \Circ S.$$

\Subheading {4.14. Properties of $\Bf{h }$}  

From the definition of $h$, we obtain the 
following.
$$h(U)\subset U.\leqno{a)} $$
\par \noindent b) If $t<j\leq u$ and $i\leq u$, 
then $h(T\cap Y_{i}\setminus U)$ contains $T\cap 
Y_{j}\setminus U$ if it intersects it, This 
happens if and only if $SY_{i}\cap Y_{j}$ 
separates the essential leaf segments in 
$\partial Y_{j}$.

\par \noindent c) For $i\leq t$, 
$h(v_{i})=v_{\tau (i)}$, and if $j\leq t$ also,
$$h(T\cap Y_{i}\setminus U)\cap Y_{j}=\phi  ,$$
unless $j=\tau (i)$, when $h^{-1}(Y_{j}\cap 
T\setminus U)\cap T\cap Y_{i}\setminus U$ is 
connected, and then, if $b$ is a component of 
$Y_{i}\setminus U\setminus \lbrace v_{i}\rbrace 
$, $h(b)$ contains exactly one component of 
$Y_{j}\setminus U\setminus \lbrace v_{j}\rbrace $ 
and intersects no other.

Now the $1$-sided Neighbourhood Theorem will 
imply the existence of $k$ with $2\leq k\leq \ell 
$ and a subgraph $T_{1}$ of $T$ with
$$T_{1}\setminus U\subset \Cup \lbrace 
W^{u}(x,h^{k},V): x\in T_{1}\setminus U\Rm {, 
}h^{i}x\in T_{1}\Rm { for all }i\Rm {, 
}h^{k}x=x\Rm {,}\leqno{d)} $$
$$V\Rm { is a }1\Rm {-sided neighbourhood of 
}x\Rm { in }T_{1}\rbrace .$$
Then the properties of $\pi $ imply that

\par \noindent e) if $h^{k}x=x$ and 
$W^{u}(x,h^{k},V)\not =\phi  $ for some $x\in 
T_{1}\cap Y_{i}\setminus U$, some $i\leq t$ and 
for some $V$, then $b\subset S^{k}b$ for each 
component $b$ of $\partial Y_{i}'\cap (\cup L)$. 
Moreover, given $i\leq t$, there must be such an 
$x$, $V$, unless every arc of $T_{1}\cap 
Y_{i}\setminus U'$ is in $h(T_{1}\setminus 
Y_{i})$, that is, unless, for any maximal arc $b$ 
of $T_{1}\cap Y_{i}$, there is an arc of 
$S(T_{1}\cap Y\setminus Y_{i})\cap Y_{i}$ running 
parallel to it throughout its length.

\Subheading {4.15. Proof of the Nonrational 
Lamination Map Theorem}  

We shall find $A_{0}$ such that exactly one 
connected leaf piece from $A_{0}$ runs parallel 
to each arc in $T_{1}\cap (\Cup 
_{i=1}^{u}Y_{i}\setminus U)$. Then we shall join  
the endpoints of the leaf pieces of $A_{0}$ by 
$S^{1}$-arcs whose union we denote by $B_{0}$, so 
that $A_{0}\cup B_{0}$ is homotopic to $T_{1}$ in 
$\Acc\ChOline \Bbd C\setminus X(\rho )$.  We 
shall have $\rho ^{-n}B_{0}\subset S^{1}$ for all 
$n$, since each component of $B_{0}$ is a union 
of a subset of $\cup Y_{i}\cap S^{1}$ and of 
intervals in $S^{1}$ between points of $Y_{j}\cap 
S^{1}$ ($j>u$). Thus, $\rho ^{-n}B_{0}\cap (\cup 
Q_{i}\cup \lbrace a_{1}..a_{q}\rbrace )=\phi  $ 
for all $n\geq 0$. The homotopy between 
$A_{0}\cup B_{0}$ and $A_{1}\cup B_{1}$ is 
obtained by constructing a homotopy between 
$A_{0}\cup B_{0}$ and $T_{1}$ which preserves 
$U$, $Y\setminus U$, constructing a similar 
homotopy between $ST_{1}$ and $A_{1}\cup B_{1}$, 
and constructing yet another similar homotopy 
between $\pi $ and the identity. So it remains to 
construct $A_{0}$. From properties b), d), e) of 
$h$ in 4.14, the proof of the theorem will be 
completed after the following lemma.

\Theorem {4.16. Lemma}  

If  $t<i\leq u$ and $x\in Y_{i}\cap 
T_{1}\setminus U,$  $h^{k}x=x,$  $h^{j}x\in 
T_{1}$ for $0\leq j\leq k,$ and 
$W^{u}(x,h^{k},V)\not =\phi  ,$ then there is a 
leaf segment $\ell $ in $Y_{i}$ separating the 
essential leaf segments in $\partial Y_{i}$ such 
that $\ell \subset S^{k}\ell ,$ and 
$S^{k}\Midvert \ell $ preserves 
orientation.\endTheorem   

\Proof {Proof}  

Let $x\in Y_{i}\cap T_{1}\setminus U$, 
$t<i\leq u$, 
$h^{k}x=x$, $h^{j}x\in Y_{i_{j}}\cap 
T_{1}\setminus U$ with $i_{0}=i_{k}=i$. Then 
$t<i_{j}\leq u$ for all $j$. Inductively, let 
$Q_{0}=Y_{i}$, and let $Q_{j+1}$ be the component 
of $SQ_{j}\cap Y_{i_{j+1}}$ in the component of 
$SY_{i_{j}}\cap Y_{i_{j+1}}$ which contains 
$Sh^{j}x$ for $0\leq j<k$. Then, for $n\geq k$, 
let $Q_{n+1}$ be the unique component of 
$SQ_{n}\cap Q_{n+1-k}$.  Then inductively, using 
property b) of $h$ in 4.13, we see that 
$Q_{nk+j}$ separates the essential leaf segments 
in $\partial Y_{i_{j}}$. Then $\cap Q_{nk}$ must 
be a leaf segment $\ell $ with $\ell \subset 
S^{k}\ell $, because otherwise it is a periodic 
gap polygon in $Y_{i}$, and there are none. If 
$W^{u}(x,h^{k},V)\not =\phi  $ for some $V$, then 
$S^{k}\Midvert \ell $ preserves orientation. This 
is because, for all $j\leq k$, the images under 
$S$ and $h$ of the endpoints of the arc of $T\cap 
Y_{i_{j}}\setminus U$ containing $h^{j}x$ lie in 
the same components of $U$, and the endpoints of 
the arc of $S^{j+1}\ell $ in $Q_{j+1}$ also lie 
in these same components of $U$.

\Subheading {4.17. End of proof of the Theorem}  
It suffices to choose $A_{0}\subset Y$ so that a 
leaf segment from $A_{0}$ intersects 
$Y_{i}\setminus U$ for all $t<i\leq u$, and some 
connected union of leaf segments in $A_{0}$ runs 
parallel to each component of $\partial Y_{i}'$, 
for all $i\leq t$. Now fix $i\leq u$, and suppose 
$W^{u}(x,h^{k},V)\not =\phi  $ for some $x\in 
T_{1}\cap Y_{i}\setminus U$, and $V$. If $i\leq 
t$ let $\partial Y_{i}'\cap (\cup L)\subset 
A_{0}$, and if $t<i\leq u$, let $A_{0}$ include a 
leaf segment $\ell $ as in the lemma, and let 
$A_{0}$ also contain a finite connected union of 
segments containing $\ell $ which is contained in 
$Y$, has endpoints in $U$, and runs parallel to 
as many arcs of $T_{1}\setminus U$ as possible. 
Then $A_{0}$ has the required properties.
\end

\pageno=43
\Subheading {Chapter 5. Proofs of the Admissible 
Boundary Theorem and the Admissibility 
Proposition}  

\Subheading {5.1}  

The theorems to be proved in this chapter are 
stated in Chapter 1 of \Cite{R2}, and will be 
restated here, but it will not be possible to 
explain all the terms involved, for which the 
reader should refer to Chapters 6 and 7 of 
\Cite{R2}, in particular, for a full account of 
invariant and parameter laminations. The main 
result of the last chapter, the Nonrational 
Lamination Map Theorem, was essentially a result 
about invariant laminations (or their lamination 
maps). The task, now, is to transfer this to a 
result about parameter laminations.

\Subheading {5.2. Invariant sets of complete 
geodesics}  
We recall, at least, that an invariant lamination 
$L$ is a geodesic lamination on $\Acc\ChOline 
\Bbd C\setminus K$, where $K$ is a Cantor subset 
of $S^{1}$, where $K=K_{r}$ depends on a 
parameter $r$ (which is an odd denominator 
rational). Also, there is a fixed lamination 
called (somewhat incongruously at this stage) 
$\Acc\ChOline L_{r}$ such that $L\cup 
\Acc\ChOline L_{r}$ is always a lamination, and 
$L$ contains all but countably many leaves of 
$\Acc\ChOline L_{r}$. The sets $S^{1}$,  $K$ and 
$\Acc\ChOline L_{r}$ are invariant under a 
critically finite degree two branched covering 
$\Acc\ChOline s=\Acc\ChOline s_{r}$ which fixes
 $\infty $. For 
precise definitions of $K$, $\Acc\ChOline L_{r}$, 
$\Acc\ChOline s$, see Chapter 1 of \Cite{R2}.

The invariance is actually $\beta $- invariance 
for some path $\beta :[0,1]\rarrow \Acc\ChOline 
\Bbd C\setminus (K\cup (\cup \Acc\ChOline 
L_{r}))$, and uses a homeomorphism $\sigma 
_{\beta }$ defined in terms of $\beta $. (See 
\Cite{R2} 1.7.) We recall that, in \Cite{R2} 6.5 
we put extra conditions on $\beta $ in order to 
ensure that, for a complete geodesic $\ell $ in a 
$\beta $-invariant lamination, $(\sigma _{\beta 
}\Circ \Acc\ChOline s)^{-1}\ell $ could be easily 
modified to a set of geodesics. We also need to 
do something similar here for finite geodesics, 
that is, geodesics with finitely many 
$S^{1}$-crossings, which may be incomplete. This 
means making modifications to $\beta $ which are 
similar to those of 6.5 of \Cite{R2}, but 
slightly less technical. We describe these now.

\Subheading {5.3. The Minimal Intersection 
Property}  

Let $\beta : [0,1]\rarrow \Acc\ChOline \Bbd 
C\setminus (K\cup (\cup \Acc\ChOline L_{r}))$ 
with $\beta (0)=\infty $ and let $A$ be a union 
of finite disjoint geodesic arcs in $\Acc\ChOline 
\Bbd C\setminus (K\cup (\cup \Acc\ChOline 
L_{r}))$ with any endpoints in $S^{1}$. Then 
after isotopy in $\Acc\ChOline \Bbd C\setminus 
(K\cup (\cup \Acc\ChOline L_{r}))$ keeping $\beta 
(0)=\infty $ and fixing $\beta (1)$,  $\beta $ 
can be assumed to have the following properties.
\par \noindent a) If $c$, $d$ are adjacent points 
of $\beta ^{-1}(S^{1})$, then $\beta (c)$, $\beta 
(d)$ lie in different components of 
$S^{1}\setminus K$, and $\beta $ is transverse to 
$S^{1}$.
\par \noindent b) No sub-path of $\beta $ can be 
homotoped into $A$ by a homotopy in $\Acc\ChOline 
\Bbd C\setminus K$ fixing the sub-path endpoints, 
$\beta $ is transverse to $A$, and the number of 
points of $(0,1)\cap \beta ^{-1}(A)$ cannot 
reduced by homotopies fixing $\beta (1)$ and 
$\beta (0)=\infty $, while keeping property a).
\medskip 
Let $\beta _{u}$ be any homotopy through paths in 
$\Acc\ChOline \Bbd C\setminus (K\cup (\cup 
\Acc\ChOline L_{r}))$ satisfying a) and b) such 
that $\beta _{0}=\beta $, $\beta _{u}(0)=\beta 
(0)=\infty $ and $\beta _{u}(1)=\beta (1)$  for 
all $u$, and such that $\beta _{u}(t)\in A$ only 
if $\beta (t)\in A$. Now we can state condition 
c).

\par \noindent c) If $\beta (t)\in A$, then 
$\beta (t)\notin S^{1}$, $\beta _{u}(t)\in A$ and
$${\Sharp}\lbrace v<t:\beta (v)\in S^{1}\rbrace 
\leq {\Sharp}\lbrace v<t:\beta _{u}(v)\in 
S^{1}\rbrace $$
for all $u$.

The conditions b), c) subject to a) ensure that 
$(\sigma _{\beta }\Circ \Acc\ChOline s)^{-1}A$ 
has the least possible number of intersections 
with $S^{1}$, up to homotopy fixing the endpoints 
of $(\sigma _{\beta }\Circ \Acc\ChOline s)^{-1}A$ 
and preserving $K$ and $\Acc\ChOline L_{r}$.

Then $(\beta ,A)$ \It{has the minimal 
intersection property } if conditions a) to c) 
hold. A \It{good  } isotopy will be one which 
preserves $S^{1}$, $K$ and $\Acc\ChOline L_{r}$.

 Note the following.

\par \noindent 1. The image of $(A,\beta )$ under 
any good isotopy also has the minimal 
intersection property.

\par \noindent 2. If $A$ can be isotoped into 
$A'$ by a good isotopy, then there is $\beta '$ 
such that $(A',\beta ')$ has the minimal 
intersection property, and such that $(A,\beta )$ 
can be isotoped into $(A',\beta ')$ by a good 
isotopy.

\par \noindent 3. For any finite set $A$ of 
finite geodesics in a $\beta $-invariant 
lamination with $\beta $ satisfying the 
conditions of 6.5 of \Cite{R2}, $(A,\beta )$ has 
the minimal intersection property. 

\Theorem {5.4. First Invariant-implies-Parameter 
Theorem}  
 In what follows, all inclusions and equalities 
are up to good isotopy. Let $(A_{0},\beta )$ have 
the minimal intersection property, with $\beta 
(1)$ in a component $\mu _{0}$ of $A_{0}$, and 
suppose that 
$$A_{0}\subset (\sigma _{\beta '}\Circ 
\Acc\ChOline s)^{-1}A_{0},$$ 
\par \noindent where $\beta '$ is an arbitrarily 
small extension (or restriction) of $\beta $. If 
$x\notin K$ is an endpoint of $A_{0}$, let $\tau 
(x)$ be the endpoint (if it exists) of $A_{0}$ 
such that $(\sigma _{\beta '}\Circ \Acc\ChOline 
s)^{-1}(\tau (x))$ is isotoped to $x$ by the good 
isotopy. Suppose also that, for any endpoint 
$x\notin K$, $\tau ^{n}(x)$ is not defined for 
all $n$. 

Then there are $\gamma $, and a $\gamma 
$-invariant lamination $L$, such that 
$A_{0}\subset L$, $\mu _{0}$ is contained in the 
minor leaf or a side of the minor gap of $L$, and 
$(A_{0},\beta ')$ or $(A_{0},\beta 
)=(A_{0},\gamma )$ up to good isotopy, depending 
on whether or not $L$ has a minor gap. 

The minor leaf of $L$ is eventually periodic, and 
strictly preperiodic if $L$ does not have a minor 
gap. Moreover, no half leaf of $L$ has infinitely 
many successive intersections with periodic 
components of $S^{1}\setminus K$.\endTheorem   

\Proof {Proof} Again, all inclusions and 
equalities are up to good isotopy. Write $\beta 
_{0}=\beta '$. Inductively, define $\beta _{n}$, 
$A_{n}$  for $n\geq 1$ with the following 
properties.

\par \noindent 1. $A_{n+1}=(\sigma _{\beta 
_{n}}\Circ \Acc\ChOline s)^{-1}A_{n}$

\par \noindent 2. $(A_{n},\beta _{n})\subset 
(A_{n+1},\beta _{n+1})$
\par \noindent 3. $(A_{n},\beta _{n})$ has the 
minimal intersection property.

We can clearly define $A_{1}$, $\beta _{1}$. 
Suppose $A_{k}$, $\beta _{k}$ have been defined 
with the correct properties for $k\leq n$and 
$n\geq 1$, and define $A_{n+1}$ by 1. Then since 
$(A_{n-1},\beta _{n-1})\subset (A_{n},\beta 
_{n})$, we have $(\sigma _{\beta _{n-1}}\Circ 
\Acc\ChOline s)^{-1}A_{n-1}\subset (\sigma 
_{\beta _{n}}\Circ \Acc\ChOline s)^{-1}A_{n}$, 
and hence $A_{n}\subset (\sigma _{\beta 
_{n}}\Circ \Acc\ChOline s)^{-1}A_{n}$.  Then by 
5.3.1 and 2 we can find a $\beta _{n+1}$ so that 
$(A_{n},\beta _{n})\subset (A_{n+1},\beta 
_{n+1})$.

Then the condition on $\lbrace \tau 
^{n}(x)\rbrace $ implies that, if $x$ is an 
endpoint of some $A_{n}$, then either $x\in K$ or 
$x$ is interior to $A_{m}$ for some $m>n$. It 
follows that $A_{n}$ converges to a lamination 
$L$, containing $A_{0}$. We can also assume that 
$\beta _{n}$ converges to $\gamma '$. Here, 
$\gamma '$ might contain points of $K$ which are 
endpoints of intervals of $S^{1}\setminus K$, but 
can be perturbed to $\gamma $ in $\Acc\ChOline 
\Bbd C\setminus (K\cup (\cup \Acc\ChOline L))$ 
such that $L$ is $\gamma $-invariant. Note that 
$A_{0}$ contains all the periodic segments in the 
full orbit of $A_{0}$, and $L$ is the closure of 
this full orbit. Hence no half leaf in $L$ has an 
infinite succession of periodic segments, as 
required. 

Let $\ell _{0}$ be the component of $A_{0}$ 
containing $\beta (1)$ and let $\ell _{i+1}$ be 
defined inductively as a component of $A_{i+1}$ 
which is a component of $(\sigma _{\beta 
_{i}}\Circ \Acc\ChOline s)^{-1}\ell _{i}$. 
Obviously, $\ell _{0}$ in $L$ is always in an 
eventually periodic leaf, and in a periodic leaf 
if and only if $\ell _{0}\subset \ell _{n}$ with 
orientation preserved for some $n$, for some 
choice of $\lbrace \ell _{i}\rbrace $. If this 
happens, no other component $\ell _{n}'$ of 
$A_{n}$ can run parallel to $\ell _{n}$ for the 
whole length of $\ell _{n}$ on the same side as 
$\beta _{n}(1)$. For if so, let $\ell _{i}'$ 
$(0\leq i\leq n)$ be defined similarly to $\ell 
_{i}$, so that $\ell _{i+1}'$ is a component of 
$(\sigma _{\beta _{i}}\Circ \Acc\ChOline 
s)^{-1}\ell _{i}'$. Suppose $\ell _{i+1}$ and 
$\ell _{i+1}'$ run parallel along their entire 
lengths. Then $\sigma _{\beta _{i}}^{-1}\ell 
_{i}$ and $\sigma _{\beta _{i}}^{-1}\ell _{i}'$ 
run parallel, and so do $\ell _{i}$ and $\ell 
_{i}'$. Then, by induction, $\sigma _{\beta 
'}^{-1}\ell _{0}$ and $\sigma _{\beta '}^{-1}\ell 
_{0}'$ run parallel along their entire lengths. 
This is impossible, because $\beta '(1)$ 
separates $\ell _{0}$ and $\ell _{0}'$. It 
follows that $\beta _{k}$ can be chosen so that 
$(A_{k},\beta _{k})$ converges to $(L,\gamma ')$, 
modifying to $(L,\gamma )$,
with $\gamma (1)$ in a minor gap of $L$. Thus, 
$L$ has all the required properties.

\Subheading {5.5}  
We recall that, in \Cite{R2}, $\Acc\SwtTilde  U$ 
denoted the universal cover of the union of 
component of $\Acc\ChOline \Bbd C\setminus (K\cup 
(\cup \Acc\ChOline L_{r}))$ containing $\infty $ 
and some  boundary leaves, with $0\in 
\Acc\SwtTilde  U$ projecting to $\infty $. If 
$\beta :[0,1]\rarrow \Acc\ChOline \Bbd C\setminus 
(K\cup (\cup \Acc\ChOline L_{r}))$ is a path with 
$\beta (0)=\infty $, then $\Acc\SwtTilde  \beta $ 
denotes the lift to $\Acc\SwtTilde  U$ of $\beta 
$ with $\Acc\SwtTilde  \beta (0)=0$. Given a 
$\beta $-invariant lamination $L$ with minor gap, 
the lift to $\Acc\SwtTilde  U$ containing 
$\Acc\SwtTilde  \beta (1)$ is also called the 
minor gap, and similarly for minor leaves: the 
minor leaf in $\Acc\SwtTilde  U$ either contains 
$\Acc\SwtTilde  \beta (1)$, or is intersected by 
$\Acc\SwtTilde  \beta $ and in the boundary of 
the minor gap, if this exists. 

We recall the Parameter Laminations Theorem (1.16 
of \Cite{R2}) which says that minor leaves and 
boundaries of minor gaps of primitive invariant 
laminations never intersect transversally. 
Primitivity is no real restriction. As in the 
statement and proof of the First 
Invariant-implies-Parameter Theorem, inclusions 
and equalities for sets of geodesic pieces are up 
to good isotopy, in the following theorem.
\Theorem {  Second Invariant-implies-Parameter 
Theorem}  

Let $(A_{0},\beta )$ and $\beta '$ be as in the 
First Invariant-implies-Parameter Theorem. Let 
$L'$ be a primitive $\gamma '$-invariant 
lamination and let $\Acc\SwtTilde  \gamma '$ be 
the lift of $\gamma '$ with $\Acc\SwtTilde  
\gamma '(0)=0$. Suppose that there is an isotopy 
preserving $K$, $\Acc\ChOline L_{r}$ and $A_{0}$ 
of $\beta '$ to $\gamma '$. Suppose the minor gap 
of $L'$ - if it exists - is not both periodic and 
nonsimply-connected. Then any transverse 
intersections between $A_{0}$ and the minor leaf 
or sides of the minor gap of $L'$ can be removed 
by good isotopy, and one of the following holds.

\par \noindent 1. Each component of $A_{0}$ is 
contained in $L'$ or a finite-sided gap of $L'$.
\par \noindent 2. There is a $\gamma $-invariant 
lamination $L$ containing $A_{0}$ and with all 
the other properties of the First 
Invariant-implies-Parameter Theorem, such that 
$L'\subset L$. We can choose $L$ so that either 
of the following holds. The minor leaves of $L$, 
$L'$ are the same, or at least one component of 
$A_{0}$ is contained in the forward orbit of the 
minor leaf or sides of the minor gap of 
$L$.\endTheorem   
\Proof {Proof}   The proof is similar to that of 
the First Invariant-implies-Parameter Theorem, 
but we start by fixing a finite set $C_{0}$ of 
finite geodesics in the forward orbit of sides of 
the minor gap of $L'$ - if this exists - or in 
the minor leaf of $L'$, such that $(C_{0},\gamma 
')$ has the minimal intersection property, 
$$C_{0}\subset (\sigma _{\gamma '}\Circ 
\Acc\ChOline s)^{-1}C_{0},$$
and the minor leaf and sides of the minor gap (if 
the minor gap exists) are in the full orbit of 
$C_{0}$. Now we construct inductively 
$A_{n}\subset A_{n+1}$, $C_{n}\subset 
C_{n+1}$, with
$$A_{n+1}=(\sigma _{\beta _{n}}\Circ \Acc\ChOline 
s)^{-1}A_{n}\Rm {, }C_{n+1}=(\sigma _{\beta 
_{n}}\Circ \Acc\ChOline s)^{-1}C_{n},$$
where the definition of $\beta _{n}$ is slightly 
different from that in the First Theorem in some 
cases, and we do not define $C_{m}$ for $m\geq n$ 
if $C_{n}$ intersects $A_{n}$ transversally. For 
$n\geq 0$, $\beta _{n}$ is always defined so that 
either $(A_{n},\beta _{n})$  or $(A_{n}\cup 
C_{n},\beta _{n})$ has the minimal intersection 
property, depending on whether or not $C_{i}$ 
intersects $A_{i}$ transversally for some $i\leq 
n$. Similarly, $\beta _{n}$ and $\beta _{n+1}$ 
have the same intersections with $A_{n}$ or 
$A_{n}\cup C_{n}$. Also, $\beta _{n}(1)$  is 
always in the minor gap of $L'$. The lamination 
$L$ is then obtained by taking the limit of 
$A_{n}$ or $A_{n}\cup C_{n}$, depending on 
whether or not $C_{n}$ intersects $A_{n}$ 
transversally. The Parameter Laminations Theorem 
(1.16 \Cite{R2}) and (7.1, 7.3 \Cite{R2}) then 
show that, in fact, $C_{n}$ does not intersect 
$A_{n}$ transversally unless the minor gap of 
$L'$ is periodic nonsimply-connected. Every leaf 
of $L$ is contained in $L'$ or a finite-sided gap 
of $L$ unless the minor gap of $L'$ is periodic 
(and hence infinite-sided).

There is a certain amount of choice in the 
construction. If some $A_{m}$ intersects $C_{m}$ 
transversally or the minor gap of $L'$, then for 
some $N$, $\beta _{n}$ must always be close to 
the same component of $C_{N}\subset C_{n}$ or 
$A_{N}\subset A_{n}$ for $n\geq N$. Otherwise, 
$\beta _{n}$ must be close to the same component 
of $C_{N}\subset C_{n}$ for $n\geq N$ for some 
$N$. We can choose the component of $A_{N}$ to 
contain a component of $A_{0}$, if the minor gap 
of $L'$ is periodic. We can satisfy these 
conditions, so we are done.
 
\Subheading {5.6. Piecewise-geodesic Fixed 
Holed-spheres}  

Recall that $\Acc\ChOline s$ is a critically 
finite branched covering, and $X(\Acc\ChOline 
s)\setminus \lbrace \infty \rbrace =X_{1}(\rho 
)=\lbrace \rho ^{i}(0):i\geq 0\rbrace $  for any 
branched covering lamination map $\rho $ (for an 
invariant lamination containing $\Acc\ChOline 
L_{r}$).

Given an invariant lamination $L$ with lamination 
map $\rho $, let $A_{0}$, $B_{0}$ be the sets 
proved to exist in the Nonrational Lamination Map 
Theorem when $L$ is primitive and $\rho $ has a 
nondegenerate Levy cycle. It is clear, by 
considering the boundary of the set $A_{0}\cup 
B_{0}$, that the converse is true: if $\rho $ is 
a critically finite branched covering, and such 
that $A_{0}$, $B_{0}$ exist, then $\rho $ has a 
nondegenerate Levy cycle.  The important thing to 
remember is that, if $m$ is the period of $0$ 
under $\Acc\ChOline s$ (and $\rho $) then all 
segments of $A_{0}$ are eventually periodic of 
eventual oriented period $\leq m$, and if 
$A_{0}\cup B_{0}$ is not connected, the eventual 
oriented period is strictly less than the period 
of any periodic critical point. (I apologise
for the change in use of $m$. We now revert to the
letter conventions established in Chapter 1 of 
\Cite{R2}. The lettering in Chapters 1 to 4 of this 
paper is internally consistent as far as I can manage
it, and that is the best I can do.) We shall call 
$A_{0}\cup B_{0}$ a \It{piecewise-geodesic fixed 
holed-sphere }, if it is connected, and, in 
general, a \It{piecewise-geodesic periodic 
holed-sphere. } Recall that $A_{n}\cup B_{n}$ is 
the component of $\rho ^{-n}(A_{0}\cup B_{0})$ 
containing $A_{0}$, and the properties satisfied 
by $A_{0}\cup B_{0}$ imply that $\Acc\ChOline 
{\cup A}$ is noncontractible in $\Acc\ChOline 
\Bbd C\setminus X_{1}(\rho )=\Acc\ChOline \Bbd 
C\setminus (X(\Acc\ChOline s)\setminus \lbrace 
\infty \rbrace )$.
 
 We also have the following, using the same 
notations as in the First 
Invariant-implies-Parameter Theorem and its 
proof.
\Theorem {Lemma} Let $(A_{0},\beta )$ and $L$ be 
as in the First Invariant-implies-Parameter 
Theorem. Let $\rho $ be the lamination map of 
$L$. Also, let $B_{0}$ be a finite union of 
intervals in $S^{1}\setminus K$ with endpoints in 
$A_{0}$, such that the forward orbits of the 
components of $S^{1}\setminus K$ intersected by 
$\beta $ are not intersected by $B_{0}$.
Then $\rho ^{-n}(B_{0})\subset S^{1}$ for all 
$n\geq 0$.

Similarly, let $(A_{0},\beta )$, $L'$, with 
lamination map $\rho '$, be as in the Second 
Invariant-implies-Parameter Theorem, with the 
same hypotheses as above. If the first conclusion 
of that theorem holds, $\rho '^{-n}(B_{0})\subset 
S^{1}$ for all $n$. If the second conclusion 
holds, $L$ can be chosen so that, in addition to 
the other properties, $\rho '^{-n}(B_{0})\subset 
S^{1}$ for all $n$.\endTheorem   
\Proof {Proof} In both the First and Second 
Theorems, $(A_{n},\beta _{n})$ has the minimal 
intersection property. It suffices that, for all 
$n$, $\beta _{n}\cap B_{n}=\phi  $, where 
$B_{n+1}=(\sigma _{\beta _{n}}\Circ \Acc\ChOline 
s)^{-1}B_{n}$. In the case of the First Theorem, 
this is obvious, by the condition on $\beta $, 
since $B_{n}=(\sigma _{\beta _{n}}\Circ 
\Acc\ChOline s)^{-1}B_{n-1}=(\sigma _{\beta 
}\Circ \Acc\ChOline s)^{-n}B_{0}=\Acc\ChOline 
s^{-n}B_{0}$. In the case of the Second Theorem, 
it is also clear, unless some component of 
$A_{n}$ is in the minor gap of $L'$. In that 
case, we consider the forward orbits of those 
components of $S^{1}\setminus K$ which have to be 
crossed from the end of $\beta $ to reach that 
component. Then there is $N>n$ such that 
$\Acc\ChOline s^{-N}B_{0}$ intersects none of 
these components of $S^{1}\setminus K$. (Note 
that, by hypothesis, all components of $B_{0}$ 
lie in strictly preperiodic components of 
$S^{1}\setminus K$.) Take $\beta _{k}$ with 
endpoint close to the same side of $L'$ for all 
$k<N$. Then, by the choice of $N$, we can choose 
$\beta _{N}$ having the same intersections with 
$A_{N-1}$ as $\beta _{N-1}$, so that $\beta _{N}$ 
does not intersect the forward orbit of $B_{N}$ 
and $\beta _{N}(1)$ is close to a component of 
$A_{N}$ in the minor gap. Then we take the end of 
$\beta _{k}$ close to this component of $A_{N}$ 
for all $k\geq N$.

\Subheading {5.7. Moving to a Primitive 
Lamination}  

In Chapter 7 of \Cite{R2} we showed how to modify 
a lamination with eventually periodic minor gap 
(or leaf) to a primitive lamination. The 
corresponding lamination maps are then 
equivalent, and if one has a nondegenerate Levy 
cycle, they both do. Thus, by the Nonrational 
Lamination Map Theorem, if one has a 
piecewise-geodesic periodic holed-sphere, they 
both do. In fact, it is relatively easy to obtain 
the holed sphere for one directly from the other 
if one remembers the definition of primitive, but 
we omit the details. 

\Subheading {5.8. Admissible Points} 

The definition of admissible given in 1.17 of 
\Cite{R2} was slightly incorrect. There are 
several possible rather similar definitions. We 
give one here which avoids certain complications 
in the subsequent proof. The definition depends 
on a parameter lamination $\Calig L=\Calig L_{r}$ 
on $\Acc\SwtTilde  U$, in which all leaves are 
either minor leaves or sides of minor gaps of 
invariant laminations containing all but 
countably many leaves of $\Acc\ChOline L_{r}$. 
\par \noindent \It{Types II, III and IV. } Let 
$f$ be a degree two branched covering with 
critical points $c_{1}$ and $c_{2}$, with $c_{1}$ 
periodic. Recall from \Cite{R2} 1.2 that $f$ is 
type II if $c_{2}$ is in the same periodic orbit 
as $c_{1}$, type III if $c_{2}$ is in the full 
orbit of $c_{1}$ but not periodic, and type IV if 
$c_{2}$ is in a distinct periodic orbit. This 
includes all critically finite degree two 
branched coverings such that every critical point 
has a periodic critical point in its forward 
orbit.

A leaf $\Acc\SwtTilde  \mu $ in $\Calig L$ is 
\It{admissible } if it has a neighbourhood $V$ 
such that, whenever $L$ is a primitive invariant 
lamination with minor gap intersecting $V$, then 
$\rho _{L}$ does not have a nondegenerate Levy 
cycle.  Note that if $\rho _{L}$ is type II, III 
or IV and has no nondegenerate Levy cycle, then 
$\rho _{L}$ is equivalent to a rational map. A 
point in $\Rm {Interior}(\Acc\SwtTilde  U)$ is 
\It{admissible } if it is either in an admissible 
leaf of $\Calig L$ or in a gap with at least one 
admissible side. A point in $\partial 
\Acc\SwtTilde  U$ is \It{admissible } if it has a 
neighbourhood whose intersection with $\Rm 
{Interior}(\Acc\SwtTilde  U)$ consists of 
admissible points. The set of admissible points 
is open by definition, and is denoted by 
$\Acc\SwtTilde  U_{ad}$.

\Subheading {5.9}  
We are now ready to prove the Admissible Boundary 
Theorem, which we restate. See also the 
introduction of \Cite{R2}. Again, a very slight 
correction has been made. Recall that if $m$ is 
the period of $r$ under $x\mapsto 2x\Rm {mod }1$, 
then $m$ is also the period of $0$ under 
$\Acc\ChOline s=\Acc\ChOline s_{r}$, and under 
lamination maps of all invariant laminations with 
minor leaves and minor gap sides in $\Calig 
L_{r}$.

\Theorem {Admissible Boundary Theorem}  Let 
$\Acc\SwtTilde  \nu $ be a leaf in $\partial 
\Acc\SwtTilde  U_{ad}$. Then either 
$\Acc\SwtTilde  \nu $ is a side of the minor gap 
$G(L)$ of a primitive invariant lamination $L$, 
or there is a sequence $\lbrace L_{n}\rbrace $ of 
primitive invariant laminations with minor gaps  
$G(L_{n})$ such that $\Acc\SwtTilde  \nu 
=\lim_{{n\rarrow \infty }} G(L_{n})$. Also, the 
following hold. Let $\rho =\rho _{L}$ or $\rho 
_{L_{n}}$ be the lamination map of $L$ or 
$L_{n}$. Then $\rho $ is type IV, not equivalent 
to a rational map, and $\infty $ is periodic of 
period $\leq m$ under $\rho $. There is a 
piecewise-geodesic fixed holed sphere $A_{0}\cup 
B_{0}$ for $L$ (or $L_{n}$) such that $A_{0}$ 
intersects sides of the minor gap of $L$ (or 
$L_{n}$).\endTheorem   
\Proof {Proof}  There are minor gaps of primitive 
invariant laminations intersecting any 
neighbourhood of $\Acc\SwtTilde  \nu $ for which 
some of the corresponding lamination maps have 
nondegenerate Levy cycles, and some do not. So at 
least $L$ or $\lbrace L_{n}\rbrace $ exist with 
minor gaps $G(L)$, $G(L_{n})$ such that either 
$\Acc\SwtTilde  \nu $ is a side of $G(L)$ or 
$\Acc\SwtTilde  \nu =\lim_{{n\rarrow \infty }} 
G(L_{n})$, and $\rho $ has a nondegenerate Levy 
cycle, where $\rho =\rho _{L}$ or $\rho 
_{L_{n}}$. In what follows, the point to remember 
is that two points on disjoint geodesics are 
close if the two geodesics have the same 
$S^{1}$-crossings for a considerable distance in 
both directions. 

Write $\Acc\SwtTilde  \mu $ for the minor leaf of 
$L_{n}$, or $\Acc\SwtTilde  \mu =\Acc\SwtTilde  
\nu $ if $\Acc\SwtTilde  \nu $ is a side of 
$G(L)$, so that $\Acc\SwtTilde  \nu $ is 
approximated by leaves $\Acc\SwtTilde  \mu $. 
 Write $L$ rather than $L_{n}$ for the moment. The 
Nonrational Lamination Map Theorem implies that 
there is a piecewise-geodesic periodic 
holed-sphere $A_{0}\cup B_{0}$ for $L$. Let 
$A_{k}\cup B_{k}$ denote the subset of $\rho 
^{-k}(A_{0}\cup B_{0})$ which is homotopic to 
$A_{0}\cup B_{0}$ in $\Acc\ChOline \Bbd 
C\setminus X(\rho )$ and such that each component 
of $A_{k}$ contains a component of $A_{0}$. (See 
4.2.) The First Invariant-implies-Parameter 
Theorem can be interpreted as follows. Here, 
$\varepsilon $, $\delta $ depend only on a 
neighbourhood $V$ of a compact subset of 
$\Acc\SwtTilde  \nu $ which $G(L)$ intersects. 
Write $L(A_{k},\beta ))$ for the lamination $L$ 
constructed from $(A_{k},\beta )$ as in the First 
Theorem. 

\It{Given $\varepsilon >0$, there is $\delta >0$ 
such that, either $\Acc\SwtTilde  \mu \cap V$ 
lies in an $\varepsilon $-neighbourhood of 
$G(L(A_{k},\beta ))$, or no $A_{k}$ passes within 
$\delta $ of $V\cap \Acc\SwtTilde  \mu $. } 

In fact, given $\varepsilon >0$, the first 
possibility must hold for some $\Acc\SwtTilde  
\mu $. Otherwise, the Second 
Invariant-implies-Parameter Theorem gives the 
following.

\It{There is $\eta >0$, which depends on $\delta 
$, and hence on $\varepsilon $, such that every 
primitive lamination $L'$ with minor gap passing 
within $\eta $ of $\Acc\SwtTilde  \mu \cap V$ 
satisfies one of the following.}

\par \noindent \It{ a) $A_{k}$ (for any $k$) is 
contained in leaves or finite sided gaps of $L'$. }
\par \noindent  \It{ b) $L'\subset L''$ with 
$A_{k}\subset L''$ such that $L''$ and $L'$ have 
the same minor leaf.} 
\par \noindent  \It{Xc) $L'$ has a nonsimply-connected 
minor gap. }
 
By Lemma 5.6 and 5.7, we can choose $L''$ so that 
$A_{k}\cup B_{k}$ is a piecewise-geodesic 
periodic holed sphere for $L''$, and then modify 
to a primitive lamination. 
(The hypothesis of 5.6 is satisfied, because
$A_{0}\cup B_{0}$ is a piecewise-geodesic periodic
holed-sphere for $L$.) If $\Acc\SwtTilde  \mu 
$ is sufficiently close to $\Acc\SwtTilde  \nu $, 
this implies $\Acc\SwtTilde  \nu $ is in the 
interior of $\Acc\SwtTilde  U\setminus 
\Acc\SwtTilde  U_{ad}$, giving a contradiction.

Now using laminations $L(A_{k},\beta )$ instead 
of the original $L_{n}$ (or $L$), the 
corresponding lamination maps are all type IV, 
not equivalent to rational maps. If the other 
property (concerning the piecewise-geodesic 
\It{fixed } holed-sphere) is not satisfied, then 
the period of $\infty $, for all the maps in the
sequence, is less then $m$, by the 
Nonrational Lamination Map Theorem. We can then 
repeat the entire process. After finitely many 
repeats, we obtain a sequence with the required 
properties.

\Subheading {5.10. Preliminaries to the 
Admissibility Proposition} 

 Every leaf intersecting the admissible set is 
contained in it. If a lamination map $\rho _{L}$ 
of a lamination $L$ is critically finite and 
equivalent to a rational map, one would expect 
all sides of the minor gap to be in the 
admissible set. Actually, if $\rho _{L}$ is 
equivalent to a polynomial this cannot quite be 
true - because then some tunings of $\rho _{L}$ 
(with minor gaps contained in the minor gap of 
$L$) will not be equivalent to rational maps. But 
if $\rho _{L}$ is equivalent to a polynomial and 
hence to $\Acc\ChOline s_{r}$ for some $r$, then 
$\rho _{L}$ is associated to the gap of $\Calig 
L_{r}$ containing $0$, which is admissible. 
Perhaps the easiest degree two rational maps to 
study are matings and captures, originally 
defined by Douady and Hubbard \Cite{D} and 
Wittner \Cite{W}. They are also described in 1.19 
of \Cite{R2}, where a subset $\Acc\SwtTilde  
U_{mc}$ of $\Acc\SwtTilde  U$ is described in 
which all matings and captures are represented, 
and a description is given (without proof) of 
$\Acc\SwtTilde  U_{ad}\cap \Acc\SwtTilde  
U_{mc}$. This set is connected and contains up to 
equivalence all matings and captures which are 
equivalent to rational maps. For matings, the 
description is the result of Tan Lei's thesis
\Cite {TL}. In 
both cases, it follows from the more general 
Admissible Boundary Theorem, or even the 
Nonrational Lamination Map Theorem. A simplifying 
factor is that, if $\Acc\SwtTilde  \mu $ and 
$\Acc\SwtTilde  \mu '$ are two minor leaves in 
$\Acc\SwtTilde  U_{mc}$ of laminations $L$ and 
$L'$, and $\Acc\SwtTilde  \mu $ separates 
$\Acc\SwtTilde  \mu '$ from $0$, then all 
periodic leaves in $L$ are  in $L'$. So now we 
state the Admissibility Proposition, in a 
slightly different, but equivalent form to that 
given in 1.17 of \Cite{R2}, given that the 
definition of admissibility has been altered in 
5.8 from that of 1.17 of \Cite{R2}.
\Theorem {The Admissibility Proposition}  Let $f$ 
be a critically finite degree two rational map 
which is not equivalent to a polynomial, mating 
or capture. Then  whenever $L$ is an invariant 
lamination with minor leaf or gap in 
$\Acc\SwtTilde  U$ and the lamination map $\rho 
_{L}$ is equivalent to $f$, the minor leaf or gap 
lies in the admissibility set $\Acc\SwtTilde  
U_{ad}$.\endTheorem   

\Subheading {5.11. The Inadmissible Type II Case} 

\Theorem {Lemma}  

Let $L$ be an invariant lamination with type II 
lamination map $\rho _{L}$ which is equivalent to 
a rational map $f.$  Let $L'$ be the lamination 
$L\cup \Acc\ChOline L_{r}$ with type IV 
lamination map $\rho '$ which is not equivalent 
to a rational map. Then $\rho '$ is equivalent to 
a mating and $\rho _{L}$ is equivalent to a 
capture. \endTheorem   
\Proof {Proof}  Let $G_{\infty }$ be the gap of 
$L$ containing $\infty $, and let $\Acc\SwtTilde  
G$ be a lift to the universal cover of 
$\Acc\ChOline \Bbd C\setminus K$, with projection 
$\pi :\Acc\ChOline {\Acc\SwtTilde  G}\rarrow 
\Acc\ChOline G_{\infty }$. Let $\Acc\SwtTilde  
\rho $ be a lift of $\rho _{L}^{m}\Midvert 
G_{\infty }$ fixing $\Acc\SwtTilde  G$, and 
extending to $\Acc\ChOline {\Acc\SwtTilde  G}$,  
if $m$ is the period of $G_{\infty }$ under $\rho 
_{L}$. Let $\Phi :\Acc\ChOline \Bbd C\rarrow 
\Acc\ChOline \Bbd C$ be the semiconjugacy 
satisfying $\Phi \Circ \rho _{L}=f\Circ \Phi $.  
Let $\varphi _{1}: \lbrace z:\Midvert z\Midvert 
\leq 1\rbrace \rarrow \Acc\ChOline 
{U_{2}(f)}$ and $\varphi _{2}:\Acc\ChOline 
{\Acc\SwtTilde  G}\rarrow \lbrace z:\Midvert 
z\Midvert \leq 1\rbrace $ satisfy
$$\varphi _{1}\Circ \varphi _{2}=\Phi \Circ \pi 
,$$
$$\varphi _{1}(z^{4})=f\Circ \varphi _{1}(z),$$
$$\varphi _{2}\Circ \Acc\SwtTilde  \rho 
(z)=(\varphi _{2}(z))^{4}$$
for all $z$. Let $\omega _{1}$, $\omega _{2}$ be 
the cube roots of unity which are the images 
under $\varphi _{2}$ of the endpoints of the lift 
to $\Acc\SwtTilde  G$ of $\ell $, where $\ell $ 
is the unique leaf of $\Acc\ChOline 
L_{r}\setminus L$ of oriented period $m$ in 
$G_{\infty }$. Let $\alpha $ be an arc in 
$\Acc\ChOline {U_{2}(f)}$ which is the image 
under $\varphi _{1}$ of an arc in $\lbrace 
z:\Midvert z\Midvert \leq 1\rbrace $ joining 
$\omega _{1}$ and $\omega _{2}$. We can choose 
$\alpha $, and can perturb $f$ to a critically 
finite branched covering $f'$, so that $f'=f$ 
outside a neighbourhood of the critical orbit of 
$f$, $f'$ preserves the set $\cup 
_{i=0}^{m-1}f'^{i}\alpha $, $f'^{i}\alpha \subset 
f'^{i}U_{2}(f)=f^{i}U_{2}(f)$ separates the two 
points of $X(f')$ in $f^{i}U_{2}(f)$, and 
$f'\simeq \rho '$. So $f'$ admits a minimal Levy 
cycle $\Gamma $, that is, one with only one 
complementary component which is not a disc. (See 
1.6 of \Cite{R2}.) Let $P$ be the closure of the 
nondisc component of $\Acc\ChOline \Bbd 
C\setminus (\cup \Gamma )$. Let $F_{0}:P\times 
[0,1]\rarrow \Acc\ChOline \Bbd C\setminus X(f')$ 
be a homotopy between $P$ and a component of 
$f'^{-1}(P)$ with $f'\Circ 
F_{0}(x,1)=F_{0}(x,0)$. (This exists for a 
minimal Levy cycle.)  We can assume 
$F_{0}(P\times [0,1])\cap (\cup _{i\geq 
0}f^{i}(U_{2}(f))\subset \cup _{i\geq 
0}f'^{i}\alpha $. Defining $F_{n}$ inductively by 
$f'\Circ F_{n+1}=F_{n}$, 
$F_{n+1}(x,0)=F_{n}(x,1)$, and taking limits, we 
can assume $\cup \Gamma \subset \cup _{i\geq 
0}\Acc\ChOline {f'^{i}\alpha }$, if we drop the 
assumption that the loops of $\Gamma $ are 
embedded and disjoint. (See Section 8 of \Cite{R1}
for a similar construction.) This set must contain at 
least one point fixed by $f'$ by the Lefschetz 
Theorem. So one endpoint of $\alpha $ must be 
fixed by $f$ and by $f'$. The endpoints of 
$\alpha $ must be distinct if $m\geq 3$, for 
topological reasons, but also if $m=2$, since 
then $f$ is $z\mapsto z^{-2}$, up to M\"obius 
conjugation. Then the second endpoint must also 
be fixed. Then we can find a circle $\Gamma '$ 
invariant up to homotopy under $f'$, and 
separating the critical orbits, which is 
equivalent to $f'$ being a mating. See Diagram 5.

\Diagram {{250}Diagram 5}

Then 
$$f\simeq \sigma _{\Acc\ChOline \eta }\Circ 
\sigma _{\gamma }\Circ f',$$
where $\gamma $, $\eta $ are loops crossing 
$\Gamma '$ once each, (where $\Acc\ChOline \eta 
$, as usual, is $\eta $ reversed), and $f'(\eta 
)$ is homotopic to $\gamma $ in $\Acc\ChOline 
\Bbd C\setminus X(f')$. Let $f_{0}$ be the 
polynomial obtained by modifying $f'$ on the 
component of $\Acc\ChOline \Bbd C\setminus \Gamma 
'$ containing $c_{2}(f')$, so that $c_{2}(f_{0})$ 
is fixed. Up to homotopy in $\Acc\ChOline \Bbd 
C\setminus X(f_{0})$, $f_{0}^{-1}$ preserves 
$\Gamma '$ and the components of $D\setminus 
(\cup (f')^{i}\alpha )$, where $D$ is the 
component of $\Acc\ChOline \Bbd C\setminus \Gamma 
$ containing $c_{1}$. Hence, $f_{0}$ is a 
polynomial (corresponding to a minimal minor leaf 
in $QML$: see 1.10 of \Cite{R2}). Then, modifying 
$\gamma $ and $\eta $ also, there are paths 
$\beta $ and $\zeta $ in $\Acc\ChOline \Bbd 
C\setminus X(f_{0})$ with $f_{0}(\zeta )$ 
homotopic in $\Acc\ChOline \Bbd C\setminus 
X(f_{0})$ and 
$$f\simeq \sigma _{\Acc\ChOline \zeta }\Circ 
\sigma _{\beta }\Circ f_{0},$$
as required.

\Subheading {5.12. Definition of $\Bf{W }_{\Bf{1 
}}$, $\Bf{W }_{\Bf{2 }}$}  

Let $f$ be a critically finite rational map with 
numbered critical points $c_{1}$, $c_{2}$. Let 
$c_{1}$ be periodic, and suppose $f$ is not type 
II. Then we define sets $W_{1}$, $W_{2}$ as 
follows. If $f$ is not type III or IV (so that, 
in particular, $c_{2}$ is not periodic) then 
$W_{1}=W_{2}$ is simply $\lbrace 
f^{i}(c_{2}):i>0\rbrace $. If $f$ is type III, 
then $W_{1}=W_{1}(x)=W_{2}$ is $\lbrace 
f^{i}(x):i\geq 0\rbrace $ for an eventually 
periodic point $x\in \partial f(U_{2})$. If $f$ 
is type IV, then $W_{2}$ is the closure of the 
(finite) union of components of $\Acc\ChOline 
\Bbd C\setminus J(f)$ which intersect the forward 
orbit of $c_{2}$, and $W_{1}$ is the smallest set 
with the following properties.

\par \noindent 1. $W_{2}\subset W_{1}$ 
\par \noindent 2. $W_{1}$ contains $\Acc\ChOline 
U$, whenever $U$ is a component of $\Acc\ChOline 
\Bbd C\setminus J(f)$ in the full orbit of 
$U_{2}(f)$ with $\Acc\ChOline U\cap W_{1}\not 
=\phi  $.

Then $W_{1}$ has finitely many components. By 5.4 
of \Cite{R2}, if $U$ and $V$ are in the same 
period $p$ orbit of components of $\Acc\ChOline 
\Bbd C\setminus J(f)$, with $U\not =V$, then 
$\Acc\ChOline U\cap \Acc\ChOline V$ consists of 
at most one point, which is fixed by $f^{p}$.

\Theorem {5.13. Lemma}  

Let $f$ be type IV. All components of 
$\Acc\ChOline W_{1}$ are simply connected, unless 
$f$ is equivalent to a polynomial or mating 
\endTheorem   

\Proof {Proof} Suppose that components of 
$\Acc\ChOline W_{1}$ are not simply connected and 
that $f$ is not a polynomial, so that $c_{2}(f)$ 
has period $p>1$. By 5.4 of \Cite{R2}, 
$\Acc\ChOline U$ is simply connected for any 
component $U$ of $\Acc\ChOline \Bbd C\setminus 
J(f)$. So the components of $W_{1}$ are larger 
than the closures of components of $\Acc\ChOline 
\Bbd C\setminus J(f)$. There is $x\in \cap 
_{j\geq 0}f^{jt}\Acc\ChOline {U_{2}(f)}$ with $x$ 
of period $t<p$. Let $\gamma $ be an embedded 
circle enclosing $x$ and lying in the union of 
the $f^{jt}U_{2}(f)$ and a neighbourhood of $x$, 
and enclosing $x$ and $f^{jt}(c_{2})$ ($j\geq 
0$). Then $\gamma $ is isotopic to a component 
$\gamma _{1}$ of $f^{-t}(\gamma )$ in 
$\Acc\ChOline \Bbd C\setminus X(f)$, with 
$f^{t}\Midvert \gamma _{1}$ degree 2. Then $f$ 
must be a tuning $g\vdash _{c_{2}}h$ of a type IV 
rational map $g$ (see \Cite{R2} 1.20), where $h$ 
is a polynomial for which the boundaries of the 
attractive basins have a common fixed point. It 
is possible to define a branched covering $f_{1}$ 
which is equivalent to $f$ and equal to $g$ on 
$J(g)$. Moreover, the Semiconjugacy Proposition 
of 4.1 of \Cite{R2} implies that there is 
$\varphi :\Acc\ChOline \Bbd C\rarrow \Acc\ChOline 
\Bbd C$ with $\varphi \Circ f_{1}=f\Circ \varphi 
$, and $\varphi $ is a limit of homeomorphisms.
Since $f_{1}\Midvert J(g)$ is expanding, each set
$\varphi ^{-1}x$, which is connected and simply-
connected, has finite intersection with $J(g)$.
 (See 4.6 of \Cite{R2}.) 
Then $\Acc\ChOline W_{1}=\varphi (\cup _{i\geq 
0}g^{i}\Acc\ChOline {U_{2}(g)})$. So components 
of $\cup _{i\geq 0}g^{i}\Acc\ChOline {U_{2}(g)}$ 
are not simply connected. Then by 5.4 
of \Cite{R2}, $g$ is a polynomial. Then $f$ is a 
mating, since the tuning of a polynomial is a 
mating.

\Theorem {5.14. Lemma}    

There is an arbitrarily small neighbourhood $W$ 
of $W_{1}$ and such that, for all $n,$ each 
component of $f^{-n}W$ is either contained in, or 
disjoint from, $W.$ \endTheorem   

\Proof {Proof}  We start with a neighbourhood 
$W'$ of $W_{1}$ such that $W'$ contains the 
components of $f^{-p}W'$ intersecting $W_{1}$ (if 
$p$ is the least integer $>0$ with 
$f^{p}W_{1}=W_{1}$). The idea is then to reduce 
$W'$ to obtain $W$, without violating the 
neighbourhood condition. We see that there is 
$\delta >0$ with the following property.

\par \noindent a) If $W''$ is a component of 
$f^{-n}W'$ for some $n$, $W''$ does not contain 
$W_{1}$ and has diameter $\varepsilon $, then 
$\Rm {dist}(W_{1},W'')>\delta \varepsilon $

Of course, this property remains true if we make 
$W'$ even smaller. Given any $\eta >0$, we can 
choose $W'$ so close to $W_{1}$ that the 
following is true.

\par \noindent b) If $W''$ is as in a), and 
$W''\cap W'\not= \phi$, then 
$\varepsilon <\eta $.
 
There is $M>0$ such that the following is true 
(and remains true if $W'$ is made smaller).

\par \noindent c) For all $z$, $w$ in the domain 
of a local inverse $S$ of $f^{n}$ defined on a 
subset of $W'\setminus W_{1}$,
$${1\over M}\leq {\Midvert S'(z)\Midvert \over 
\Midvert S'(w)\Midvert }\leq M.$$
It follows that if $W'$ is close enough to 
$W_{1}$, then given $\nu >0$,

\par \noindent d) if $W'$, $W''$ are as in a) and 
$W'\cap W''\not =\phi  $ and $W_{n}'$, $W_{n}''$ 
are preimages under $f^{n}$ with $W_{n}'\cap 
W_{n}''\not =\phi  $, then $W_{n}''$ is contained 
in the $\Rm {diameter}(W_{n}')\times \nu 
$-neighbourhood of $W_{n}'$ and has diameter 
$<\Rm {diameter}(W_{n}')\times \nu $.

 Now we inductively define $T_{0}=W'$, and 
$T_{n+1}$ is $T_{n}\setminus R_{n}$, where 
$R_{n}$ is the union of those components of 
$f^{-k}T_{n}$ (for all $k$) which intersect, but 
are not contained in, $T_{n}$. We take $W=\cap 
_{n}T_{n}$. Then $W$ has the required properties 
- provided $W'$ is close enough to $W_{1}$ - 
because every component of $R_{n+1}$ intersects a 
component of $R_{n}$. Hence we see inductively 
that every component of $R_{n}$ is distance $\leq 
$
$$\sum _{i=1}^{n}\nu ^{i}\times \Rm 
{diameter}(C)$$
from a component $C$ of $R_{0}$ which is itself 
distance $>\delta \times \Rm {diameter}(C)$ from 
$W_{2}$, by a). So we only need to take 
$\displaystyle{{\nu \over 1-\nu }<\delta }$. By 
its definition $T_{n}$ contains the component of 
$f^{-p}(T_{n})$ containing $W_{1}$, so $W$ 
contains the component of $f^{-p}W$ containing
$W_{1}$.

\Subheading {5.15. Proof of the Admissibility 
Proposition}  

Let $f$ be a critically finite rational map which 
is not equivalent to a polynomial, mating or 
capture. Let $L$ be an invariant lamination with 
minor leaf or gap in $\Acc\SwtTilde  U$, such 
that the lamination map $\rho $ is equivalent to 
$f$. By 5.11, if $f$ is type II, the lamination 
map for the lamination $L\cup \Acc\ChOline L_{r}$ 
is equivalent to a rational map of type IV. So 
now we can assume that $f$ is not type II. Then 
the Semiconjugacy Proposition 4.1 of \Cite{R2} 
gives a semiconjugacy $\Phi :\Acc\ChOline \Bbd 
C\rarrow \Acc\ChOline \Bbd C$ such that $\Phi 
\Circ \rho =f\Circ \Phi $, $\Phi $ is a limit of 
homeomorphisms and collapses each leaf and finite 
sided gap of $L$ to a point. If $f$ is type III, 
let $x$ be the image under $\Phi $ of the minor 
leaf of $L$. Then let $W_{1}$ be as in 5.12, with 
$W_{1}=W_{1}(x)$ if $f$ is type III. By 5.13, 
$\Acc\ChOline W_{1}$ is simply-connected. Let $W$ 
be as in 5.14, and we can choose $W$ close enough 
to $W_{1}$ so that $\Acc\ChOline W$ is 
contractible in $\Acc\ChOline \Bbd C\setminus 
\lbrace f^{i}(c_{1}):i\geq 0\rbrace $. Put 
$$V=\pi ^{-1}((\Acc\ChOline \Bbd C\setminus 
(K\cup (\cup \Acc\ChOline L_{r})))\cap \Phi 
^{-1}W)\subset \Acc\SwtTilde  U,$$
 where $\pi $ is the universal covering map. Then 
$\pi (\Acc\ChOline V)$ is contractible in 
$\Acc\ChOline \Bbd C\setminus \lbrace 
\Acc\ChOline s^{i}(0):i\geq 0\rbrace $ and is a 
neighbourhood of the minor leaf of $L$, and of 
the minor gap if this exists, such that $\pi (V)$ 
contains any component of $\rho ^{-n}(V)$ it 
intersects, for any $n\geq 0$. We claim that the 
proof of the Admissibility Proposition is 
completed by the following lemma, using the 
notation established here.

\Theorem {Lemma}  Let an integer $m>0$ be given. 
Let $A$ denote the set of leaves of $L$ of period 
$\leq m$. Let $V'\subset \vee $ be a 
neighbourhood of the minor leaf or minor gap. Let 
$L'$ be any primitive invariant lamination with 
minor gap intersecting $V'$. Let $A'$ be the set 
of leaves of $L'$ of period $\leq m$. Then if 
$V'$ is sufficiently small, every leaf of $A'$ is 
in $A$, or a finite-sided gap of $L$ bounded by 
leaves of $A$, or in $\pi (V)$.

\endTheorem   

We see that the lemma implies the proposition as 
follows. Let $m$ be the period of $r$ under 
$x\mapsto 2x\Rm { mod }1$, so that $m$ is the 
period of $0$ under $\Acc\ChOline 
s_{r}=\Acc\ChOline s$, $\rho $, $\rho '$, and, 
indeed, under $\tau $, where $\tau $ is any 
lamination map of an invariant lamination with 
minor leaf in $\Calig L_{r}$.  Suppose the 
lamination map $\rho '$ of $L'$ is a critically 
finite branched covering. By the Nonrational 
Lamination Map Theorem, $\rho '$ does not have a 
nondegenerate Levy cycle provided $\Acc\ChOline 
{A'}$ is contractible in $\Acc\ChOline \Bbd 
C\setminus X(\Acc\ChOline s)$ (to finitely many 
points). But this is true, because $\Acc\ChOline 
{A\cup \pi (V)}$ is contractible in $\Acc\ChOline 
\Bbd C\setminus X(\Acc\ChOline s)$, since $A$ 
maps to points under $\Phi $ and $\Acc\ChOline 
{\pi (V)}$ maps to $\Acc\ChOline W$.

\Subheading {5.16. Proof of Lemma 5.15}  

Let  $C$ be an orbit in $A'$. We have $C=\cup 
_{n\geq 0}C_{n}$, where $C_{n}\subset 
C_{n+1}\subset \rho '^{-1}C_{n}$ and the number 
of segments in $C_{0}$ can be bounded 
independently of $L'$. (The bound depends on 
$V'$.) By taking $V'$ sufficiently small, given 
$t$, we can assume that any geodesic of $L$ of 
$\geq t$ segments intersecting $\pi (V')$ cannot 
be isotoped into $\partial \pi (V)$ by good 
isotopy. We can also take $V'$ small enough that  
any such $C_{0}$ does not lift to lie strictly 
between the minor leaves of gaps of $L$, $L'$, up 
to good isotopy. Take any $k$ such that $C_{i}$ 
does not intersect $\pi (V')$ for $i<k$. By the 
Second Invariant-implies-Parameter Theorem (with 
$L$ replacing $L'$), $C_{k}$  does not intersect 
$L$ transversally. So either $C_{k}\subset \pi 
(V)$, or $C_{k}\cap \pi (V)=\phi  $, and then 
$C_{k}$ lies in leaves or finite-sided gaps of 
$L$ up to good isotopy. But $\rho ^{-1}=\rho 
'^{-1}$ outside $\pi (V)$, and $\rho ^{-k}(\pi 
(V))=\rho '^{-k}(\pi (V))$ for all $k\geq 0$. So 
either $C\subset \cup _{n\geq 0}\rho '^{-n}C_{0}$ 
and $C\subset \pi (V)$ or, for all $k$, $C_{k}$ 
lies in leaves or finite-sided gaps of $L$ up to 
good isotopy, and hence so does $C$.
\end

\pageno=57

\References {
References
}

\Benchmark 
\Cite{D} Douady, A.: Syst\`emes dynamiques 
holomorphes. S\'eminaire Bourbaki 1983, 
\It{Ast\'erisque 105/6 } (1983), 39-63.

\Benchmark 
\Cite{D-H1} Douady, A. and Hubbard, J.H.: Etudes 
dynamiques des polyn\^omes complexes, avec la 
collaboration de P. Lavaurs, Tan Lei, P. 
Sentenac. Parts I and II, Publications 
Math\'ematiques d'Orsay, 1985.

\Benchmark 
\Cite{D-H2} Douady, A. and Hubbard, J.H.: 
It\'erations des polyn\^omes quadratiques 
complexes. \It{C.R.Acad.Sci. Paris, S\'erie I, 
t.294 } (1982), 123-126.

\Benchmark 
\Cite{D-H3} Douady, A. and Hubbard, J.H.: A proof 
of Thurston's topological characterization of 
rational functions. Mittag-Leffler preprint, 
1985.

\Benchmark 
\Cite{F-L-P} Fahti, A., Laudenbach, F., and 
Po\'enaru, V.: Travaux de Thurston sur les 
Surfaces. \It{Ast\'erisque 66-67 } (1989).

\Benchmark 
\Cite{Fr} Franks, J.: Anosov Diffeomorphisms. 
\It{Proc. Symposia in Pure Math., 14 } (1968), 
61-93.

\Benchmark 
\Cite{R1} Rees, M.: Components of degree two 
hyperbolic rational maps. \It{Invent. Math., } 
100 (1990), 357-382.

\Benchmark 
\Cite{R2} Rees, M.: A Partial Description of the 
Parameter Space of Rational Maps of Degree Two. 
Part 1. To appear in \It{Acta Math. }

\Benchmark 
\Cite{S} Shishikura, M.: In preparation.

\Benchmark 
\Cite{T} Thurston, W.P.: On the Combinatorics of 
Iterated Rational Maps. Preprint, Princteton 
University and I.A.S., 1985.

\Benchmark 
\Cite{TL} Tan Lei: Accouplements des polyn\^omes 
complexes. Th\`ese, Universit\'e de Paris-Sud, 
Orsay, 1987.

\Benchmark 
\Cite{W} Wittner, B.: On the bifurcation loci of 
rational maps of degree two. Thesis, Cornell 
University, 1988.\endReferences 
\end